\documentclass[preprint,12pt]{elsarticle}
\usepackage{amsmath, amssymb, graphicx, natbib}
\usepackage[margin=2.5cm]{geometry}
\usepackage{xurl}
\usepackage{hyperref}
\usepackage{xcolor}
\usepackage{float}
\usepackage{amsmath}
\usepackage{amssymb}
\usepackage{amsthm}
\usepackage{amsfonts}
\usepackage{algorithm}
\usepackage{algpseudocode}

\usepackage{booktabs}
\usepackage{url}
\usepackage{multirow}
\usepackage{afterpage}
\usepackage{placeins}
\usepackage{caption}
\usepackage{subcaption}

\usepackage{bm}
\usepackage{mathtools}
\usepackage{geometry}
\usepackage{microtype}

\journal{Journal of Computational Physics}

\begin{document}

\begin{frontmatter}

\title{Structure-preserving Fourier Neural Operators for Long-time Cahn--Hilliard Dynamics under Coarse Temporal Supervision}

\fntext[note1]{Corresponding author. Email: yzhou@colostate.edu}
\author{Yingli Li}
\address{Department of Mathematics, Colorado State University, Fort Collins, CO 80523, USA}
\author{Y.C. Zhou \fnref{note1}}
\address{Department of Mathematics, Colorado State University, Fort Collins, CO 80523, USA}

\date{\today}

\begin{abstract}
    Accurate long-time prediction of Cahn--Hilliard dynamics, particularly in the late-stage coarsening regime, is computationally demanding because it requires simulations over extended physical time. Various machine learning approaches, such as physics-informed neural networks and operator learning, have been developed to speed up the simulations through efficient surrogate evaluations. While the learned map can be accurate over a single learned interval, small prediction error will accumulate when the map is repeatedly applied in a long-time prediction, leading to non-dissipative free-energy evolution, large excursions beyond a nominal phase-field range, and degraded coarsening fidelity. In this work, we identify the structure function for the coarsening dynamics of the Cahn--Hilliard equation and propose a two-stage structure-preserving Fourier Neural Operator (FNO) learning framework. Stage~I introduces a bound-conforming FNO (bcFNO) with a soft, magnitude-dependent amplitude penalty. Stage~II retains the bound-conforming objective and supplements it with structure-function regularization, yielding a structure-preserving FNO (spFNO) that corrects late-stage coarsening statistics. Numerical experiments on the various spatial resolutions and reference-solver time steps suggest that the bcFNO stabilizes the long-time prediction by suppressing large amplitude excursions, without changing the low online cost of the FNO. Furthermore, the structure-function term consistently reduces discrepancies in the normalized structure function and provides further correction in associated late-stage coarsening quantities, including the $L^3(t)$ growth trend and the scaling collapse. These results indicate that our multi-stage implementation of physical principles in operator learning could be applied to other multiscale systems with statistical scaling, such as the functionalized Cahn--Hilliard equations and turbulence, provided that the state constraints and statistical observables are adapted to the governing dynamics.
\end{abstract}
\end{frontmatter}

\section{Introduction}\label{sec:intro}
The Cahn--Hilliard equation is a canonical phase-field model for phase separation and interfacial pattern formation in binary systems \cite{Cahn1958,Cahn1961}.
Originally introduced to describe spinodal decomposition \cite{Cahn1958,Grant1993},
it has since been widely used in materials science \cite{Zhu2001,Fife2003}, polymer physics \cite{Smolders1971,Anderson2002}, multiphase-flow modeling \cite{Badalassi2003},
and other interfacial-dynamics problems involving topological changes such as merging, pinch-off, and domain annihilation \cite{Tierra2015}. Its widespread application has motivated the development of accurate and efficient numerical methods for solving the equation, including finite-difference \cite{LezamaAlvarez2013}, finite-element \cite{Banas2008,Cherfils2014}, and spectral discretizations \cite{Lecoq2011,Christlieb2013}. Common temporal discretization strategies include semi-implicit and fully implicit schemes \cite{bosch2014fast, chen1998applications}, convex-splitting methods \cite{eyre1998unconditionally}, operator-splitting methods \cite{elliott1989second, cheng2015fast}, and exponential time-differencing schemes \cite{cox2002exponential, kassam2005fourth}.
More recent energy-stable formulations based on scalar auxiliary variable (SAV) and invariant energy quadratization (IEQ) approaches enable linearized time integration while preserving a discrete modified energy-dissipation law \cite{shen2018scalar,yang2020convergence}.

Despite these advances, high-fidelity simulations of the Cahn--Hilliard equation remain computationally demanding. The fourth-order diffusion operator makes the dynamics stiff, whereas accurate resolution of diffuse interfaces requires fine spatial grids. For standard explicit discretizations of the fourth-order term, stability typically imposes \(\delta t\lesssim C h^4\). Mesh refinement would greatly increase the number of steps required to reach a fixed physical time. Although semi-implicit, implicit, and energy-stable schemes alleviate this stability restriction, their time steps must still be sufficiently small to resolve rapid interfacial relaxation and nonlinear transient dynamics, and implicit treatments require the solution of linear or nonlinear systems at every step. The coexistence of rapid early-stage phase separation and slow late-stage coarsening consequently entails a large number of time steps over long-time simulations, making the generation of accurate high-fidelity reference solutions computationally expensive, particularly in high dimensions. This cost will grow further when many long-time trajectories are required across multiple initial conditions or parameters. These challenges motivate scientific machine learning (SciML) surrogate models that complement conventional numerical solvers by learning approximations of the evolution map from high-fidelity data.  The model trained offline can be used to generate long-time rollouts at a substantially lower online cost.

SciML studies that take the Cahn--Hilliard equation or related phase-field systems as central test problems have developed along several directions. A major line of work concerns physics-informed neural networks (PINNs) for the stiff, nonlinear fourth-order equation. Adaptive spatio-temporal sampling was introduced to improve the accuracy and efficiency of PINN solvers for Allen--Cahn and Cahn--Hilliard equations \cite{wight2020solving}, while sequential training over successive time segments was used to address the loss of accuracy of PINNs for nonlinear, higher-order time-dependent equations \cite{mattey2022novel}. Subsequent studies considered more robust training formulations \cite{zhang2023robust}, mass-preserving adaptive PINNs for Cahn--Hilliard equations with different bulk potentials, singularities, and three-dimensional settings \cite{huang2024mass}, causal training for one-dimensional Cahn--Hilliard dynamics \cite{hu2025improved}, and residual-based adaptive weighting and a priori error analysis \cite{zhang2026priori}.
These studies primarily aim to improve the accuracy and trainability of continuous space-time PINN approximations, with mass conservation being the physical property explicitly enforced in several Cahn--Hilliard formulations.

In another major line of SciML studies, operator learning is used to approximate PDE solutions or evolution operators from data.
Convolutional neural operators formulate these maps through continuous convolutional layers \cite{raonic2023convolutional}, whereas multiwavelet neural operators decompose operator kernels at several spatial scales \cite{gupta2021multiwavelet}. Multipole graph neural operator uses a multilevel graph construction to capture interactions over different spatial ranges \cite{li2020multipole}, and attention-based operator encodes the solution field in a reduced space before learning its evolution \cite{li2022transformer}.
For Cahn--Hilliard dynamics, Oommen et al. \cite{oommen2022learning} combined a convolutional autoencoder with DeepONet to learn the evolution of two-phase microstructures in a latent space \cite{kontolati2024learning}, while equivariant multiscale neural operators have been proposed to improve the representation of fine-scale structures and generalization in solution prediction \cite{xue2025equivariant}. Other data-driven work has focused on discovering Cahn--Hilliard-based phase-field dynamics from data \cite{kiyani2022machine}, predicting properties from Cahn--Hilliard-generated microstructures \cite{nguyen2024efficient}, or constructing surrogates for nonlocal Cahn--Hilliard models \cite{geng2025end}. Extensions to coupled phase-field systems include Allen--Cahn and Cahn--Hilliard \cite{chen2025pf, gangmei2025learning, chen2025sharp}, Cahn--Hilliard--Hele--Shaw \cite{fang2025simulating}, and Navier--Stokes--Cahn--Hilliard--heat transfer system \cite{lv2025simulation}, as well as fractional Cahn--Hilliard variants \cite{kang2025solving}.

Collectively, existing SciML studies of Cahn--Hilliard have emphasized PDE-residual accuracy, fieldwise prediction error, mass conservation, and computational acceleration. However, many studies employ simplified test problems, including deterministic initial conditions, prescribed trajectories, short prediction intervals, or controlled phase-field test problems. Such tests are useful for comparing field-wise approximation error, but they do not contain the essential features of the Cahn--Hilliard equation as a phase separation model, including late-stage coarsening, growth of the characteristic length scale and the scaling collapse of the structure function. 

Among these SciML approaches, neural operators are particularly attractive because they approximate mappings between function spaces and, once trained, can be evaluated efficiently across different initial conditions. Among various operator learning approaches, the Fourier neural operator (FNO) is particularly well suited to the periodic Cahn--Hilliard setting, since its Fourier layers efficiently represent global convolutions on periodic domains  \cite{Li2020}. 
However, accurate predictions on one-step supervised pairs in these SciML studies, FNO included, do not guarantee reliable long-time predictions. During recursive inference, small local prediction errors can accumulate and be amplified, eventually producing substantial deviations from the reference trajectory. For phase-field dynamics, such deviations may manifest as large amplitude excursions beyond the nominal phase-field range, non-dissipative free-energy behavior, or simply an energy blowout.

To mitigate these difficulties, we identify physical features of the Cahn--Hilliard dynamics in addition to the governing equation to guide the learning process. In the binary phase-field formulation, the minima of the bulk free-energy density are located at $\phi=\pm 1$. This inspires us to adopt $[-1,1]$ as a nominal phase-field reference range. The resulting bound-conforming FNO (bcFNO) training objective includes a soft, magnitude-dependent penalty that increases for excursions beyond this range without prohibiting small departures or imposing a hard projection. Thus, $[-1,1]$ is used as a physically motivated reference rather than being assumed to be a universal invariant of the Cahn--Hilliard equation.

Controlling nominal phase-field admissibility alone does not suffice to establish faithful recovery of coarsening morphology.
A closer look at the Cahn--Hilliard dynamics is needed to identify more physical features. 
This dynamics can be broadly divided into four regimes: early linear-instability growth, structure formation, late-stage coarsening, and finite-size saturation \cite{konig2021two}. Late-stage coarsening is of particular interest because it evolves over long physical-time scales and is governed by collective domain growth and rearrangement rather than short-lived local transients. As coarsening proceeds, the characteristic domain scale increases and the interfacial density decreases \cite{Chakrabarti1993}. To study this regime, we consider source-free Cahn--Hilliard trajectories initialized from independently sampled random initial perturbations. These initial fields generate distinct coarsening morphologies, and the learned evolution map is evaluated on  held-out long-time predictions through the coarsening regime.

More importantly, late-stage Cahn--Hilliard coarsening exhibits dynamic scaling, under which suitably rescaled structural statistics approach an approximately time-independent form \cite{Zhu1999,Konig2021}. In diffusion-controlled Cahn--Hilliard coarsening, the characteristic length is commonly expected to follow $L(t)\sim t^{1/3}$ at late times \cite{lifshitz1961kinetics,zhu1999coarsening}.
A faithful long-time surrogate should therefore recover both the characteristic domain-growth law and the statistical scaling behavior of the evolving morphology. Furthermore, since these statistics can be interpreted as indicators of faithful Cahn--Hilliard dynamics only when the prediction remains physically admissible, we retain the bound-conforming objective when introducing the structure-function regularization, yielding a structure-preserving FNO (spFNO). This added supervision evaluates whether explicit control of the evolving spectral morphology improves the quantitative 
characterization of the coarsening dynamics and scaling, when the bound-conforming objective of bcFNO is retained as the base training objective.

The main contribution of this paper is the development of a two-stage structure-preserving FNO framework for long-time Cahn--Hilliard prediction under coarse temporal supervision. Stage~I introduces bcFNO, which incorporates a bound-conforming objective to stabilize physically plausible long-time rollouts by discouraging
departures from the nominal phase-field range. Building on this formulation, Stage~II retains the bound-conforming objective and adds structure-function supervision, yielding spFNO. This nested formulation separates the control of energy stability from the correction of late-stage coarsening statistics. On the various spatial resolutions, temporal discretizations, training seeds, and temporal-sampling cases, the bound-conforming component provides the major stabilization of the predicted energy landscape, whereas the structure-function supervision further improves agreement in the normalized structure function \(s(k,t)\) in the predicted solutions. These results demonstrate that the nominal phase-field admissibility is an essential precondition for meaningful long-time coarsening characterization, while the explicit structure-function supervision provides a targeted correction of selected coarsening statistics.

The challenges we undertake in this work are not unique to Cahn--Hilliard dynamics. In many multiscale nonlinear systems, accurate short-time predictions do not guarantee physically meaningful long-time behavior. For example, turbulence is commonly characterized by inertial-range energy spectra and velocity-increment structure functions, whose scaling laws reflect cascade dynamics rather than diffusion-driven coarsening \cite{kolmogorov1941dissipation,frisch1996turbulence}. Recent SciML work has used these statistical quantities to evaluate generated two-dimensional turbulent fields, showing scaling exponents consistent with the expected Kolmogorov behavior \cite{whittaker2024turbulence}. These relevant studies suggest a broader design principle for operator learning: reliable long-time learning should combine constraints that promote physically interpretable states with statistical observables tailored to the governing equation and dynamical regime. 

The rest of this paper is organized as follows. Section~\ref{sec:ch_physics} introduces the Cahn--Hilliard equation, its continuous late-stage dynamic scaling, and the corresponding discrete spectral quantities needed in the numerical experiments. Section~\ref{sec:num} describes the reference solver, data-generation procedure, and FNO architecture. Section~\ref{sec:sp_fno} presents the bound-conforming objective and the structure-function regularization. Section~\ref{sec:experiments} reports numerical experiments on a Cahn--Hilliard-like equation with analytical solution and on general Cahn--Hilliard coarsening with random initial conditions. We then analyze the solutions of bcFNO and spFNO by examining the long-time stability, coarsening diagnostics, sampling sensitivity, and online costs. Finally, Section~\ref{sec:conclusion} summarizes the numerical results, discusses limitations, and indicates future research directions.

\section{Cahn--Hilliard Equation and Dynamic Scaling} \label{sec:ch_physics}

\subsection{Phase-Field Formulation and Cahn--Hilliard Equation}
We consider the Cahn--Hilliard equation, a phase-field model describing the evolution of a conserved order parameter $\phi(\mathbf{r},t)$
defined on a d-dimensional ($d\ge 2$) periodic domain $\Omega \subset \mathbb R^d$.
The evolution of $\phi(\mathbf{r},t)$ is formulated as the
$H^{-1}$ gradient flow governed by the free-energy functional
\begin{equation}
E[\phi]
=
\int_{\Omega}
\left(
F(\phi)
+
\frac{\kappa}{2} |\nabla \phi|^2
\right)\, \mathrm{d}\mathbf{r},
\end{equation}
where the parameter $\kappa > 0$ controls the interfacial energy and the diffuse interface width, and $F(\phi)$ is the bulk free-energy density, taken as the classical double-well potential
\begin{equation*}
F(\phi) = \frac{1}{4}{(\phi^2 - 1)}^2,
\end{equation*}
whose local minima at $\phi = \pm 1$ denote the coexisting stable bulk phases.

We consider the chemical potential $\mu$ that is defined as the variational derivative of the free-energy functional with respect to the order parameter
\[
\mu = \frac{\delta E}{ \delta \phi}.
\]
Evaluating the first variation of the free-energy functional and applying integration by parts under periodic boundary conditions yields
\begin{equation*}
\mu =  F'(\phi) - \kappa \nabla^2 \phi = \phi^3 - \phi - \kappa \nabla^2 \phi.
\end{equation*}
Combining the chemical potential with the $H^{-1}$ gradient-flow dynamics
\begin{equation*}
\frac{\partial \phi}{\partial t}
=
 \nabla \cdot \left(M\nabla \mu\right),
\label{equ:mu}
\end{equation*}
yields the classical Cahn--Hilliard equation,
\begin{equation}
\frac{\partial \phi}{\partial t}
=
M \nabla^2
\left(
\phi^3 - \phi - \kappa \nabla^2 \phi
\right),
\label{equ:ch}
\end{equation}
where $M>0$ is the constant mobility. This is a nonlinear fourth-order parabolic PDE that conserves the total mass. 
Furthermore, the $H^{-1}$ gradient-flow structure implies the energy dissipation law
\begin{equation}
\frac{dE}{dt}
=
-\int_{\Omega}
M|\nabla\mu|^2\,d\mathbf r
\le0,
\label{eq:energy_law}
\end{equation}
which shows that the total free energy decreases monotonically over time.
The dissipation and conservation laws constrain the global evolution of the Cahn--Hilliard system.
To quantify the evolving domain scale, statistical morphology, and late-stage dynamic scaling, we introduce correlation-based spectral quantities and their derived coarsening statistics.

\subsection{Correlation Functions and Dynamic Scaling}
\label{subsec:cor_scaling}

To characterize the evolving morphologies generated by Eq.~\eqref{equ:ch}, we introduce statistical observables rooted in spatial correlation analysis. 
Under periodic boundary conditions, integrating Eq.~\eqref{equ:ch} over the domain and applying the divergence theorem yields
\[
\frac{\mathrm d}{\mathrm dt}\int_\Omega \phi(\mathbf r,t)\,\mathrm d\mathbf r =\int_{\partial \Omega} \nabla \mu \cdot \boldsymbol{n}= 0,
\]
indicating the conservation of the total order parameter. The fluctuation field is defined as 
\begin{equation*}
\psi(\mathbf r,t):= \phi(\mathbf r,t) - \bar\phi,\quad \text{where}\quad
\bar\phi = \frac1{|\Omega|}\int_\Omega \phi(\mathbf r,t)\,\mathrm d\mathbf r.
\label{eq:fluctuation}
\end{equation*}
Assuming statistical homogeneity, the equal-time two-point
correlation function  is
\begin{equation}
G(\mathbf r,t)
=
\langle
\psi(\mathbf r_0+\mathbf r,t)
\psi(\mathbf r_0,t)
\rangle ,
\label{eq:G_cont}
\end{equation}
whose Fourier transform is called the structure function
\begin{equation*}
S(\mathbf k,t)
=
\int_\Omega
G(\mathbf r,t)
e^{-i\mathbf k\cdot\mathbf r}
\,d\mathbf r.
\label{eq:S_cont}
\end{equation*}
For statistically isotropic morphologies,
$G$ and $S$ depend only on
$r=|\mathbf r|$
and
$k=|\mathbf k|$,
respectively.
Hence we write
\[
G(\mathbf r,t)\equiv G(r,t),
\qquad
S(\mathbf k,t)\equiv S(k,t).
\]

In the late-stage coarsening regime, the morphology becomes statistically self-similar, and the dynamical scaling hypothesis states that all statistical properties depend on time only through a single characteristic length scale $L(t)$. Consequently, the correlation function assumes the scaling form
\begin{equation*}
G(r,t) = G_s\!\left(\frac{r}{L(t)}\right).
\label{eq:G_scaling}
\end{equation*}
To derive the scaling behavior of the structure function, we introduce the scaled variable $\mathbf u = \mathbf r/L(t)$ and the Fourier-space vector $\mathbf q = \mathbf k L(t)$. 
Then

\[
S(\mathbf k,t)
=
\int
G_s\!\left(\frac{|\mathbf r|}{L(t)}\right)
e^{-i\mathbf k\cdot\mathbf r}
\,d\mathbf r
=
L(t)^d
\int
G_s(|\mathbf u|)
e^{-i\mathbf q\cdot\mathbf u}
\,d\mathbf u .
\]
Since the system is statistically isotropic, the integral depends only on the dimensionless scalar variable
\[
q=|\mathbf q|=kL(t).
\]
We therefore introduce the universal scaling function
\[
\mathcal F(q)
=
\int
G_s(|\mathbf u|)
e^{-i\mathbf q\cdot\mathbf u}
\,d\mathbf u,
\]
which yields
\begin{equation}
S(k,t)
=
L(t)^d
\mathcal F\!\bigl(kL(t)\bigr).
\label{eq:S_scaling}
\end{equation}

Evaluating Eq.~\eqref{eq:G_cont} at zero separation (\(\mathbf r=\mathbf0\)) gives the variance of the fluctuation field,
\begin{equation}
\mathrm{var}(t)
=
G(\mathbf0,t)
=
G(0,t)
=
\langle {\psi}^2\rangle .
\label{eq:variance_cont}
\end{equation}
In the late-stage scaling regime, the bulk phases approach
\(\phi\approx\pm1\), while the interfacial volume fraction becomes asymptotically negligible. 
Consequently, the variance approaches a constant value.
It is therefore convenient to introduce the normalized structure function
\begin{equation}
s(k,t) = \frac{S(k,t)}{\mathrm{var}(t)}.
\label{eq:s_def_theory}
\end{equation}
Since normalization by the asymptotically constant variance does not alter the dependence on the scaling variable \(kL(t)\), the normalized structure function satisfies
\begin{equation}
s(k,t) = L(t)^d\,\widetilde{\mathcal F}\!\bigl(kL(t)\bigr),\qquad
\widetilde{\mathcal F}(q) = \frac{\mathcal F(q)}{G_s(0)}.
\label{eq:s_scaling}
\end{equation}

To characterize the location of spectral weight as a single scale statistic, we define the characteristic wavenumber as the first spectral moment of $s(k,t)$:
\begin{equation}
k_1(t) = \frac{\displaystyle\int k\,s(k,t)\,\mathrm dk}{\displaystyle\int s(k,t)\,\mathrm dk}.
\label{eq:k1_def}
\end{equation}
Inserting Eq.~\eqref{eq:s_scaling} into Eq.~\eqref{eq:k1_def} and changing variable to $q = kL(t)$ gives
\[
k_1(t) = \frac{1}{L(t)}\,\frac{\int q\widetilde{\mathcal F}(q)\,dq}{\int \widetilde{\mathcal F}(q)\,dq}
     \propto L(t)^{-1}.
\]
Hence
\begin{equation*}
k_1(t) \propto L(t)^{-1},
\label{eq:k1_scaling}
\end{equation*}
demonstrating that the first spectral moment provides an estimate of the characteristic inverse length scale.
Accordingly, throughout this work we use
\begin{equation}
L(t)\propto k_1^{-1}(t)
\label{eq:L_def}
\end{equation}
as the operational measure of the characteristic length; 
the proportionality constant is immaterial because only scaling relations are considered.

Using Eq.~\eqref{eq:k1_scaling}, the scaling variable $kL(t)$ is proportional to $k/k_1$.
Consequently, Eq.~\eqref{eq:s_scaling} can be rewritten as
\begin{equation}
k_1^{d}(t)\,s(k,t) = \mathcal G\!\left(\frac{k}{k_1(t)}\right),
\label{eq:collapse}
\end{equation}
where $\mathcal G$ is a time-independent scaling function. Hence, collapse of \(k_1^d(t) s(k,t)\) against \(k/k_1(t)\) tests whether the
late-stage spectra are governed by a single evolving length scale \cite{furukawa1985dynamic}.

While Eq.~\eqref{eq:collapse} characterizes the statistical self-similarity of the evolving morphology, the evolution of the characteristic length governs the coarsening kinetics. Specifically, in the diffusion-controlled late-stage regime of the constant-mobility Cahn--Hilliard equation, domain growth follows a Lifshitz--Slyozov mechanism \cite{lifshitz1961kinetics}, as derived below.

In the sharp-interface limit, the Gibbs-Thomson relation states that the chemical potential is proportional to the local interfacial curvature \cite{gibbs1878equilibrium}.
For a statistically self-similar morphology characterized by a single length scale $L(t)$, the typical interfacial curvature scales as $L^{-1}(t)$.
The Gibbs-Thomson relation therefore implies the chemical-potential scaling
\begin{equation*}
  \mu \sim  \frac{\sigma}{L(t)},
  \label{eq:gibbs}
\end{equation*}
where $\sigma$ represents the surface tension. Since the chemical
potential varies over the characteristic length scale \(L(t)\),
\begin{equation*}
|\nabla \mu| \sim \frac{\mu}{L(t)}
\sim \frac{\sigma}{L(t)^2}.
\label{eq:grad_mu}
\end{equation*}
By Fick's law,
\begin{equation*}
\mathbf{J}=-M\nabla\mu,
\end{equation*}
the net diffusive-flux magnitude therefore satisfies
\begin{equation}
  |\mathbf{J}| \sim \frac{M\sigma}{L(t)^2}.
  \label{eq:J1}
\end{equation}

Mass conservation relates the normal interface velocity to the net diffusive flux across the interface. 
A normal interface displacement $dL(t)$ changes the volume of a characteristic d-dimensional domain by $dV\sim L(t)^{d-1}dL(t)$.
The corresponding change in its conserved order parameter is therefore 
\[
\Delta\phi\, dV=\Delta\phi L(t)^{d-1}dL(t),
\]
where $\Delta\phi=\phi_+-\phi_-$.

\begin{figure}[!ht]
  \centering
  \includegraphics[width=\linewidth]{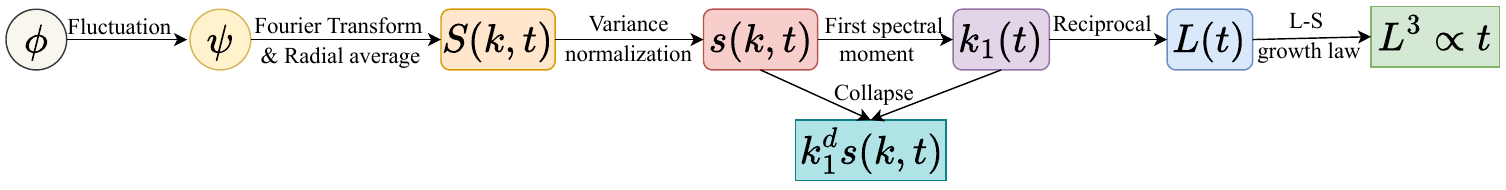}
 \caption{Logical chain of the theoretical analysis linking phase-field fluctuations to spectral statistics, the characteristic length, the macroscopic coarsening law, and dynamic-scaling collapse.}\label{fig:coarsening_evaluation_chain}
\end{figure}

The diffusive flux transfers order parameter through an interface of area 
$A(t)\sim L(t)^{d-1}$ at the rate $ |\mathbf{J}| A(t)$. Hence
\begin{equation*}
 |\mathbf{J}| A(t)\sim  \Delta \phi \frac{dV(t)}{dt}.
\end{equation*} 
After canceling the common interfacial-area factor, we obtain
\begin{equation}
 |\mathbf{J}| \sim  \Delta \phi \frac{dL(t)}{dt}.
  \label{eq:J2}
\end{equation} 
Substituting Eq.~\eqref{eq:J1} into Eq.~\eqref{eq:J2} yields 
\begin{equation}
  \frac{dL(t)}{dt}\sim\frac{M\sigma}{\Delta\phi}\frac{1}{L(t)^2} = \frac{C}{L(t)^2},
  \label{eq:dLdt}
\end{equation}
where $C=\frac{M\sigma}{\Delta\phi}$ is a constant.
Integrating Eq.~\eqref{eq:dLdt} yields
\begin{equation}
  L^3(t) - L^3(t_{0}) \sim 3C(t-t_{0}).
  \label{eq:LS_law}
\end{equation}
Therefore, $L^3(t)$ is expected to be approximately affine in time, $L^3(t) \approx at+b$, over a late-stage coarsening regime. This approximation is independent of dimensionality of the Cahn--Hilliard equation for $d \ge 2$. 
We use the growth law as a key benchmark for evaluating the long-time dynamics predicted by learned surrogate models.
The sequence connecting the phase field to this macroscopic growth evaluation is summarized in Fig.~\ref{fig:coarsening_evaluation_chain}.

\subsection{Discrete Spectral Representation of the Scaling Framework}
\label{subsec:discrete}

In order to apply the continuous scaling framework developed in Sec.~\ref{subsec:cor_scaling} to numerical simulations, we now formulate the statistical observables in a discrete spectral framework and calculate them using the numerical solutions of Eq.~\eqref{equ:ch}.
The computational domain is discretized on a uniform periodic $N^d$ lattice with unit spacing ($\Delta x=1$ in all directions), and all spectral quantities are computed using the Fast Fourier Transform (FFT).

The discrete Fourier transform of the fluctuation field is defined as 
\begin{equation*}
\widehat{\psi}(\mathbf k,t) = \sum_{\mathbf r} \psi(\mathbf r,t)\,e^{-i\mathbf k\cdot\mathbf r},
\label{eq:fft_def}
\end{equation*}
where the summation extends over all grid points in the computational domain. For statistically homogeneous fields, the discrete Wiener-Khinchin theorem establishes the equivalence between the discrete Fourier transform of the two-point correlation function and the power spectrum of the fluctuation field \cite{wiener1930generalized,khintchine1934korrelationstheorie}.
This avoids explicit evaluation of the two-point correlation function while remaining mathematically equivalent.
Accordingly, the structure function is computed directly from the squared Fourier amplitudes,
\begin{equation*}
S(\mathbf k,t) = \frac1{N^d}\big\langle |\widehat{\psi}(\mathbf k,t)|^2\big\rangle.
\label{eq:S_discrete}
\end{equation*}
Since only a single statistically homogeneous realization is available during evaluation, the structure function is estimated from the corresponding periodogram. 
Assuming statistical isotropy, radial shell averaging is subsequently performed to obtain the one-dimensional spectrum.

Applying Parseval's identity to the discrete Fourier transform yields
\[
\sum_{\mathbf r}|\psi(\mathbf r,t)|^2 = \frac1{N^d}\sum_{\mathbf k}|\widehat{\psi}(\mathbf k,t)|^2.
\]
Consequently, the variance of the fluctuation field, corresponding to the zero-separation correlation, is given by
\begin{equation}
\mathrm{var}(t) = \frac1{N^d}\sum_{\mathbf r}\psi(\mathbf r,t)^2
               = \frac1{N^{2d}}\sum_{\mathbf k}|\widehat{\psi}(\mathbf k,t)|^2.
\label{eq:variance}
\end{equation}
Consistent with the normalized structure function introduced in Eq.~\eqref{eq:s_def_theory}, the discrete structure function is defined as
\begin{equation}
s(\mathbf k,t) = \frac{S(\mathbf k,t)}{\mathrm{var}(t)},
\label{eq:s_discrete}
\end{equation}
which satisfies
\[
\sum_{\mathbf k}s(\mathbf k,t)=N^d,
\]
by construction.

Assuming statistical isotropy, the two-dimensional spectrum is reduced to a one-dimensional radial spectrum by averaging Fourier modes with similar wavenumber magnitude.
Specifically, the radial shell associated with the scalar wavenumber $k$ is defined by
\[
k-\frac{\Delta k}{2}
\le
|\mathbf k|
<
k+\frac{\Delta k}{2},
\]
where $\Delta k$ is the width of the radial wavenumber bin.
The shell-averaged normalized structure function is therefore defined by 
\begin{equation}
s(k,t) = \frac1{\mathcal N(k)}\sum_{|\mathbf k|\in k\text{-shell}} s(\mathbf k,t),
\label{eq:s_shell}
\end{equation}
where $\mathcal N(k)$ denotes the number of Fourier modes contained in the corresponding shell.
The shell-averaged spectrum $s(k,t)$ constitutes the discrete counterpart of the continuous normalized structure function introduced in Sec.~\ref{subsec:cor_scaling}.
It serves as the primary spectral descriptor from which all subsequent statistical quantities are derived.

The discrete characteristic wavenumber is computed from the first spectral moment of $s(k,t)$
\begin{equation}
k_1(t) = \frac{\displaystyle\sum_k k\,s(k,t)}{\displaystyle\sum_k s(k,t)}.
\label{eq:k1_discrete}
\end{equation}
The characteristic length is then defined consistently with Eq.~\eqref{eq:L_def} as
\begin{equation}
  L(t) = \frac{1}{k_1(t)}.
\label{eq:L_discrete}
\end{equation}
For late-stage coarsening, the Lifshitz--Slyozov law $L(t)\sim t^{1/3}$ implies an
approximately affine dependence of $L^3(t)$ on $t$ over a finite fitting
window. We therefore fit
\begin{equation}
L^3(t)\approx at+b,
\label{eq:discrete_L3_linear}
\end{equation}
where $a$ is the fitted late-time slope and $b$ is the intercept.

The normalized structure function defined in Eq.~\eqref{eq:s_shell} is used to examine dynamic scaling with
the parameter
\begin{equation}
\mathcal{S}_{c}(q;t)
=
k_1(t)^d s(k,t),
\qquad
q=\frac{k}{k_1(t)} .
\label{eq:collapse_discrete}
\end{equation}
If dynamic scaling holds, spectra evaluated at different late times collapse approximately onto a common curve when expressed in terms of $q$ and $\mathcal{S}_{c}$ \cite{Chakrabarti1993,Zhu1999}.

\begin{figure}[!ht]
  \centering
  \begin{subfigure}[t]{0.32\textwidth}
    \centering
    \includegraphics[width=\linewidth]{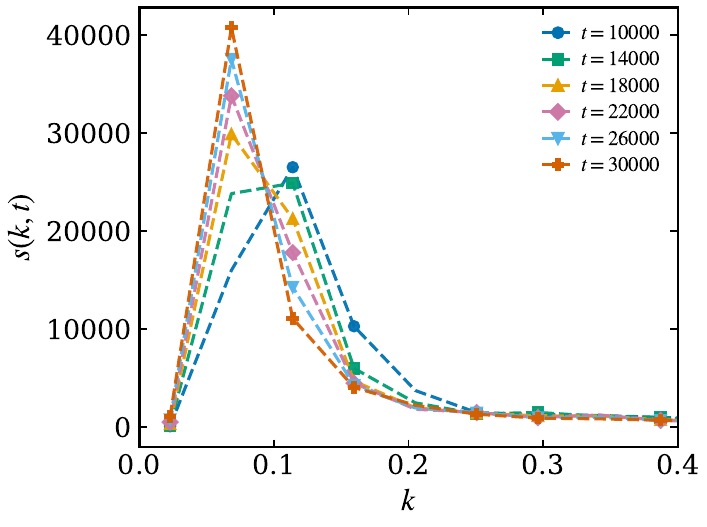}
  \end{subfigure}
  \hfill
  \begin{subfigure}[t]{0.32\textwidth}
    \centering
    \includegraphics[width=\linewidth]{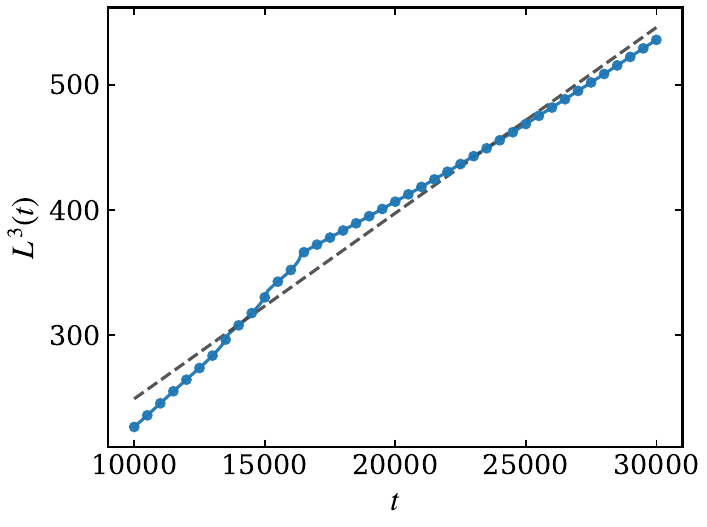}
  \end{subfigure}
  \hfill
  \begin{subfigure}[t]{0.32\textwidth}
    \centering
    \includegraphics[width=\linewidth]{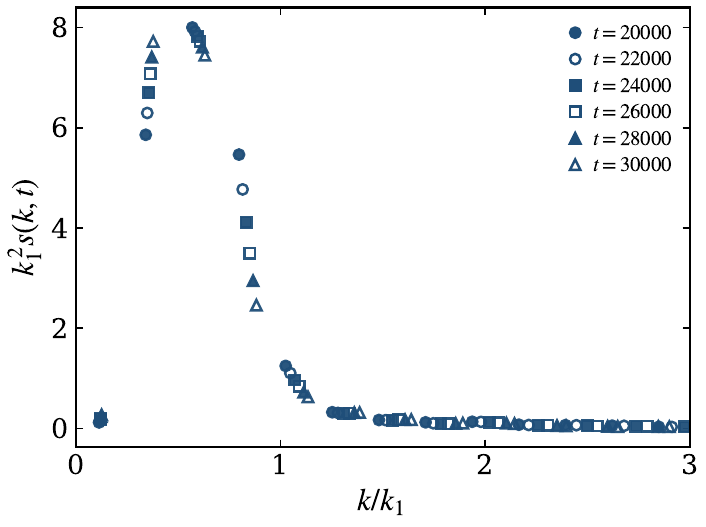}
  \end{subfigure}
 \caption{Late-stage coarsening statistics: $s(k,t)$, $L^3(t)$, and the dynamic-scaling collapse of the reference trajectory. The dashed line is a least-squares fit of $L^3(t)$ over the late-stage regime (middle). In the collapse plot (right), the zero mode is isolated and  interpolation is restricted to $q_{\min}>0$.}
  \label{fig:theory_scaling_summary} 
\end{figure}

Figure~\ref{fig:theory_scaling_summary} illustrates these discrete evaluations for the 2-D Cahn--Hilliard equation. The normalized radial spectrum shifts toward lower wavenumbers (left), indicating growth of the characteristic domain scale. 
An affine function is fitted to \(L^3(t)\) over the late-stage coarsening regime (middle).
The approximately collapsed rescaled spectra support
late-stage statistical self-similarity (right).
For the discrete collapse plot, the zero mode is shown separately as $(q,\mathcal{S}_{c})=(0,0)$. Because the mean-subtracted field has no zero Fourier mode, the first positive-wavenumber sample occurs at
\[
q_{\min}(t)=\frac{k_{\min}}{k_1(t)},
\qquad
k_{\min}=\frac{2\pi}{N}.
\]
Thus, no spectral samples exist between the isolated origin and $q_{\min}(t)$. Interpolation is performed only over the available positive-wavenumber samples, so the blank interval adjacent to the origin is a discrete-resolution effect and does not represent a missing portion of the collapse curve.

These discrete quantities provide the statistical diagnostics used to compare the reference trajectory with the learned neural operators throughout this work. 
In particular, the normalized structure function determines the characteristic length, the late-stage coarsening kinetics, and the dynamic-scaling collapse, and thereby motivates the statistical constraint introduced in Sec.~\ref{sec:sp_fno}.

\section{Numerical Methods for Data Generation and Training}
\label{sec:num}

\subsection{Numerical Solver}\label{subsec:num_solver}

All training and evaluation trajectories are generated using a first-order
semi-implicit Fourier pseudospectral discretization of the continuous
2-D Cahn--Hilliard equation~\eqref{equ:ch} with periodic boundary conditions.
The mobility is set to $M=1$ and the gradient-energy coefficient to
$\kappa=0.5$. Initial conditions are generated by perturbing a homogeneous
zero-mean state with Gaussian noise of standard deviation $0.01$. The
nonlinear term $f(\phi)=\phi^3-\phi$ is evaluated pointwise in physical space,
whereas spatial derivatives are evaluated in Fourier space.

To address the stiffness associated with the fourth-order diffusion operator, the linear biharmonic term is treated implicitly, whereas the nonlinear term is treated explicitly. 
Let $\delta t$ denote the internal solver time step and $\phi^n$ the numerical solution at time $t_n$. 
The semi-implicit Euler discretization is
\begin{equation*}
\frac{\phi^{n+1}-\phi^n}{\delta t}
=
M \nabla^2 f(\phi^n)
-
M\kappa \nabla^4 \phi^{n+1}.
\end{equation*}
Equivalently,
\begin{equation}
\left(I + \delta t\,M\kappa \nabla^4\right)\phi^{n+1}
=
\phi^n + \delta t\,M \nabla^2 f(\phi^n).
\label{eq:discrete_ch}
\end{equation}

Applying the discrete Fourier transform to Eq.~\eqref{eq:discrete_ch}, with
$\mathcal{F}(\nabla^2)=-K^2$,
$\mathcal{F}(\nabla^4)=K^4$, and
$K^2=k_x^2+k_y^2$, gives the Fourier-space update
\begin{equation}
\widehat{\phi}^{\,n+1}(\mathbf{k})
=
\frac{
\widehat{\phi}^{\,n}(\mathbf{k})
-
\delta t\,M K^2 \widehat{f(\phi^n)}(\mathbf{k})
}{
1 + \delta t\,M\kappa K^4
}.
\label{equ:solver}
\end{equation}
Equation~\eqref{equ:solver} is evaluated using the fast Fourier transform \cite{soares2023}.

\subsection{Coarse-Supervision Data Generation}
\label{subsec:data_gen}

These experiments vary the internal solver time step used to generate the numerical reference trajectories, while keeping the initial conditions and spatial grid fixed. 
The internal time step is
\[
\delta t\in\{1,\;0.5,\;0.2,\;0.1\}.
\]
For the four values of $\delta t$, states are stored every
$m\in\{20,\;40,\;100,\;200\}$  solver steps, giving
\begin{equation*}
    \Delta T=m\,\delta t=20.
    \label{eq:coarse_interval}
\end{equation*}
Hence, all models are trained and evaluated with the same interval between consecutive predicted states, \(\Delta T=20\), while the reference trajectories use different internal solver time steps.

We use $\Delta T=20$ throughout. 
Among the tested solver settings, $\delta t=1$ is the largest step that produced stable reference trajectories under the present discretization, whereas $\delta t=0.1$ is the smallest.
For each $(N,\delta t)$, we generate seven trajectories from
independently sampled initial conditions. 
Six are used to form training pairs, and the seventh is used only for long-time rollout evaluation. 
Neither its initial condition nor any state on that trajectory is included in training.
For the held-out reference trajectory, the coarsening evaluation mask is constructed from the saved reference states and stored.
The same mask is then applied to the neural-operator prediction. 
We use the resulting coarse-step state pairs directly as the training data. 
The following section specifies the neural-operator architecture.

\subsection{Fourier Neural Operator}
\begin{figure}[!ht]
  \centering
  \includegraphics[width=\linewidth]{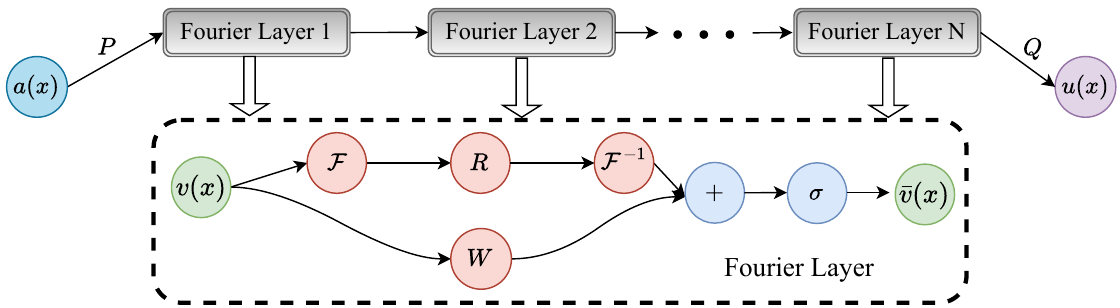}
 \caption{Schematic of the Fourier neural operator. The network consists of a lifting operator $P$, a stack of Fourier layers, and a projection operator $Q$.}\label{fig:fno}
\end{figure}

The Fourier neural operator (FNO)~\cite{Li2020} is used to parameterize the coarse time-step map over the interval $\Delta T$.
It learns nonlinear mappings between function spaces while capturing long-range, nonlocal dependencies through kernels parameterized in the frequency domain.
A schematic illustration of the architecture is provided in Fig.~\ref{fig:fno}.

As illustrated in Fig.~\ref{fig:fno}, the lifting operator $P$, implemented as a neural network, first lifts the input channels to a higher-dimensional latent representation. 
This latent representation is processed by a sequence of Fourier layers, which constitute the core of the model. 
In each layer, the integral kernel operator is parameterized in Fourier space, where nonlocal dependencies are learned via element-wise complex multiplications over the retained Fourier modes. Specifically, the layer consists of two parallel branches. 
The spectral branch applies a learnable linear transform $R$ to the truncated low-frequency modes of the Fourier-transformed input, followed by an inverse Fourier transform $\mathcal{F}^{-1}$. 
In the parallel branch, a pointwise linear transform $W$ is applied to the latent features at each spatial grid point. 
The outputs of the two branches are summed and passed through a nonlinear activation $\sigma$. 
The layer-$l$ update is thus given by 
\begin{equation}
v_{l+1}(x)
=
\sigma\left(
W v_l(x)
+
\mathcal{F}^{-1}
\big(
R \cdot \mathcal{F}(v_l)
\big)(x)
\right).
\end{equation}
By stacking multiple such layers, the network learns a nonlinear mapping between function spaces. Finally, the projection operator $Q$ decodes the final latent state to yield the predicted solution at the next time step, denoted by $\phi^{\mathrm{pred}}_{n+1}$.
The one-step state loss used to train this map is specified with the remaining regularization terms in Sec.~\ref{sec:sp_fno}.

\section{Structure-Preserving Fourier Neural Operator}\label{sec:sp_fno}

The FNO parameterizes the learned CH evolution map over the coarse interval $\Delta T$.
Building on this parameterization, we introduce a two-stage structure-preserving training procedure.  Stage~I first trains a bound-conforming FNO by augmenting this loss with a soft bound penalty. Stage~II then continues from the Stage~I model and adds structure-function regularization to correct the late-stage coarsening statistics.

\subsection{Bound-Conforming Penalty}

For the quartic double-well potential used here, the homogeneous states $\phi=\pm1$ are two minima.
In general, the Cahn--Hilliard equation does not satisfy a maximum principle that would preserve $|\phi|\le1$.
Accordingly, $|\phi|\le1$ is used as a nominal phase-field range, rather than as a universal invariant of the Cahn--Hilliard dynamics.

During long-time rollout, the standard FNO predictions can depart from this nominal range.
To penalize these departures, we introduce a soft pointwise penalty
\begin{equation}
L_{\mathrm{bc}}
=
\frac{1}{N^2}
\sum_{i=1}^{N}\sum_{j=1}^{N}
\left[
{\operatorname{softplus}(\phi(x)-1;\beta)}^{2}
+
{\operatorname{softplus}({-1}-\phi(x);\beta)}^{2}
\right].
\end{equation}
Here,
\begin{equation}
\operatorname{softplus}(z;\beta)
=
\frac{1}{\beta}\log\!\left(1+\exp(\beta z)\right),
\end{equation}
with $\beta=5$ in all experiments.

\subsection{Structure-Function Regularization}

The spFNO is initialized from the trained bcFNO and is optimized with an additional structure-function term for samples in the coarsening regime. 
Let $s^{\mathrm{pred}}(k_b,t)$ and $s^{\mathrm{ref}}(k_b,t)$
denote the normalized structure functions of the predicted and reference phase-fields, respectively, computed as described in Section~\ref{subsec:discrete}.
The normalized structure function $s(k,t)$ retains neither Fourier phase nor directional spectral information, its regularization is used together with the state loss.

\subsubsection{Spectral-shape term}

Let $b=1,\ldots B$ index the retained radial wavenumber bins, with $B=64$ in the experiments.
For a training sample at time $t$, define
\begin{equation*}
\bar{s}^{\mathrm{pred}}(k_b,t)
=
\frac{s^{\mathrm{pred}}(k_b,t)}
{\operatorname{stopgrad}\!\left[
\max_{1\leq j\leq B}s^{\mathrm{pred}}(k_j,t)
\right]+\varepsilon},
\qquad
\bar{s}^{\mathrm{ref}}(k_b,t)
=
\frac{s^{\mathrm{ref}}(k_b,t)}
{\max_{1\leq j\leq B}s^{\mathrm{ref}}(k_j,t)+\varepsilon},
\end{equation*}
where $\operatorname{stopgrad}$ denotes stop-gradient.
The shape term is 
\begin{equation}
\ell_{\mathrm{shape}}(t)
=
\frac{1}{\sqrt{B}}
\sum_{b=1}^{B}
\left(
\bar{s}^{\mathrm{pred}}(k_b,t)
-
\bar{s}^{\mathrm{ref}}(k_b,t)
\right)^2.
\end{equation}

\subsubsection{Peak-bin and low-intensity-bin terms}

Let
\begin{equation*}
s_{\max}^{\mathrm{ref}}(t)
=
\max_{1\leq b\leq B}s^{\mathrm{ref}}(k_b,t).
\end{equation*}
The peak and low-intensity bins are defined from the reference $s(k,t)$ as
\begin{equation*}
\mathcal{P}(t)
=
\left\{
b:\;
s^{\mathrm{ref}}(k_b,t)>
f_{\mathrm{peak}}s_{\max}^{\mathrm{ref}}(t)
\right\},
\end{equation*}
and
\begin{equation*}
\mathcal{N}_0(t)
=
\left\{
b:\;
s^{\mathrm{ref}}(k_b,t)<
\tau_{\rho}s_{\max}^{\mathrm{ref}}(t)
\right\},
\qquad
\mathcal{N}(t)
=
\mathcal{N}_0(t)\setminus\mathcal{P}(t).
\end{equation*}
The corresponding terms are
\begin{equation}
\ell_{\mathrm{peak}}(t)
=
\frac{1}{|\mathcal{P}(t)|}
\sum_{b\in\mathcal{P}(t)}
\left(
\frac{s^{\mathrm{pred}}(k_b,t)}
{s^{\mathrm{ref}}(k_b,t)}
-1
\right)^2,
\end{equation}
and
\begin{equation}
\ell_{\mathrm{null}}(t)
=
\frac{1}{|\mathcal{N}(t)|}
\sum_{b\in\mathcal{N}(t)}
\frac{
\left(
s^{\mathrm{pred}}(k_b,t)-
s^{\mathrm{ref}}(k_b,t)
\right)^2
}{
\max\!\left(
s^{\mathrm{ref}}(k_b,t),
\eta s_{\max}^{\mathrm{ref}}(t)
\right)
}.
\end{equation}

For one training sample, the structure-function loss is
\begin{equation}
\ell_{\mathrm{sk}}(t)
=
w_{\mathrm{shape}}\ell_{\mathrm{shape}}(t)
+
w_{\mathrm{peak}}\ell_{\mathrm{peak}}(t)
+
w_{\mathrm{null}}\ell_{\mathrm{null}}(t).
\end{equation}
Each term is evaluated separately for each training sample and then averaged over the mini-batch.

\subsection{Overall Training Objective}

For a mini-batch of $R$ training pairs, indexed by $r=1,\ldots, R$, let $t_{r}$ denote the input time of sample $r$. 
The state loss for one sample is
\begin{equation}
\ell_{\mathrm{state}}
=
\frac{1}{N^2}
\sum_{i=1}^{N}\sum_{j=1}^{N}
\left(
\phi^{\mathrm{pred}}_{t+\Delta T,i,j}
-
\phi^{\mathrm{ref}}_{t+\Delta T,i,j}
\right)^2,
\end{equation}
Let
\begin{equation*}
\chi_{\mathrm{c}}(t)
=
\begin{cases}
1, & t \text{ belongs to the coarsening regime},\\
0, & \text{otherwise}.
\end{cases}
\label{eq:coarsening_gate}
\end{equation*}
The mini-batch objective is
\begin{equation}
L_{\mathrm{total}}
=
\frac{1}{R}
\sum_{r=1}^{R}
\left[
\ell_{\mathrm{state}}(t_r)
+
\lambda_{\mathrm{bc}}\ell_{\mathrm{bc}}(t_r)
+
\chi_{\mathrm{c}}(t_r+\Delta T)
\lambda_{\mathrm{sk}}\ell_{\mathrm{sk}}(t_r+\Delta T)
\right].
\end{equation}
The state and bound terms are evaluated for every sample.
The structure-function term is evaluated only when \(\chi_{\mathrm{c}}(t+\Delta T)=1\). 
For all spFNO models, this term is omitted during the first training epoch and included from the second epoch onward. 
The coarsening indicator is used during training only and is not required for long-time rollout. Table~\ref{tab:loss_cases} summarizes these objectives and the samples on which each term is active.

The loss weights, spectral thresholds, and sampling settings are determined using the training trajectories described in Sec.~\ref{subsec:data_gen}. 
The held-out numerical trajectory is excluded from hyperparameter selection.
Table~\ref{tab:training_parameters} lists the selected values and search ranges.
The resulting objective combines state-level accuracy with amplitude regularization and statistical consistency.

\begin{table}[!t]
\centering
\caption{Training objectives and loss-activation rules.}
\label{tab:loss_cases}
\small
\begin{tabular}{p{0.08\linewidth}p{0.25\linewidth}p{0.55\linewidth}}
\toprule
Model & Objective & Active term \\
\midrule
FNO & $L_{\mathrm{state}}$ & State loss: all samples \\
bcFNO & $L_{\mathrm{state}}+\lambda_{\mathrm{bc}}L_{\mathrm{bc}}$
& State and bound penalty: all samples \\
spFNO &
$L_{\mathrm{state}}+\lambda_{\mathrm{bc}}L_{\mathrm{bc}}
+\chi_{\mathrm{c}}\lambda_{\mathrm{sk}} L_{\mathrm{sk}}$
& State and bound terms: all samples; \newline Structure-function term: coarsening only after epoch 1 \\
\bottomrule
\end{tabular}
\end{table}

\section{Experiments}\label{sec:experiments}

\subsection{Experimental Setup}\label{subsec:exp_setup}

We use two numerical experiments that share the same FNO architecture, the coarse prediction interval, and the same bound-conforming formulation, but differ in their governing equations, reference trajectories, and evaluations.

The first problem is a forced Cahn--Hilliard equation with an exact solution.
Its prescribed trajectory gives a unique pointwise reference for direct evaluation of accumulated field errors.
This problem is used primarily to assess the bound-conforming penalty in long-time field prediction.
The comparison therefore focuses on the FNO and bcFNO, with relative field errors and the phase-field amplitude.

The second problem is a source-free Cahn--Hilliard equation.
Each training or inference trajectory is initialized from an independently generated random phase-field perturbation. 
These random initial fields produce distinct coarsening morphologies and the held-out trajectories are generated independently of those used for training.
Small errors in interface position can then lead to substantial differences in domain location and topology during long-time prediction. 
Consequently, pointwise field discrepancies alone do not characterize late-stage coarsening fidelity.
For this problem, we evaluate the predicted trajectories using the free-energy evolution, the phase-field amplitude, the normalized structure function, the characteristic length scale, and the dynamic scaling.
The FNO, bcFNO, and spFNO are compared over long-time predictions.

Both problems use two-dimensional periodic fields.
The manufactured reference trajectory is evaluated directly from its analytical solution, whereas the source-free reference trajectories are generated using the numerical solver described in Sec.~\ref{subsec:num_solver}.
The data generation, model comparisons, and analysis of these quantities are presented in the following subsections.

\begin{figure}[!t]
  \centering
  \includegraphics[width=0.8\textwidth]{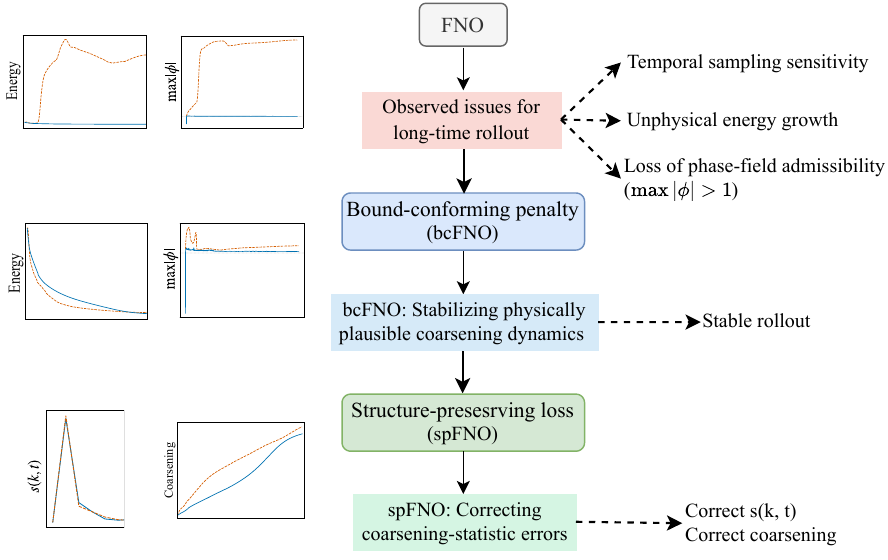}
   \caption{Two-stage training framework for Cahn--Hilliard coarsening dynamics. Stage~I trains a bound-conforming FNO (bcFNO), and Stage~II incorporates structure-function regularization to obtain a structure-preserving FNO (spFNO).}
  \label{fig:framework}
\end{figure}

Figure~\ref{fig:framework} summarizes the two-stage learning framework and its experimental evaluation.
The manufactured problem provides a fieldwise comparison of the Stage~I bound-conforming penalty, while the source-free coarsening problem is used to evaluate the complete framework, including the Stage~II structure-function regularization.

\FloatBarrier
\subsection{Forced Cahn--Hilliard Equation with an Exact Solution}
\label{subsec:manufactured_solution}

In source-free Cahn--Hilliard coarsening, a long-time FNO prediction can develop  domain distributions that are spatially displaced from the numerical reference.
A long-time field error then reflects both errors in the learned evolution and mismatch between the two morphologies.
To provide a uniquely prescribed reference state $\phi^\star(\cdot,t_n)$ at every iteration time $t_n$ for direct field-level comparison,
we first use a forced Cahn--Hilliard problem with a prescribed exact trajectory.

This problem is used only for direct field-level comparison. 
The prescribed forcing does not describe spontaneous phase separation or emergent Cahn--Hilliard coarsening. Its spectral evolution is therefore not interpreted using coarsening statistics or a scaling law.
The manufactured trajectory remains within a range narrower than the nominal interval $[-1,1]$. 
We nevertheless use the same bcFNO penalty as in the source-free experiments, without modifying its definition for the forced problem. Thus, the comparison tests whether the penalty suppresses large amplitude excursions generated during long-time prediction even when the exact trajectory itself does not approach the nominal range boundary. 
The interval $[-1,1]$ remains a soft reference in this test.
On the periodic square $\mathbb{T}_m^2=[0,m)^2$, we consider
\begin{equation}
  \frac{\partial\phi}{\partial t}
  =
  M\Delta\left(\phi^3-\phi-\kappa\Delta\phi\right)
  +S(\mathbf{x},t).
  \label{eq:forced_ch}
\end{equation}
For a prescribed smooth field $\phi^\star(\mathbf{x},t)$, the source term is computed analytically as
\begin{equation*}
  S^\star(\mathbf{x},t)
  =
  \partial_t\phi^\star(\mathbf{x},t)
  -
  M\Delta\left[
    (\phi^\star)^3-\phi^\star-\kappa\Delta\phi^\star
  \right],
  \label{eq:manufactured_source}
\end{equation*}
so that $\phi^\star$ is an exact solution of
Eq.~\eqref{eq:forced_ch}.
Its complete Fourier construction is given in~\ref{app:manufactured_trajectory}.

We use $m=256$, $M=1$, $\kappa=0.5$, and a uniform $128\times128$ grid ($\Delta x=2$). 
The trajectory is sampled analytically at the coarse interval $\Delta T=20$. 
With the source evaluated at the two Gauss--Legendre nodes $c_{1,2}=1/2\pm\sqrt{3}/6$,  the learned map is
\begin{equation*}
\left[
\phi^\star(\cdot,t_n),
S^\star(\cdot,t_n+c_1\Delta T),
S^\star(\cdot,t_n+c_2\Delta T)
\right]
\longmapsto
\phi^\star(\cdot,t_n+\Delta T).
\label{eq:manufactured_map}
\end{equation*}
During prediction, the current prediction is combined with the prescribed source values at the corresponding nodes.
We compute the relative field error at prediction time $t_n$ by
\begin{equation*}
e_{\rm field}(t_n)
=
\frac{\|\phi_{\rm pred}(\cdot,t_n)-\phi^\star(\cdot,t_n)\|_2}
{\|\phi^\star(\cdot,t_n)\|_2}.
\label{eq:manufactured_field_error}
\end{equation*}

\begin{figure}[!t]
  \centering
  \setlength{\tabcolsep}{1pt}
  \begin{tabular}{@{}ccccc@{}}
    \includegraphics[width=0.19\linewidth]{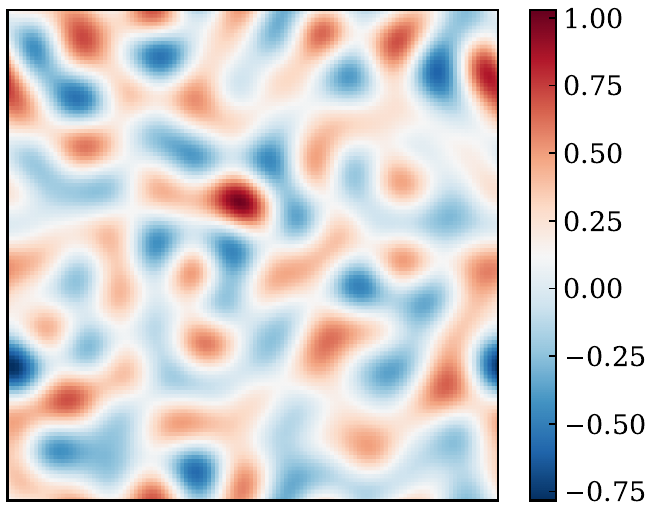} &
    \includegraphics[width=0.19\linewidth]{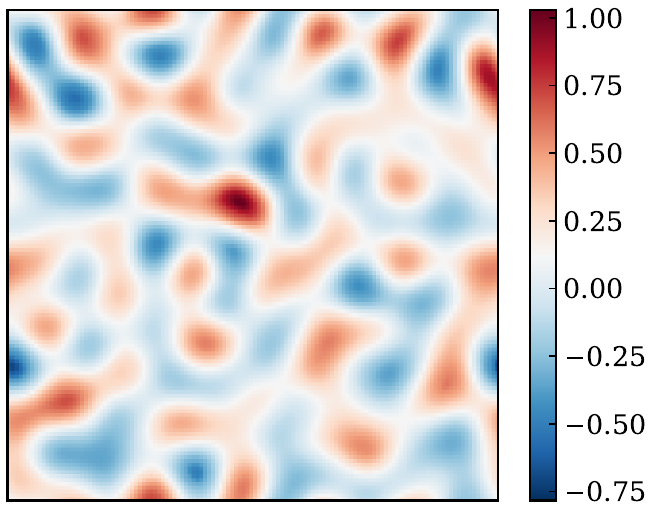} &
    \includegraphics[width=0.19\linewidth]{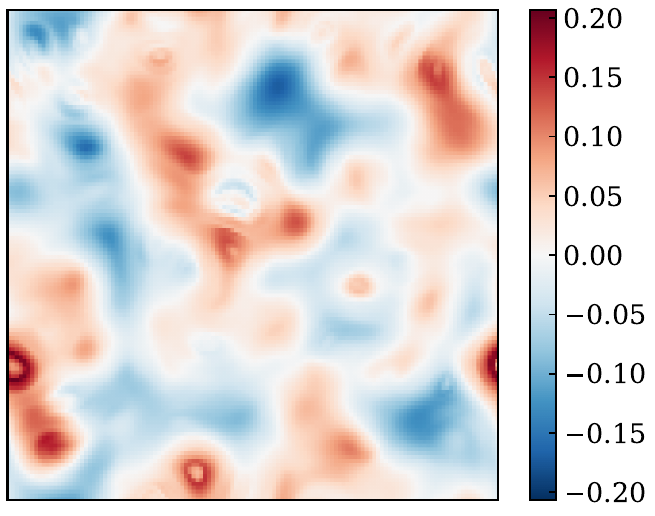} &
    \includegraphics[width=0.19\linewidth]{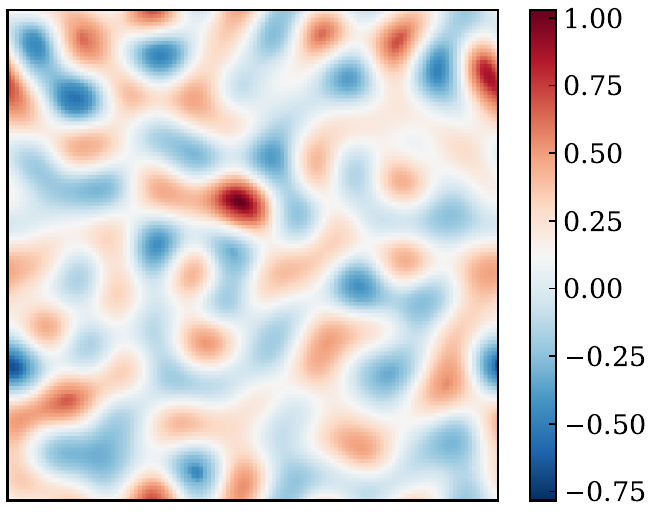} &
    \includegraphics[width=0.19\linewidth]{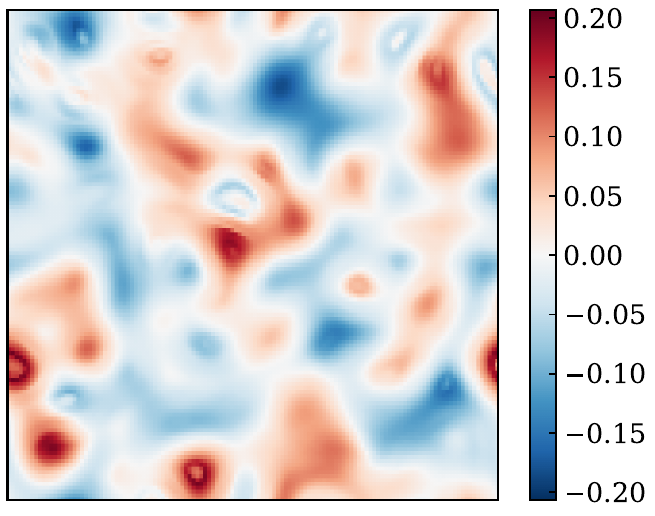} \\
    \includegraphics[width=0.19\linewidth]{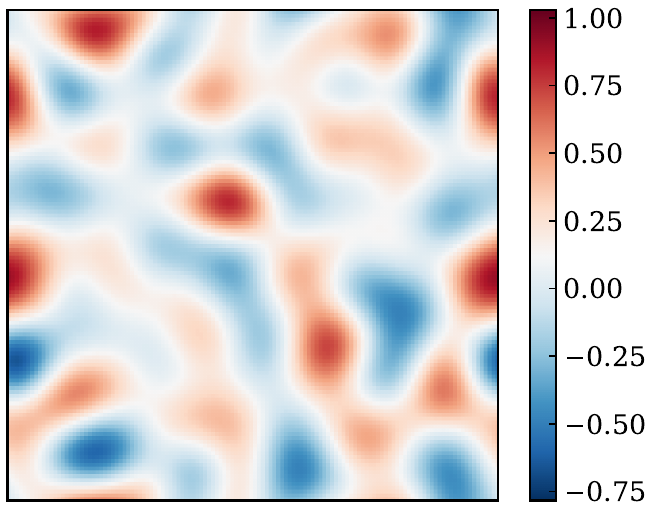} &
    \includegraphics[width=0.18\linewidth]{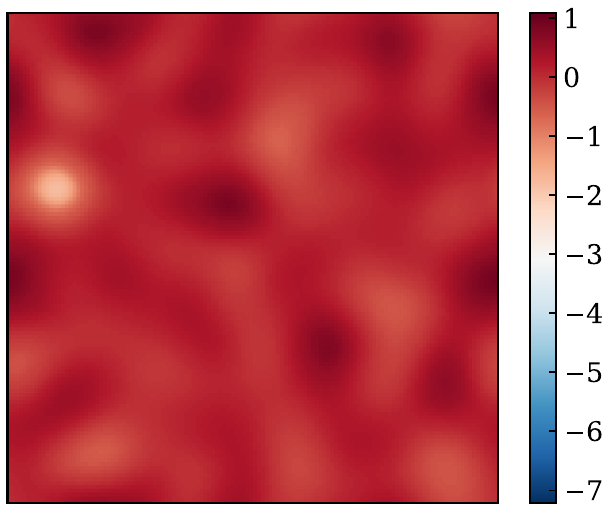} &
    \includegraphics[width=0.185\linewidth]{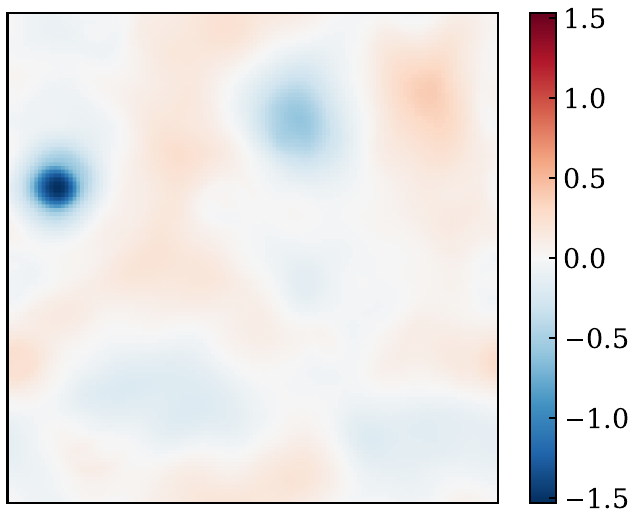} &
    \includegraphics[width=0.19\linewidth]{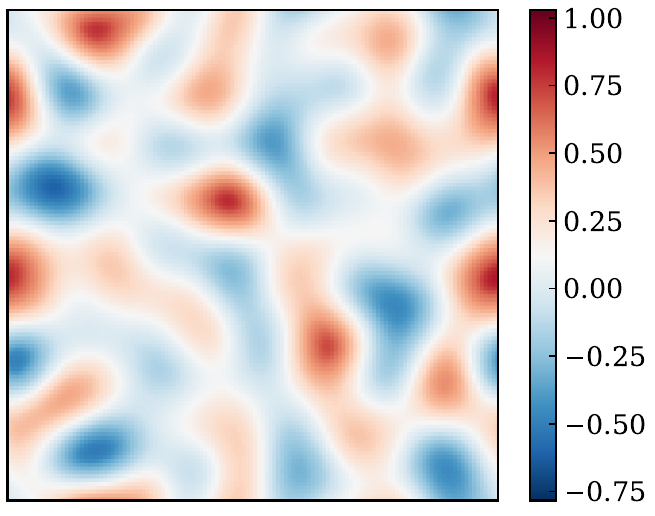} &
    \includegraphics[width=0.185\linewidth]{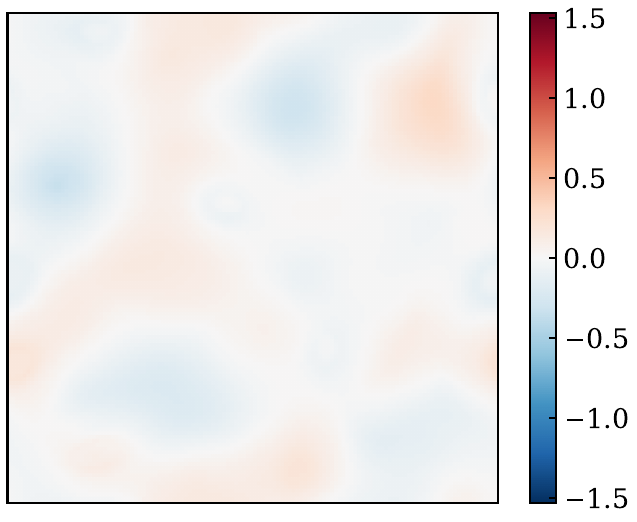} \\
    \includegraphics[width=0.19\linewidth]{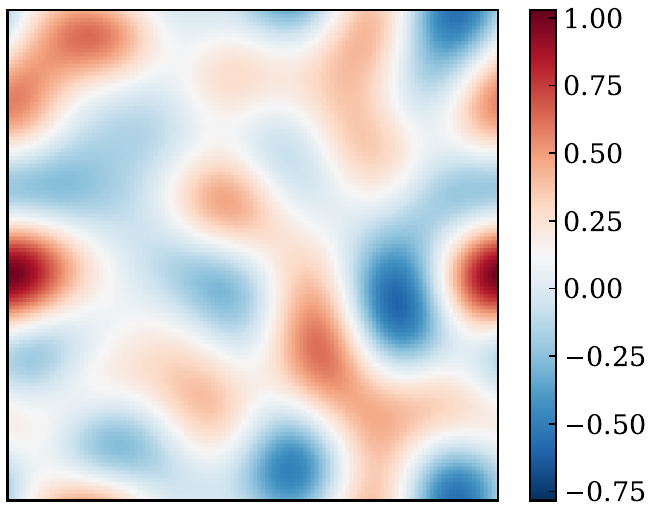} &
    \includegraphics[width=0.18\linewidth]{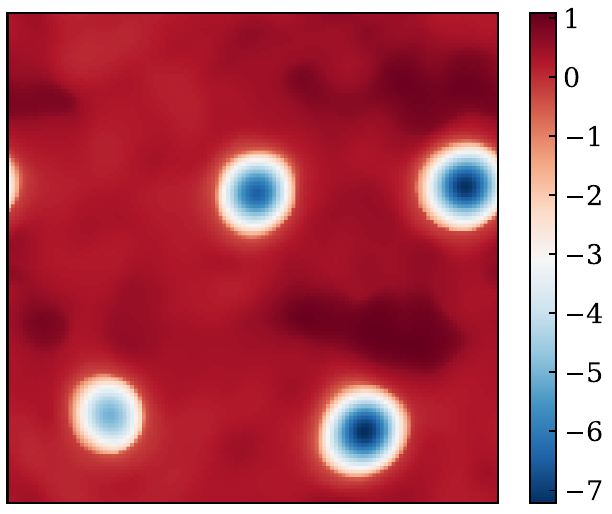} &
    \includegraphics[width=0.18\linewidth]{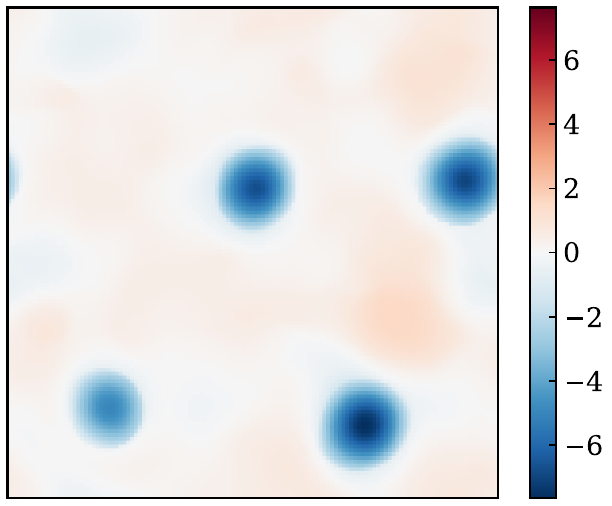} &
    \includegraphics[width=0.19\linewidth]{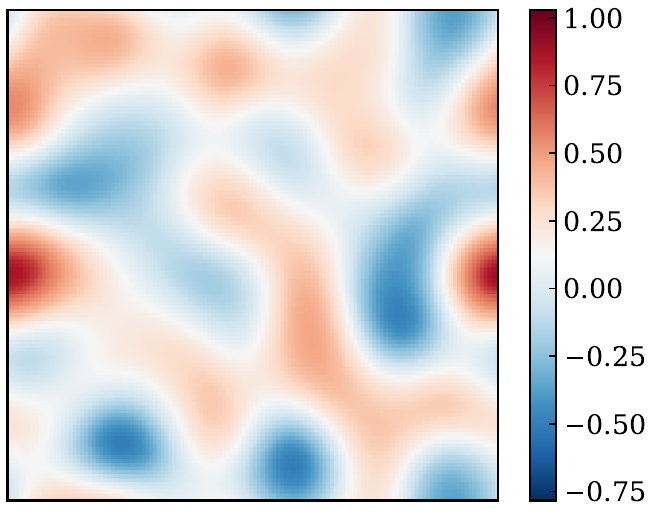} &
    \includegraphics[width=0.18\linewidth]{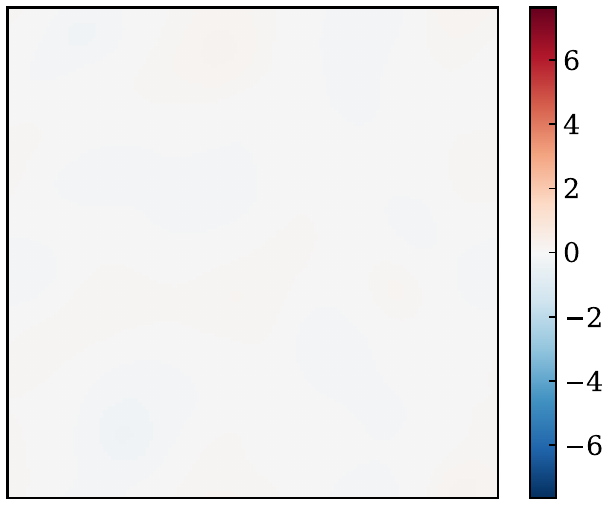}
  \end{tabular}
  \caption{Exact and predicted fields for the held-out manufactured trajectory at $t=2000$, $10000$, and $30000$ (rows) with $N=128$.
  From left to right, the columns show the exact field, the FNO field and its pointwise error, and the bcFNO field and its pointwise error.}
  \label{fig:manufactured_snapshots}
\end{figure}

Figure~\ref{fig:manufactured_snapshots} compares the exact fields, the two predicted fields, and their pointwise errors.
At $t=2000$, both models retain the main spatial pattern.
By $t=10000$, the FNO exhibits large amplitude excursions and visible morphological distortion; by $t=30000$, these errors dominate the predicted field. 
The bcFNO suppresses the amplitude growth and remains closer to the exact trajectory, although a coherent spatial displacement is still visible at late times.
At $t=30000$, the relative field errors are approximately $4.9$ for FNO and $0.28$ for bcFNO.

\FloatBarrier

\subsection{Source-Free Cahn--Hilliard Coarsening Experiments}
\label{subsec:source_free_setup}

We next consider long-time predictions of source-free Cahn--Hilliard coarsening with the FNO, bcFNO, and spFNO.
The numerical solver uses $\delta t\in\{1,0.5,0.2,0.1\}$ at $N\in\{128,256\}$, yielding eight $(N,\delta t)$ cases.
All models predict at the fixed interval $\Delta T=20$ and are compared with reference trajectories at the same physical times.
Table~\ref{tab:fno_training_settings} summarizes the data-generation, architecture, and optimization settings.
\begin{table}[!t]
\centering
\caption{Data, architecture, and optimization settings for the
source-free CH coarsening experiments.}
\label{tab:fno_training_settings}
\small
\begin{tabular}{p{0.3\linewidth}p{0.68\linewidth}}
\toprule
Parameter & Value \\
\midrule
Training trajectories
  & 6 independent initial-condition realizations \\
Evaluation trajectory
  & 1 independent held-out realization \\
Prediction interval
  & $\Delta T=20$ \\
Internal solver step
  & $\delta t\in\{1,0.5,0.2,0.1\}$ \\
Architecture
  & Two-dimensional FNO; 4 Fourier layers; width 64; GELU activation\\
Fourier modes, $N=128$
  & $16\times16$; 2,410,177 parameters \\
Fourier modes, $N=256$
  & $48\times48$; 19,711,681 parameters \\
Optimizer
  & Adam; learning rate $10^{-3}$ \\
Batch size
  & 16 \\
Training schedule
  & 10 sampling blocks $\times$ 3 inner epochs \\
Resolution protocol
  & Independent training for each $(N,\delta t)$ \\
\bottomrule
\end{tabular}
\end{table}

The $N=128$ and $N=256$ models use different Fourier-mode truncations and have different parameter counts. 
Comparisons across resolutions therefore use independently trained models and do not constitute zero-shot resolution transfer.
At a fixed $(N,\delta t)$, FNO, bcFNO, and spFNO share the same underlying architecture and differ primarily in their training objectives. 
Training pairs with start time $t< t_{\mathrm{cut} } = 500$ are excluded.
This removes the early spinodal decomposition from the training set.
The remaining pairs are drawn from a common physical-time domain using prescribed allocations over global physical time intervals. 
In each sampling block, an allocation specifies the number of pairs selected from its interval.
Within an interval, candidate start times are chosen deterministically in increasing physical time order, giving approximately uniform coverage of that interval.
The allocations therefore set the time distribution of the one-step loss.
Each input--target pair represents the fixed-interval transition from $t$ to $t+\Delta T$. 
Pairs are used independently during one-step training.

We show the free energy $E(t)$, the maximum field magnitude $\|\phi(\cdot,t)\|_\infty$, the characteristic wavenumber $k_1(t)$, the cubic characteristic length $L^3(t)$, and the normalized structure function $s(k,t)$. 
The spatial mean $\overline{\phi}(t)$ is also monitored as deviations from mass conservation.
Coarsening statistics are computed on a contiguous reference-defined coarsening interval. 
This interval is identified once from the held-out numerical-reference trajectory using the characteristic length and its local growth rate; the selection rule is given in~\ref{app:regime_classification}. 
The same interval is then used for the numerical reference, bcFNO, and spFNO predictions. 
The fitted values are therefore late-stage coarsening comparisons with the corresponding numerical-reference trajectory.
Detailed results are presented for the representative cases $(N,\delta t)=(128,0.2)$ and $(256,0.1)$;  the remaining cases are reported in~\ref{app:robustness}.

\subsection{Sampling Sensitivity in FNO}\label{subsec:fno_sensitivity}

To investigate the effect of temporal sampling on long-time prediction, we construct two temporal sampling schemes, denoted Sampling~A and Sampling~B, by assigning different relative sampling weights to prescribed physical-time windows. These weights determine the relative contribution of one-step training pairs from each window to the empirical training loss.
Sampling~A and Sampling~B are defined separately for each resolution and do not refer to a common pair of window weights across resolutions.
Table~\ref{tab:sampling_schedules} lists the normalized window weights for the two representative cases. 
Within each case, the training trajectories, network architecture, optimization procedure, and coarse prediction interval are fixed.

\begin{figure}[!t]
  \centering
  \setlength{\tabcolsep}{8pt}
  \renewcommand{\arraystretch}{1.08}
  \begin{tabular}{@{}c@{\hspace{0.06\textwidth}}c@{}}
    \begin{subfigure}[t]{0.46\textwidth}
      \centering
      \includegraphics[width=\linewidth]{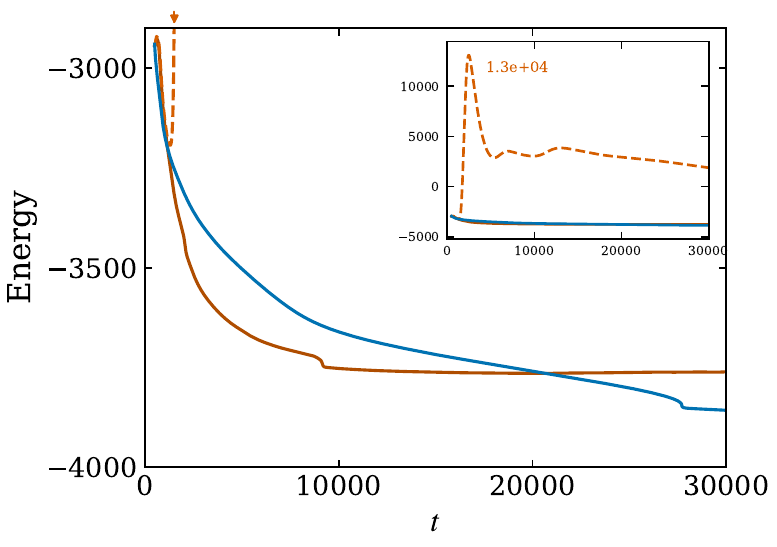}
    \end{subfigure} &
    \begin{subfigure}[t]{0.46\textwidth}
      \centering
      \includegraphics[width=\linewidth]{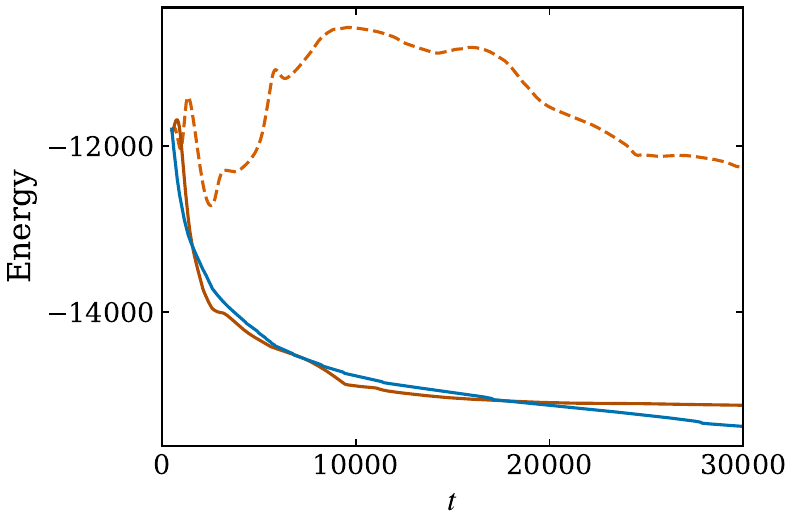}
    \end{subfigure} \\[0.5em]
    \begin{subfigure}[t]{0.46\textwidth}
      \centering
      \includegraphics[width=\linewidth]{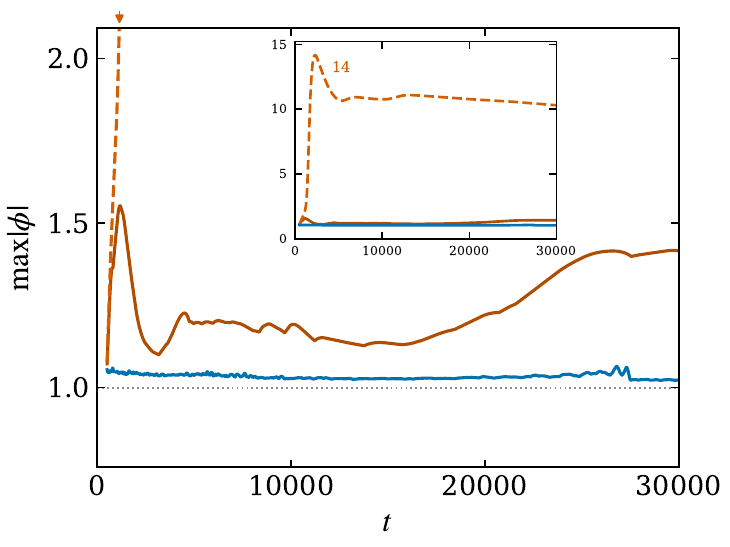}
    \end{subfigure} &
    \begin{subfigure}[t]{0.46\textwidth}
      \centering
      \includegraphics[width=\linewidth]{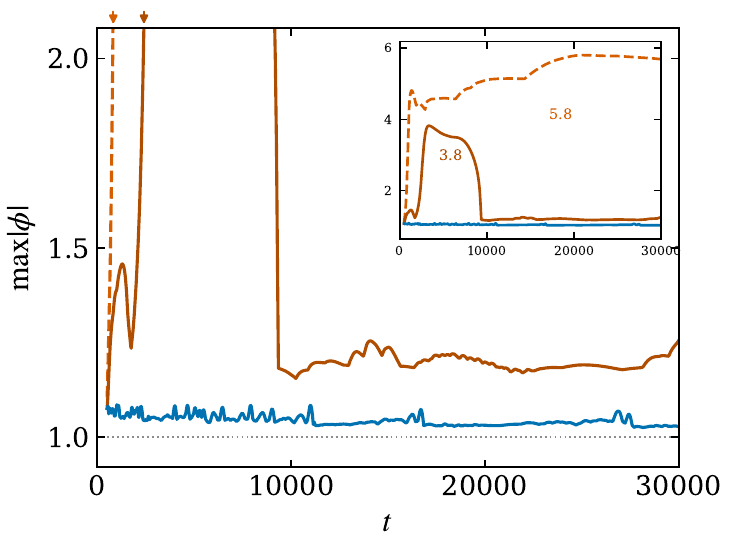}
    \end{subfigure}
  \end{tabular}
    \caption{Sensitivity of long-time FNO rollouts to physical-time sampling. Free energy $E(t)$ (top) and maximum field magnitude $\max|\phi|$ (bottom) are shown for
    $(N,\delta t)=(128,0.2)$ (left) and $(256,0.1)$ (right).
    Blue and orange denote the reference trajectories and FNO predictions, respectively; solid and dashed lines correspond to Sampling~A and Sampling~B.}
  \label{fig:fno_sampling_sensitivity}
\end{figure}

Figure~\ref{fig:fno_sampling_sensitivity} compares the resulting long-time trajectories.
For both representative cases, the FNO trained with Sampling~A remains bounded over the prediction time. In contrast, Sampling~B produces large excursions beyond the nominal phase-field range, followed by departures from the corresponding reference trajectories.
The free-energy trajectories exhibit the same qualitative contrast.

\begin{table}[!t]
\centering
\caption{Normalized physical-time bucket probabilities for the temporal sampling distributions used in the two representative cases.}
\label{tab:sampling_schedules}
\small
\begin{tabular}{llcccc}
\toprule
$(N,\delta t)$
& Schedule
& $[500,2000]$
& $[2000,10000]$
& $[10000,20000]$ 
& $[20000,30000]$ \\
\midrule
$(128,0.2)$
& Sampling~A & 15\% & 25\% & 35\% & 25\% \\
& Sampling~B & 20\% & 20\% & 35\% & 25\% \\
\addlinespace
$(256,0.1)$
& Sampling~A & 30\% & 20\% & 30\% & 20\% \\
& Sampling~B & 33.3\% & 22.2\% & 27.8\% & 16.7\% \\
\bottomrule
\end{tabular}
\end{table}

These results indicate that the long-time prediction of the FNO is sensitive to the temporal sampling of the one-step training data.
Because the sampling distributions are case-specific, the result does not compare a fixed set of weights across resolutions.
Instead, the transition from bounded to unstable predictions shows that one-step fitting and temporal redistribution of training samples alone do not reliably control long-time amplitude growth.
This observation motivates the bound-conforming formulation introduced next.

\subsection{Effect of the Bound-Conforming Penalty on Long-Time Prediction}\label{subsec:exp_bd}

\begin{figure}[ht]
  \centering
  \begin{subfigure}[t]{0.24\textwidth}
    \centering
    \includegraphics[width=\linewidth]{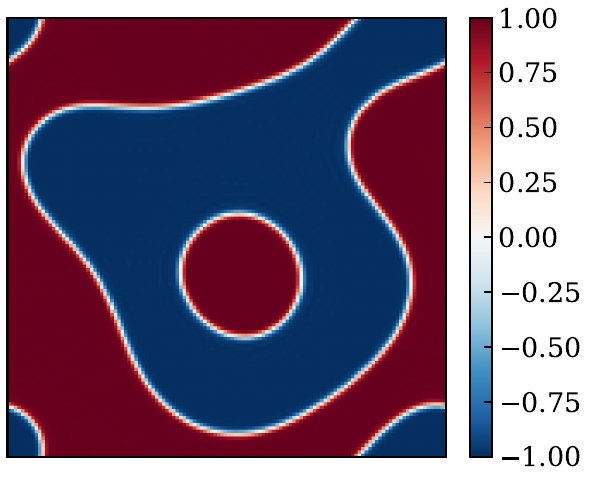}
  \end{subfigure}%
  \hspace{1.5em}%
  \begin{subfigure}[t]{0.23\textwidth}
    \centering
    \includegraphics[width=\linewidth]{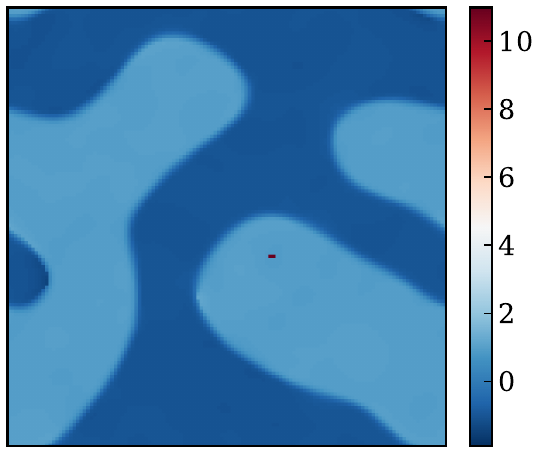}
  \end{subfigure}%
  \hspace{1.5em}%
  \begin{subfigure}[t]{0.24\textwidth}
    \centering
    \includegraphics[width=\linewidth]{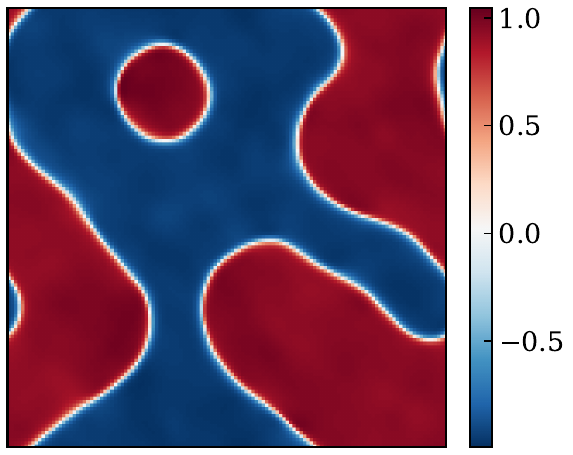}
  \end{subfigure}
  \par\vspace{0.1em}
  \begin{subfigure}[t]{0.24\textwidth}
    \centering
    \includegraphics[width=\linewidth]{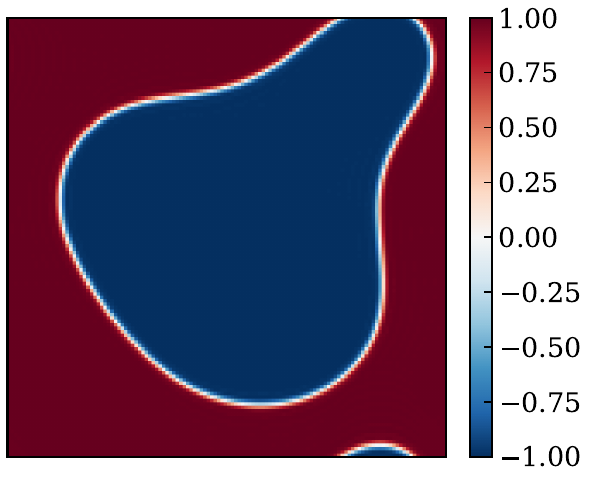}
  \end{subfigure}%
  \hspace{1.5em}%
  \begin{subfigure}[t]{0.23\textwidth}
    \centering
    \includegraphics[width=\linewidth]{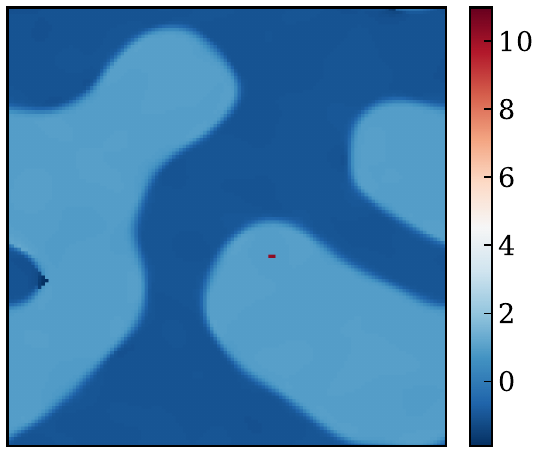}
  \end{subfigure}%
  \hspace{1.5em}%
  \begin{subfigure}[t]{0.24\textwidth}
    \centering
    \includegraphics[width=\linewidth]{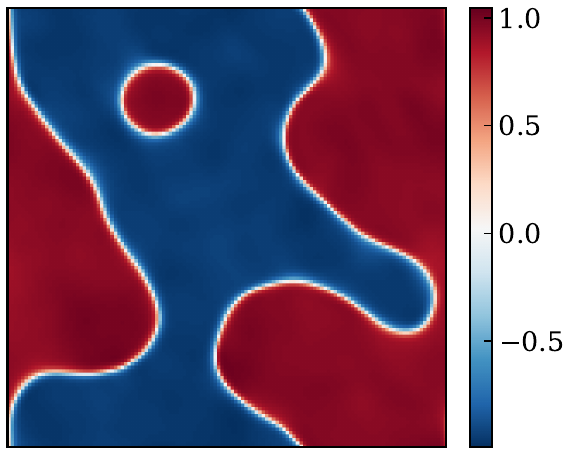}
  \end{subfigure}
  \caption{Spatial phase-field trajectories at $t=15000$ (top row) and $t=29000$ (bottom row) for the reference trajectory (left), FNO (middle), and bcFNO (right). $N=128$, $\delta t=0.2$.}
  \label{fig:bound_morphology}
\end{figure}

We examine the effect of the bound-conforming penalty on
long-time prediction. 
For $(N,\delta t)=(128,0.2)$, FNO and bcFNO use the same architecture, training data, temporal sampling strategy, and optimization budget, and differ only in the inclusion of the bound-conforming penalty. 
The case-specific penalty weights are listed in Table~\ref{tab:training_parameters}.
Figure~\ref{fig:bound_morphology} shows the numerical reference, FNO, and bcFNO phase fields at $t=15000$ and $t=29000$, from the same initial condition. 
The displayed fields are used to show phase-field amplitudes and domain morphology; they are not intended as a pointwise comparison with the reference.
The FNO develops increasingly large departures from the nominal phase-field range, accompanied by distorted domain patterns. 
The bcFNO keeps these amplitude departures smaller and retains a phase-separated domain pattern over the displayed late-stage interval.

\begin{figure}[!t]
  \centering
  \begin{subfigure}[t]{0.48\textwidth}
    \centering
    \includegraphics[width=\linewidth]{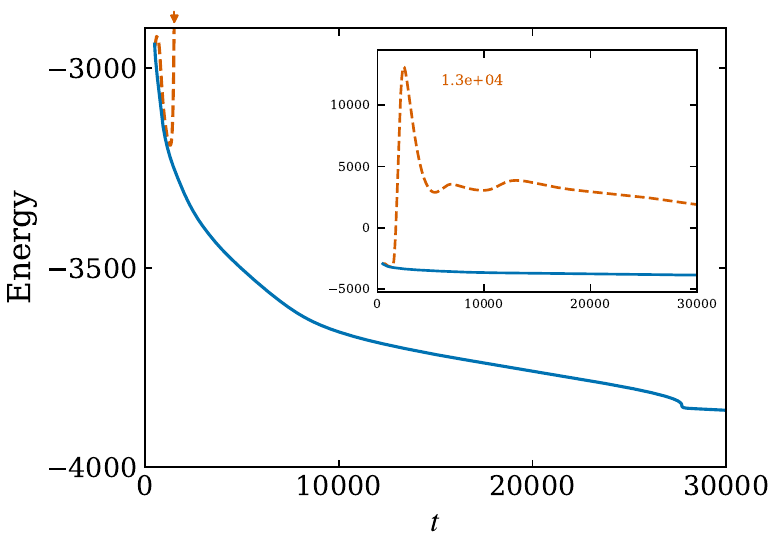}
  \end{subfigure}
  \hfill
  \begin{subfigure}[t]{0.48\textwidth}
    \centering
    \includegraphics[width=\linewidth]{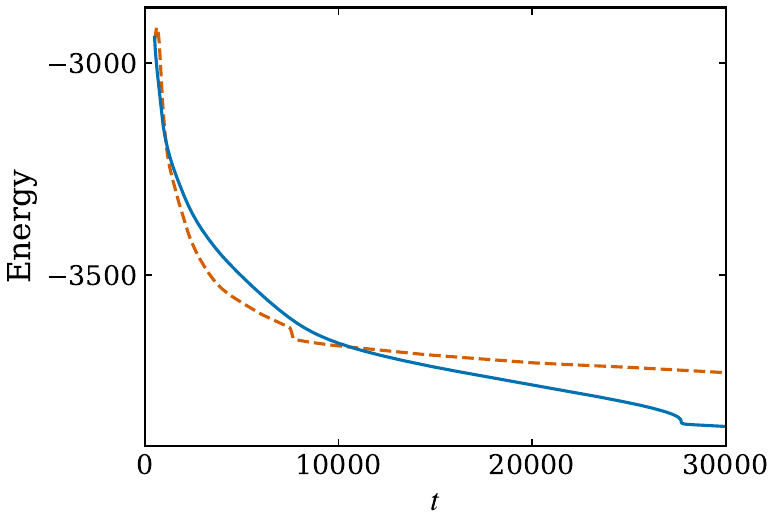}
  \end{subfigure}
  \par\vspace{0.6em}
  \begin{subfigure}[t]{0.48\textwidth}
    \centering
    \includegraphics[width=\linewidth]{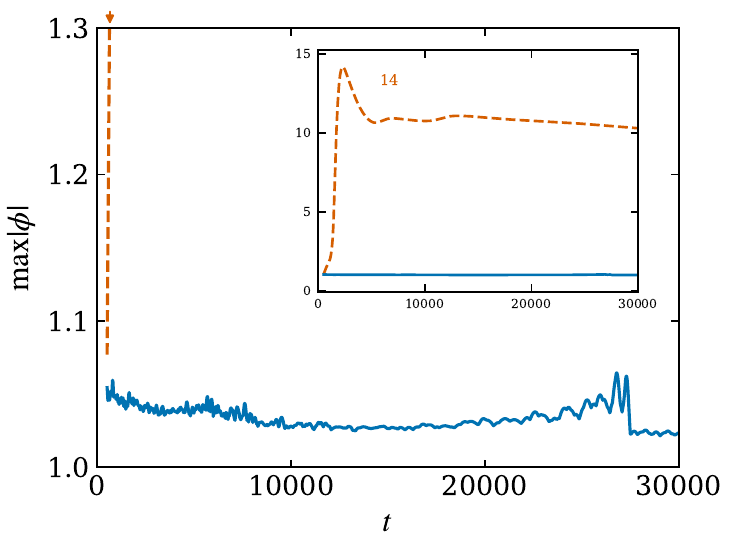}
  \end{subfigure}
  \hfill
  \begin{subfigure}[t]{0.48\textwidth}
    \centering
    \includegraphics[width=\linewidth]{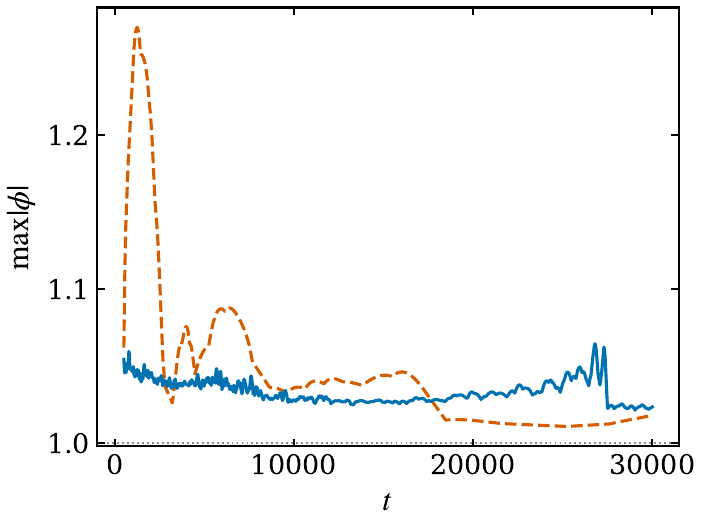}
  \end{subfigure}
     \caption{Evolution of the free energy $E(t)$ (top row) and $\max|\phi|$ (bottom row) for FNO (left) and bcFNO (right). Blue and orange denote the reference trajectory and model prediction,
  respectively. $N=128$, $\delta t=0.2$.}\label{fig:bd_prediction}
\end{figure}

Figure~\ref{fig:bd_prediction} compares the free energy $E(t)$ and
maximum field magnitude $\|\phi(\cdot,t)\|_\infty$ along an independent evaluation trajectory.
The FNO exhibits a pronounced late-time increase in free energy together with persistent departures from the nominal phase-field range.
The bcFNO suppresses the large amplitude excursions and
produces a free energy evolution that remains closer to the reference trajectory.
Additional results for remaining cases are provided in~\ref{app:fno_failures}
(Figs.~\ref{fig:app_fno_failure_energy_N128},~\ref{fig:app_fno_failure_energy_N256},~\ref{fig:app_fno_phi_max_N128}, and~\ref{fig:app_fno_phi_max_N256}). 
The multi-seed results in~\ref{app:multiseed} further show that large FNO amplitude excursions recur across the tested training seeds, whereas bcFNO suppresses such excursions under the same training procedure.

The connection between amplitude control and the subsequent spectral analysis follows from Eqs.~\eqref{eq:variance}--\eqref{eq:s_discrete}.
Since the structure function is computed from the fluctuation field $\psi=\phi-\bar\phi$, the zero Fourier mode of $\psi$ vanishes.
With the discrete Fourier-transform normalization used here, Parseval's identity gives
\begin{equation}
\frac{1}{N^4}
\sum_{\mathbf k\ne\mathbf 0}
\left|\widehat{\psi}(\mathbf k,t)\right|^2
=\frac{1}{N^2}
\sum_{\mathbf r}
\left|\psi(\mathbf r,t)\right|^2=
\frac{1}{N^2}
\sum_{\mathbf r}
|\phi(\mathbf r,t)-\bar\phi|^2
\leq
\|\phi\|_{\infty}^2-|\bar\phi|^2.
\label{eq:bd_spectral_bound}
\end{equation}
When the spatial mean remains nearly constant, suppression of large phase-field amplitude  reduces the upper bound on the total power in the nonzero Fourier modes.
This relation provides a spectral interpretation of the amplitude control provided by bcFNO.

Equation~\eqref{eq:bd_spectral_bound} constrains only the aggregate nonzero-mode power. 
It does not specify its distribution across wavenumbers and therefore cannot determine the characteristic wavenumber $k_1(t)$ or other coarsening statistics.
Those statistics depend on the normalized structure function. 
The spFNO regularization introduced next acts on this spectral distribution after the large amplitude excursions have been controlled.

\FloatBarrier

\subsection{Structure-Function Correction of Coarsening Dynamics}\label{subsec:exp_sk}

After establishing that bcFNO suppresses the large amplitude excursions of FNO, we consider the effect of adding structure-function regularization to the stabilized bcFNO.
The spFNO acts directly on the normalized structure function  $s(k,t)$.
We first compare the correction of $s(k,t)$ and then consider the coarsening statistics derived from it.

\begin{figure}[!t]
  \centering
  \begin{minipage}{0.24\textwidth}
    \centering
    \includegraphics[width=\linewidth]{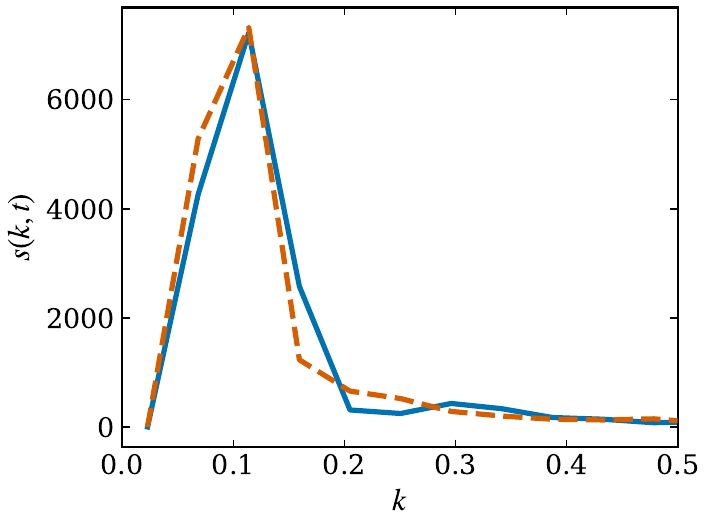}
  \end{minipage}
  \hfill
  \begin{minipage}{0.24\textwidth}
    \centering
    \includegraphics[width=\linewidth]{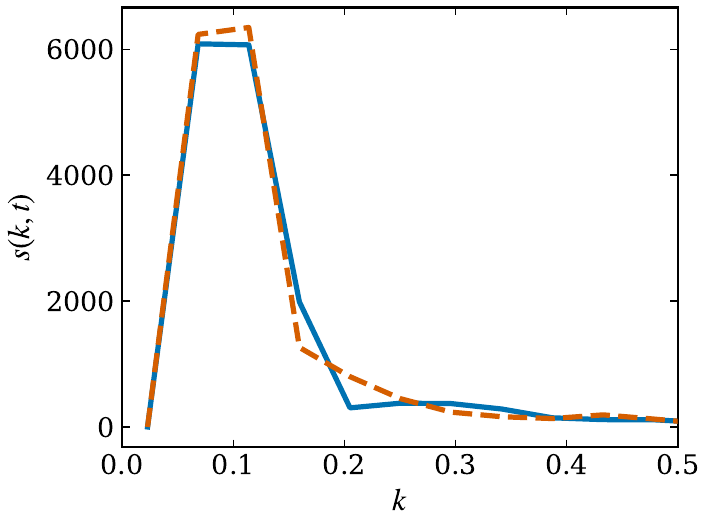}
  \end{minipage}
  \hfill
  \begin{minipage}{0.24\textwidth}
    \centering
    \includegraphics[width=\linewidth]{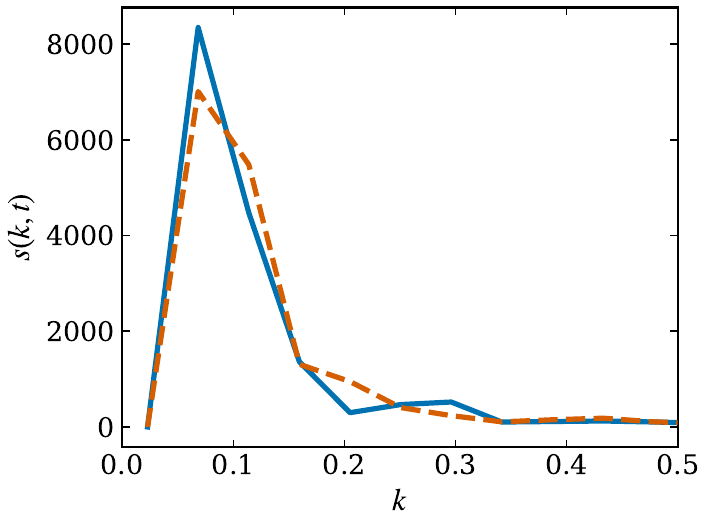}
  \end{minipage}
  \hfill
  \begin{minipage}{0.24\textwidth}
    \centering
    \includegraphics[width=\linewidth]{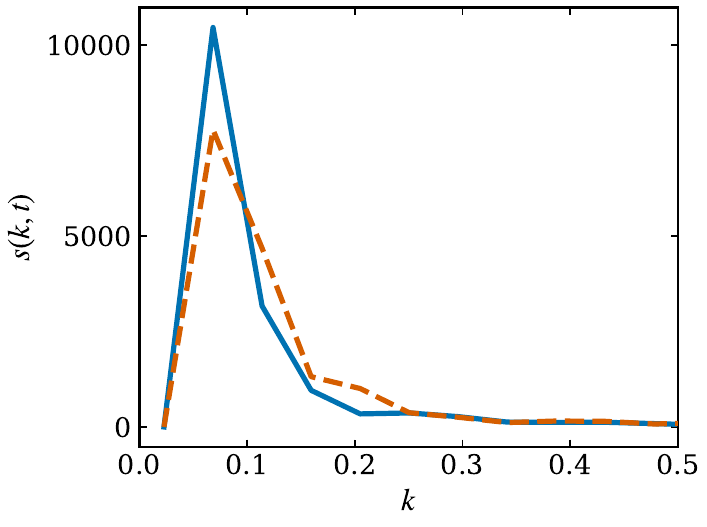}
  \end{minipage}
  \par\vspace{0.5em}
  \begin{minipage}{0.24\textwidth}
    \centering
    \includegraphics[width=\linewidth]{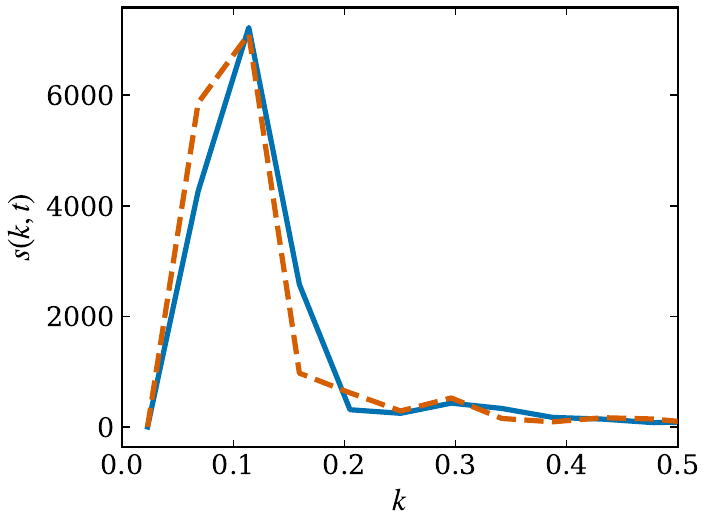}
  \end{minipage}
  \hfill
  \begin{minipage}{0.24\textwidth}
    \centering
    \includegraphics[width=\linewidth]{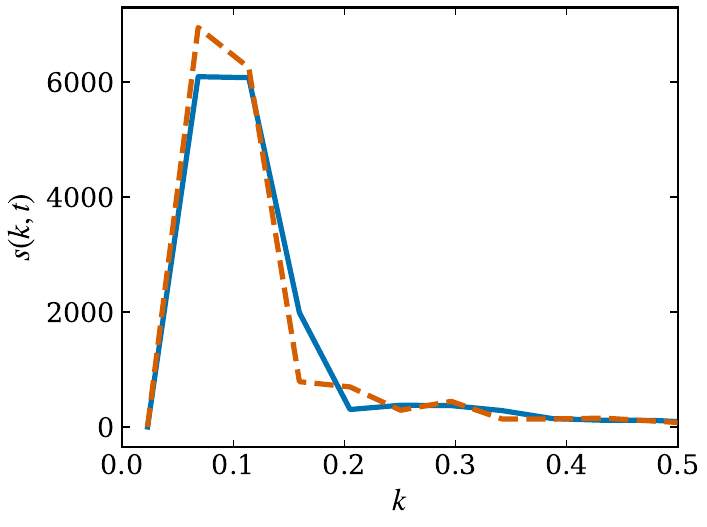}
  \end{minipage}
  \hfill
  \begin{minipage}{0.24\textwidth}
    \centering
    \includegraphics[width=\linewidth]{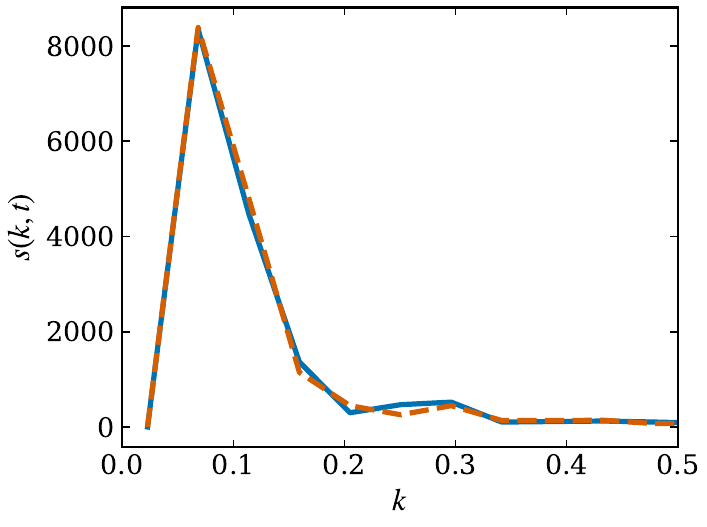}
  \end{minipage}
  \hfill
  \begin{minipage}{0.24\textwidth}
    \centering
    \includegraphics[width=\linewidth]{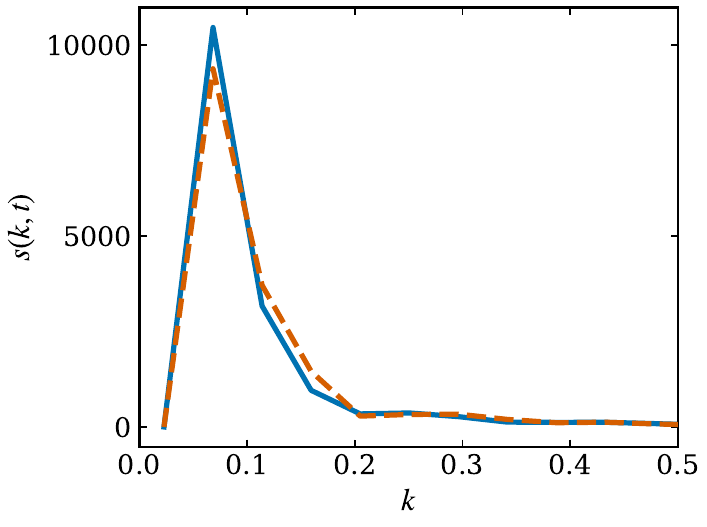}
  \end{minipage}
  \caption{Late-stage normalized structure functions $s(k,t)$ for bcFNO (top) and spFNO (bottom) at $t=17000$, $21000$, $25000$, and $29000$ from left to right.
  Blue and orange denote the reference trajectory and model prediction,
  respectively. $N=128$, $\delta t=0.2$.}
  \label{fig:sk_recovery_N128}
\end{figure}

\begin{figure}[!t]
  \centering
  \begin{minipage}{0.24\textwidth}
    \centering
    \includegraphics[width=\linewidth]{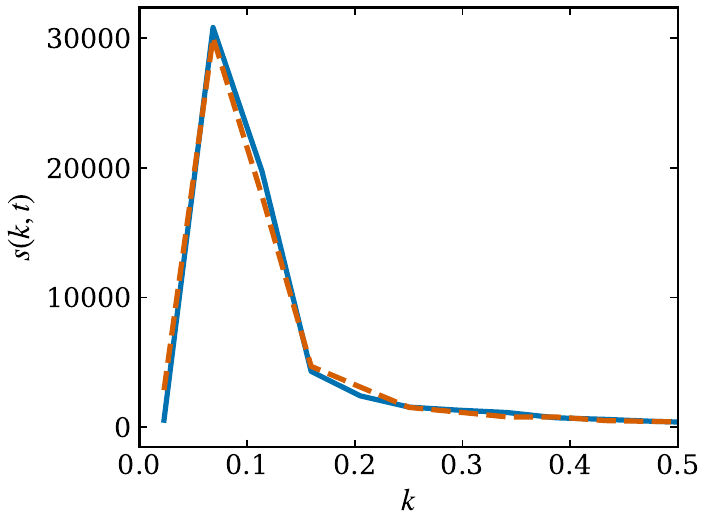}
  \end{minipage}
  \hfill
  \begin{minipage}{0.24\textwidth}
    \centering
    \includegraphics[width=\linewidth]{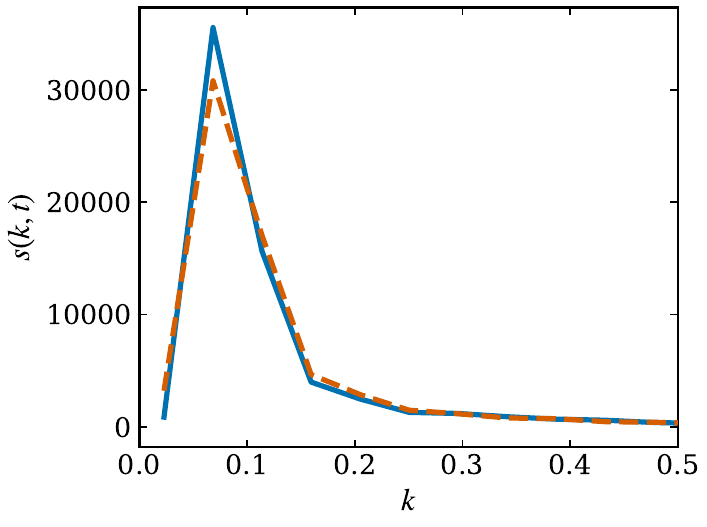}
  \end{minipage}
  \hfill
  \begin{minipage}{0.24\textwidth}
    \centering
    \includegraphics[width=\linewidth]{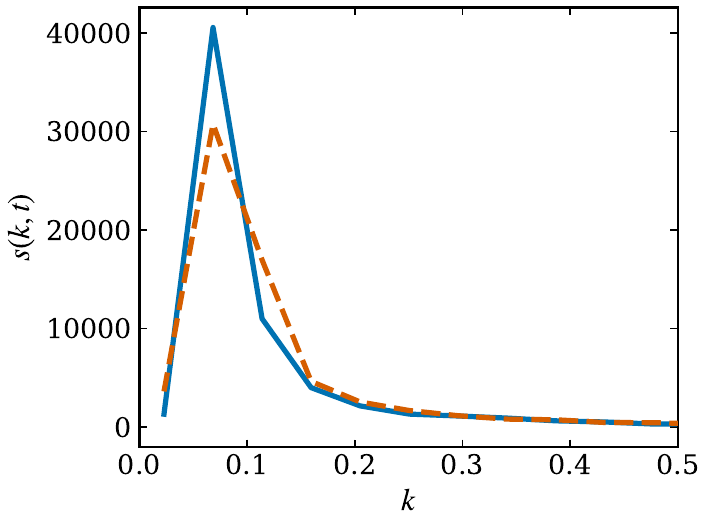}
  \end{minipage}
  \hfill
  \begin{minipage}{0.24\textwidth}
    \centering
    \includegraphics[width=\linewidth]{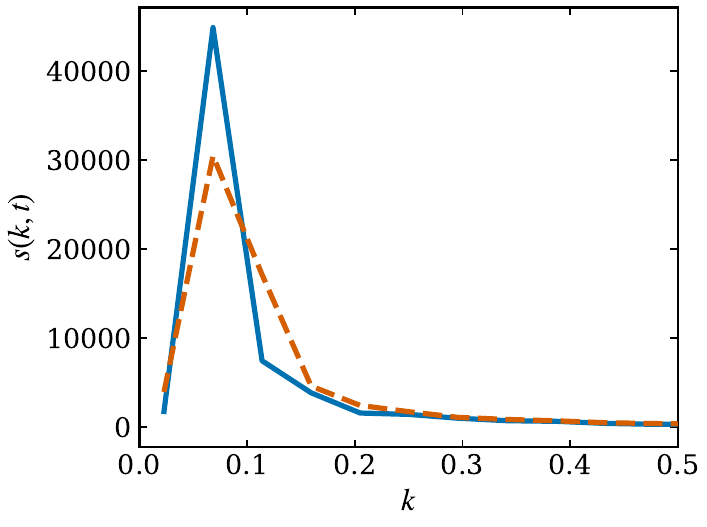}
  \end{minipage}
  \par\vspace{0.5em}
  \begin{minipage}{0.24\textwidth}
    \centering
    \includegraphics[width=\linewidth]{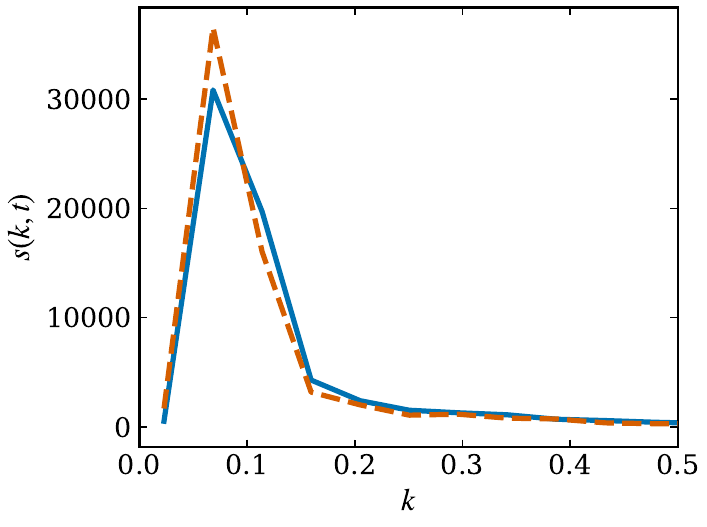}
  \end{minipage}
  \hfill
  \begin{minipage}{0.24\textwidth}
    \centering
    \includegraphics[width=\linewidth]{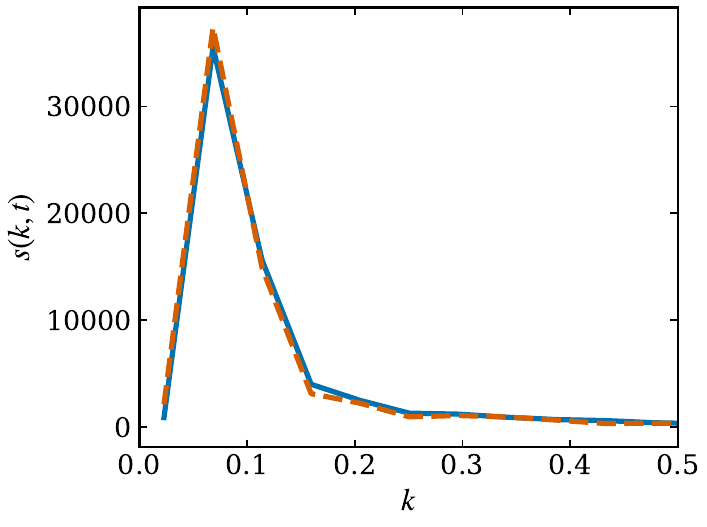}
  \end{minipage}
  \hfill
  \begin{minipage}{0.24\textwidth}
    \centering
    \includegraphics[width=\linewidth]{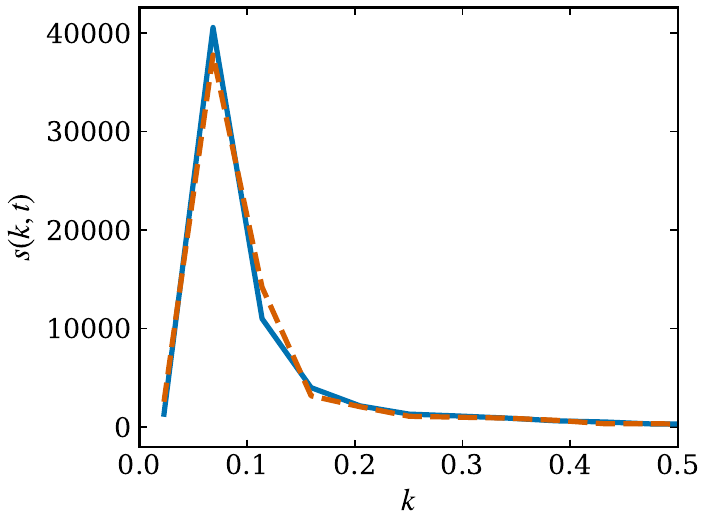}
  \end{minipage}
  \hfill
  \begin{minipage}{0.24\textwidth}
    \centering
    \includegraphics[width=\linewidth]{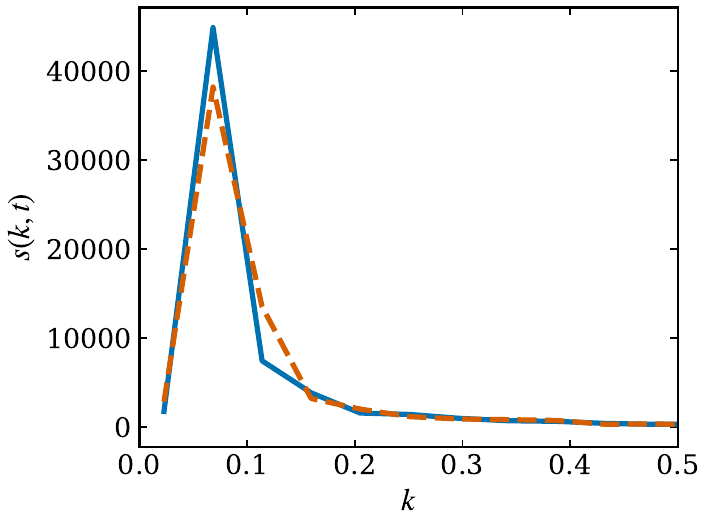}
  \end{minipage}
  \caption{Late-stage normalized structure functions $s(k,t)$ for bcFNO (top) and spFNO (bottom) at $t=17000$, $21000$, $25000$, and $29000$ from left to right. Blue and orange denote the  reference trajectory and model prediction, respectively. $N=256$, $\delta t=0.1$.}
  \label{fig:sk_recovery_N256}
\end{figure}

\begin{figure}[!t]
  \centering
  \begin{minipage}{0.24\textwidth}
    \centering
    \includegraphics[width=\linewidth]{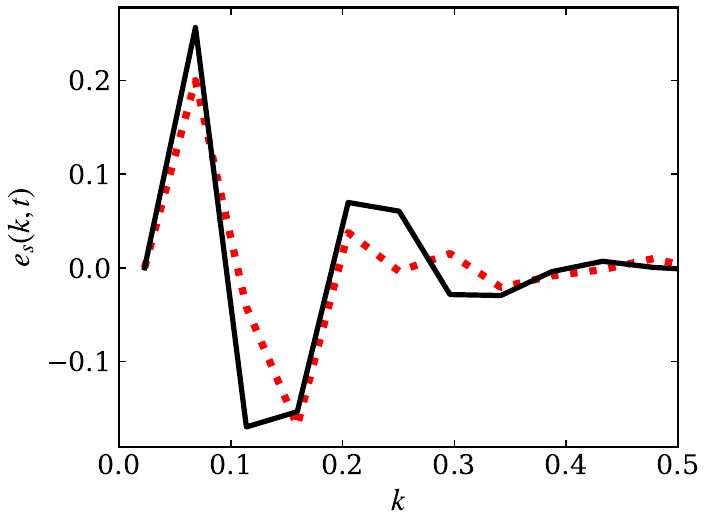}
  \end{minipage}
  \hfill
  \begin{minipage}{0.24\textwidth}
    \centering
    \includegraphics[width=\linewidth]{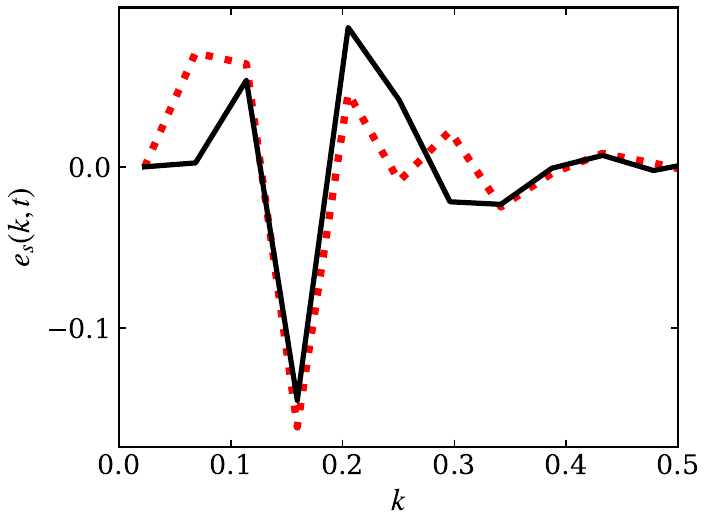}
  \end{minipage}
  \hfill
  \begin{minipage}{0.24\textwidth}
    \centering
    \includegraphics[width=\linewidth]{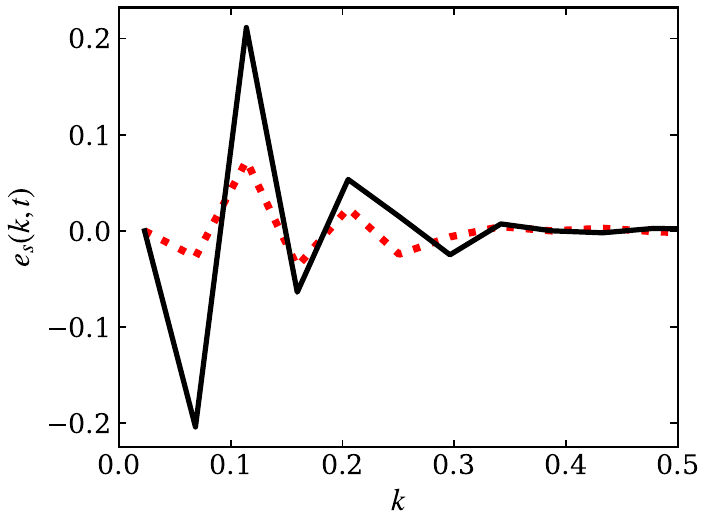}
  \end{minipage}
  \hfill
  \begin{minipage}{0.24\textwidth}
    \centering
    \includegraphics[width=\linewidth]{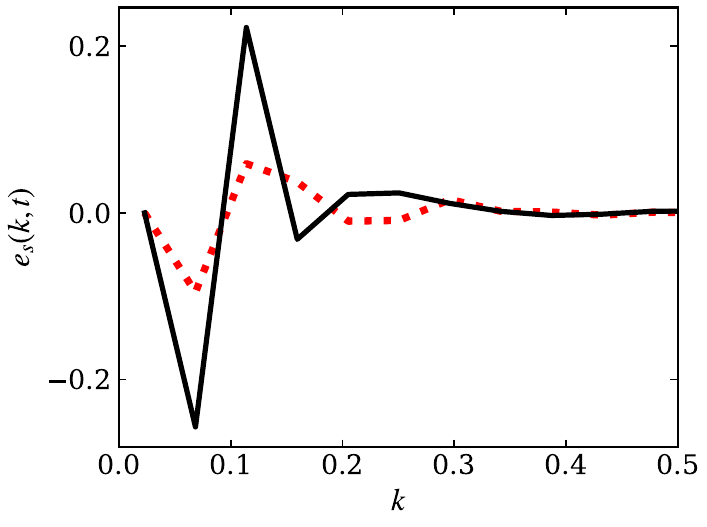}
  \end{minipage}
  \caption{Peak-normalized structure-function discrepancy $e_s(k,t)$ for
  bcFNO (solid) and spFNO (dashed) at $t=17000$, $21000$, $25000$, and
  $29000$ from left to right.
  $N=128$, $\delta t=0.2$.}
  \label{fig:sk_recovery_error_N128}
\end{figure}

\begin{figure}[!t]
  \centering
  \begin{minipage}{0.24\textwidth}
    \centering
    \includegraphics[width=\linewidth]{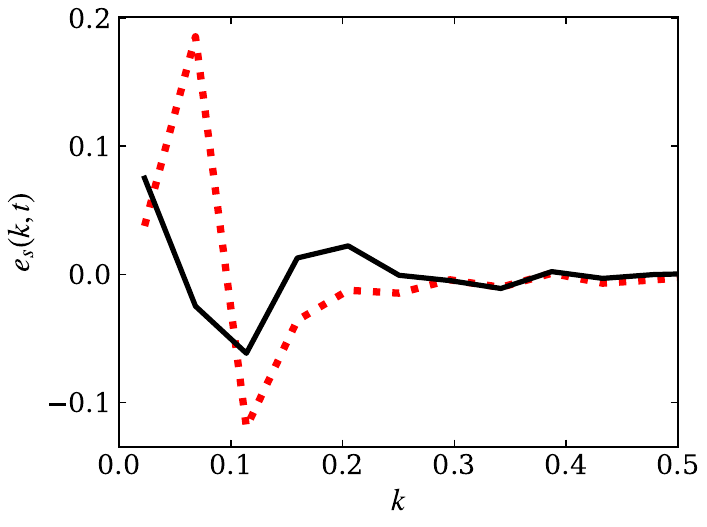}
  \end{minipage}
  \hfill
  \begin{minipage}{0.24\textwidth}
    \centering
    \includegraphics[width=\linewidth]{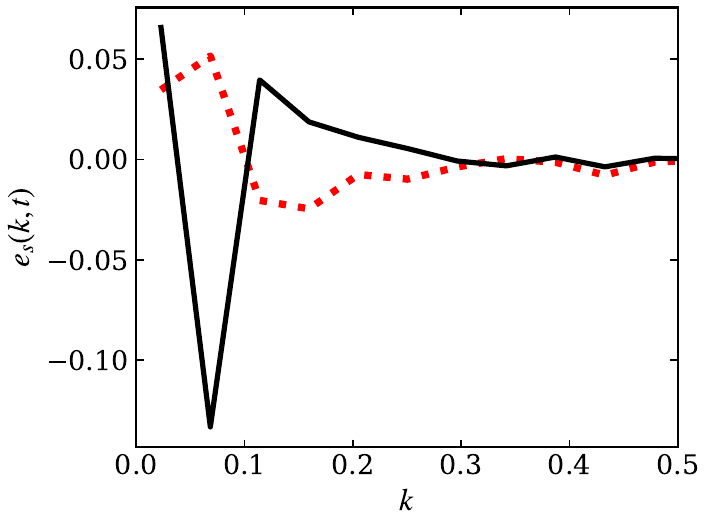}
  \end{minipage}
  \hfill
  \begin{minipage}{0.24\textwidth}
    \centering
    \includegraphics[width=\linewidth]{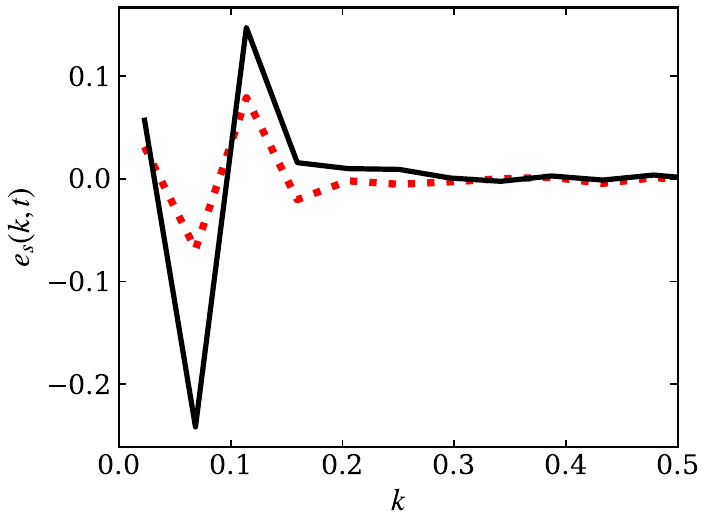}
  \end{minipage}
  \hfill
  \begin{minipage}{0.24\textwidth}
    \centering
    \includegraphics[width=\linewidth]{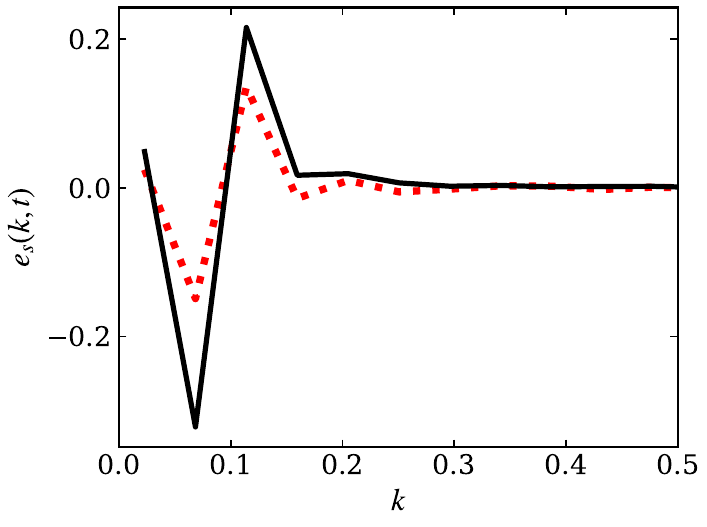}
  \end{minipage}
  \caption{Peak-normalized structure-function discrepancy $e_s(k,t)$ for bcFNO (solid) and spFNO (dashed) at $t=17000$, $21000$, $25000$, and $29000$ from left to right.
  $N=256$, $\delta t=0.1$.}
  \label{fig:sk_recovery_error_N256}
\end{figure}

Figures~\ref{fig:sk_recovery_N128} and~\ref{fig:sk_recovery_N256} compare the predicted and reference normalized structure functions at four late-stage times. 
For $(N,\delta t)=(128,0.2)$ and $(256,0.1)$, the largest differences between bcFNO and the reference occur near the low-wavenumber peak. This peak corresponds to the evolving domain scale.
At the displayed times, spFNO reduces these differences near the peak and in its neighboring wavenumber range.
To display the differences on a common scale, we use the peak-normalized signed discrepancy
\begin{equation}
e_s(k,t)
=
\frac{
s^{\mathrm{pred}}(k,t)-s^{\mathrm{ref}}(k,t)
}{
\|s^{\mathrm{ref}}(\cdot,t)\|_\infty+\epsilon_s
},
\qquad
\epsilon_s=10^{-12}.
\label{eq:sk_scaled_discrepancy}
\end{equation}
Figures~\ref{fig:sk_recovery_error_N128} and~\ref{fig:sk_recovery_error_N256} 
show that the bcFNO discrepancies are concentrated around the coarsening peak. In the same range, the spFNO discrepancies are smaller.

We next examine coarsening quantities calculated from $s(k,t)$. Over the 
coarsening interval, we fit
\[
L^3(t)= at+ b,
\qquad
k_1(t)\sim t^\beta.
\]
Here, $a$ is the fitted slope of $L^3(t)$ and $\beta$ is the fitted exponent of the first spectral moment. 
Each predicted fit is compared with the fit from the reference trajectory over the same interval.
For the rescaled spectra, we use
\begin{equation}
\mathcal{R}_{\mathrm{collapse}}
=
1-
\frac{
\left\|
\overline{\mathcal S}_{c}^{\,\mathrm{pred}}
-\overline{\mathcal S}_{c}^{\,\mathrm{ref}}
\right\|_2
}{
\left\|
\overline{\mathcal S}_{c}^{\,\mathrm{ref}}
\right\|_2
},
\label{eq:collapse_score}
\end{equation}
where $\overline{\mathcal S}_{c}(q)$ is the mean rescaled spectrum over the late-stage frames. The score equals $1$ when the two mean curves are identical.

\begin{table}[!t]
\centering
\caption{
Late-stage coarsening results for the two representative cases.
Fits use the coarsening interval identified from the numerical reference.
}
\label{tab:late_coarsening_recovery}
\small
\begin{tabular}{llrrr}
\toprule
$(N,\delta t)$ & Metric & Reference & bcFNO & spFNO \\
\midrule
$(128,0.2)$
& $a$ in $L^3(t)=at+b$
& $0.020$ & $0.0075$ & $0.011$ \\
& $\beta$ in $k_1(t)\sim t^{\beta}$
& $-0.33$ & $-0.13$ & $-0.16$ \\
& $\mathcal R_{\mathrm{collapse}}$
& 1.00 & 0.73 & 0.88 \\
\addlinespace
$(256,0.1)$
& $a$ in $L^3(t)=at+b$
& $0.020$ & $0.0083$ & $0.012$ \\
& $\beta$ in $k_1(t)\sim t^{\beta}$
& $-0.28$ & $-0.12$ & $-0.16$ \\
& $\mathcal R_{\mathrm{collapse}}$
& 1.00 & 0.91 & 0.91 \\
\bottomrule
\end{tabular}
\end{table}

\begin{figure}[ht]
  \centering
  \begin{subfigure}[t]{0.32\textwidth}
    \centering
    \includegraphics[width=\linewidth]{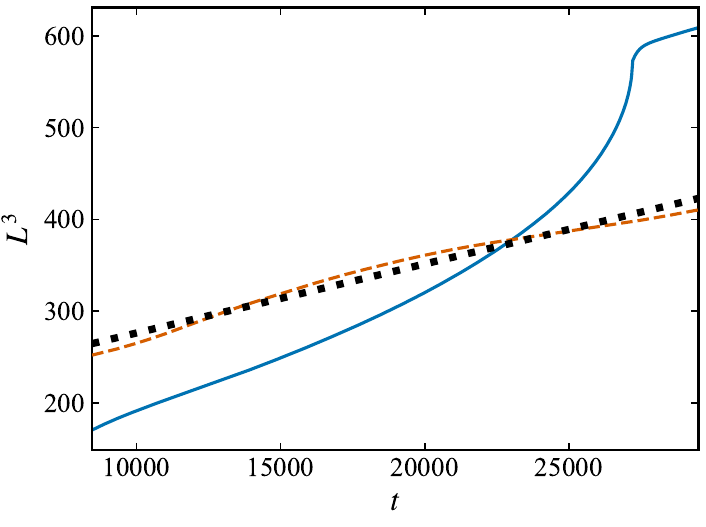}
  \end{subfigure}
  \hfill
  \begin{subfigure}[t]{0.32\textwidth}
    \centering
    \includegraphics[width=\linewidth]{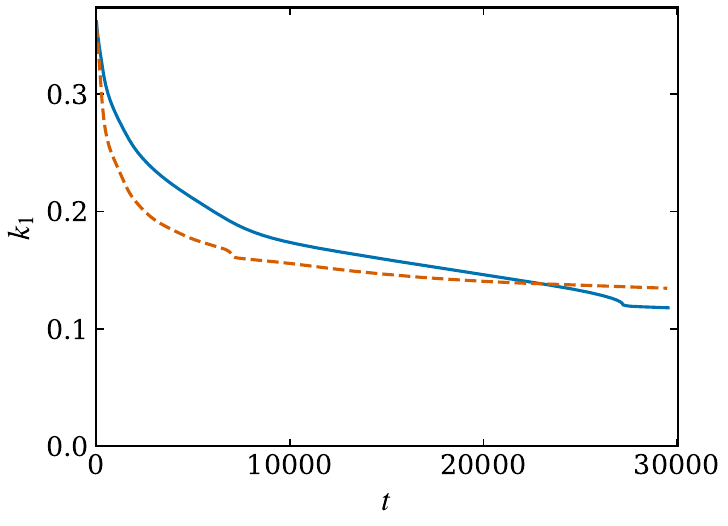}
  \end{subfigure}
  \hfill
  \begin{subfigure}[t]{0.32\textwidth}
    \centering
    \includegraphics[width=\linewidth]{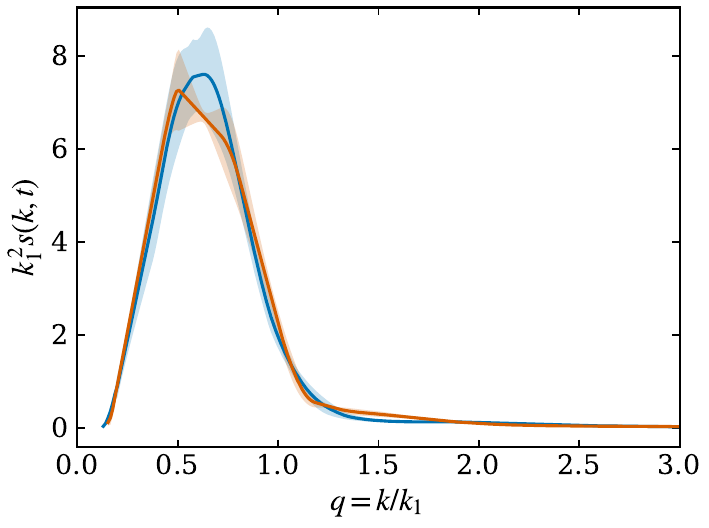}
  \end{subfigure}
  \par\vspace{0.1em}
  \begin{subfigure}[t]{0.32\textwidth}
    \centering
    \includegraphics[width=\linewidth]{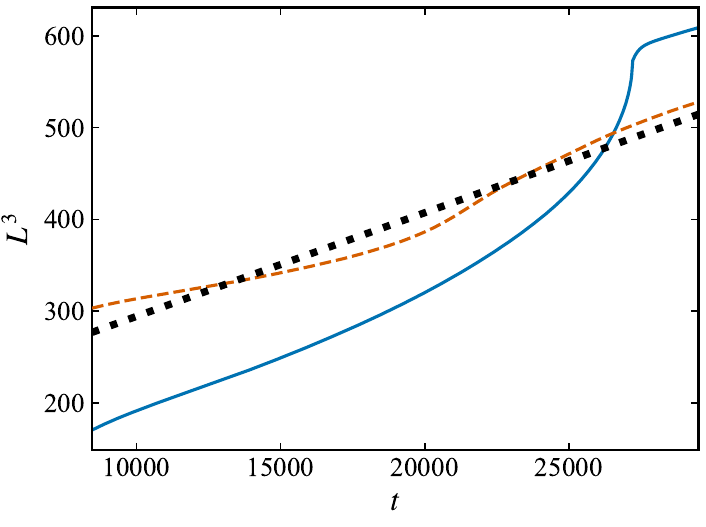}
  \end{subfigure}
  \hfill
  \begin{subfigure}[t]{0.32\textwidth}
    \centering
    \includegraphics[width=\linewidth]{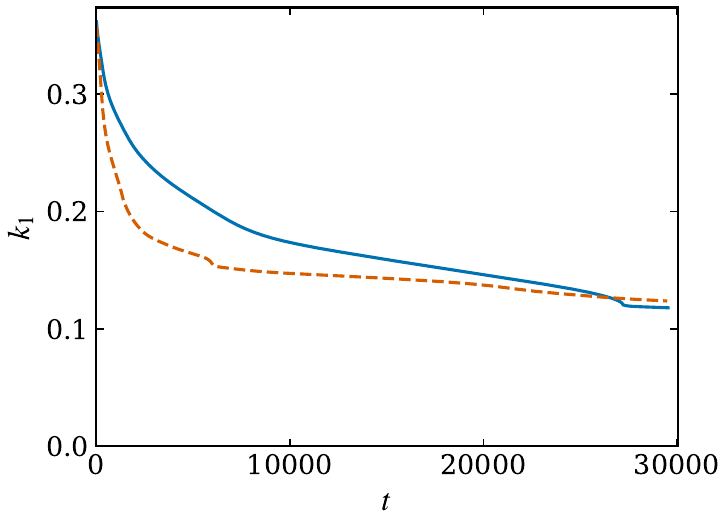}
  \end{subfigure}
  \hfill
  \begin{subfigure}[t]{0.32\textwidth}
    \centering
    \includegraphics[width=\linewidth]{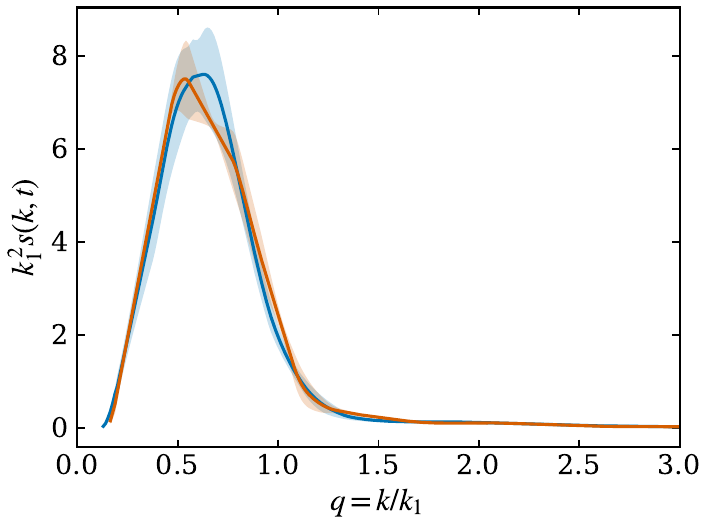}
  \end{subfigure}
\caption{Late-stage coarsening measures, $L^3(t), k_1(t)$, and the dynamic-scaling collapse, for bcFNO (top row) and spFNO (bottom row).
Blue and orange curves denote respectively the reference trajectory and model predictions.  
Black dotted lines in the $L^3(t)$ panels are linear fits. 
In the collapse panels, curves and shading denote binwise means and ranges. 
$N=128$, $\delta t=0.2$.}
  \label{fig:sk_recovery_curves}
\end{figure}

\begin{figure}[!t]
  \centering
  \begin{subfigure}[t]{0.32\textwidth}
    \centering
    \includegraphics[width=\linewidth]{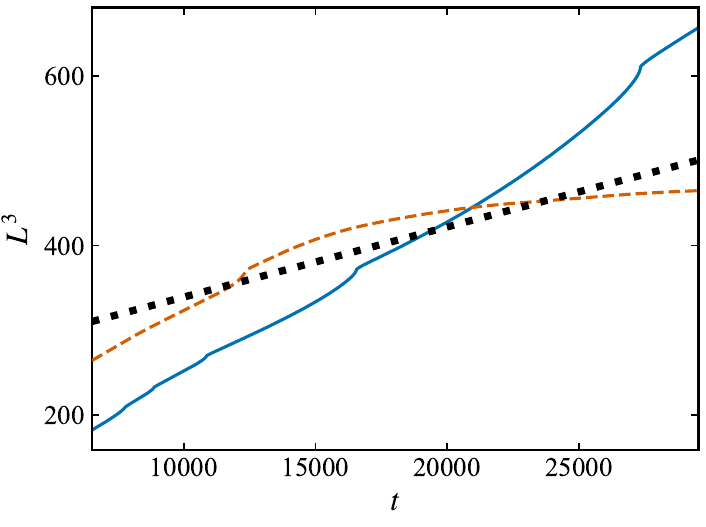}
  \end{subfigure}
  \hfill
  \begin{subfigure}[t]{0.32\textwidth}
    \centering
    \includegraphics[width=\linewidth]{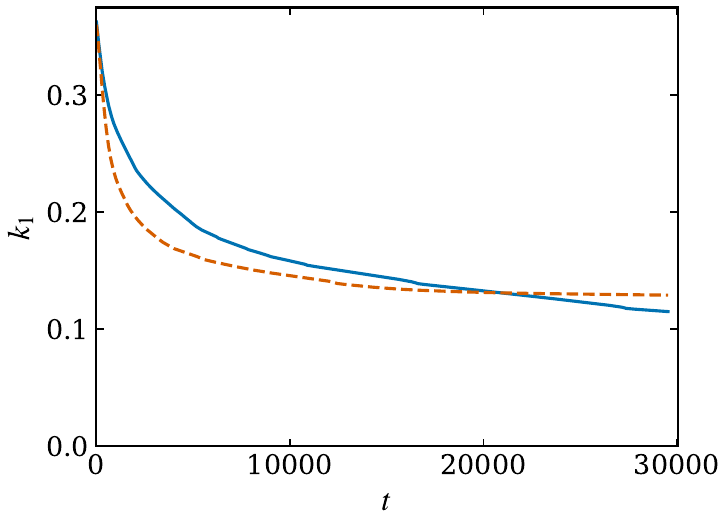}
  \end{subfigure}
  \hfill
  \begin{subfigure}[t]{0.32\textwidth}
    \centering
    \includegraphics[width=\linewidth]{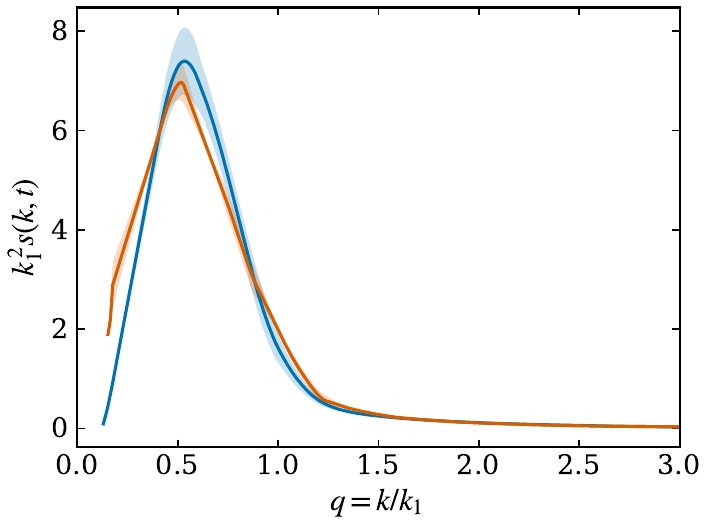}
  \end{subfigure}
  \par\vspace{0.1em}
  \begin{subfigure}[t]{0.32\textwidth}
    \centering
    \includegraphics[width=\linewidth]{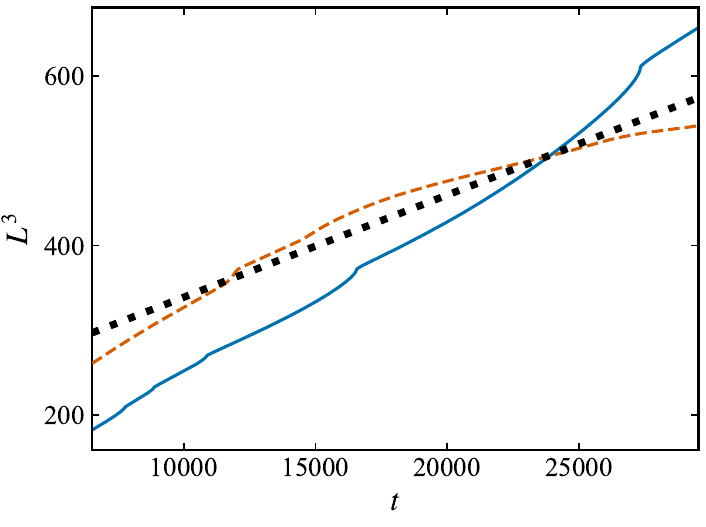}
  \end{subfigure}
  \hfill
  \begin{subfigure}[t]{0.32\textwidth}
    \centering
    \includegraphics[width=\linewidth]{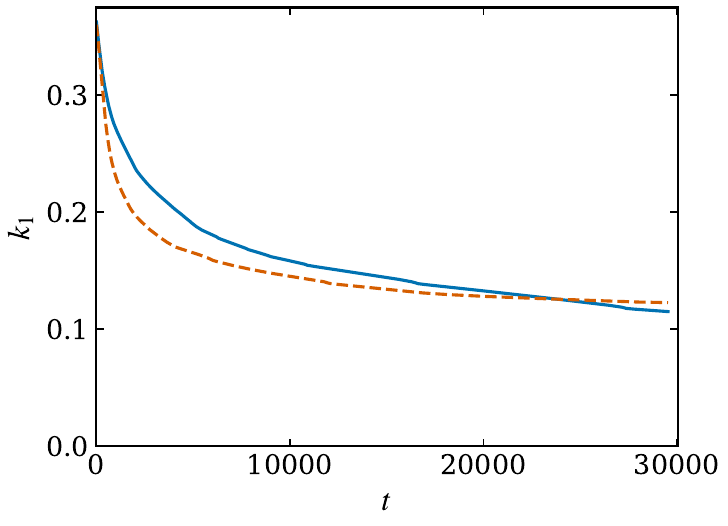}
  \end{subfigure}
  \hfill
  \begin{subfigure}[t]{0.32\textwidth}
    \centering
    \includegraphics[width=\linewidth]{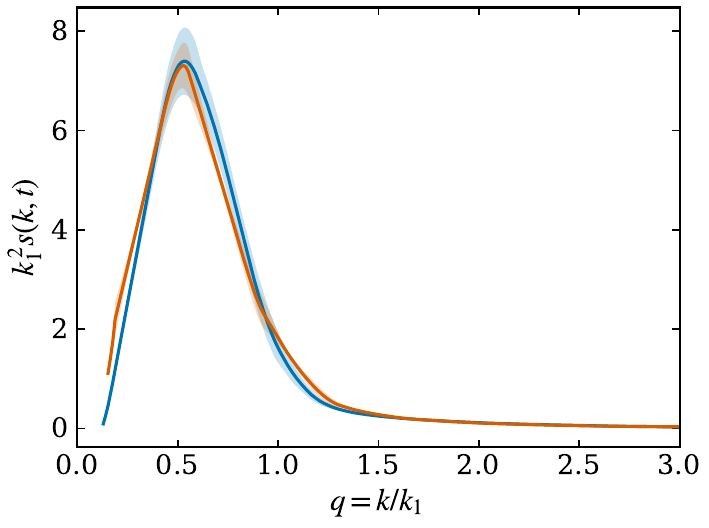}
  \end{subfigure}
  \caption{Late-stage coarsening measures, $L^3(t)$, $k_1(t)$, and the dynamic-scaling collapse, for bcFNO (top row) and spFNO (bottom row). Blue and orange curves denote the reference trajectory and model predictions, respectively. Black dotted lines in the $L^3(t)$ panels are linear fits. $N=256$, $\delta t=0.1$.}
  \label{fig:sk_recovery_curves_N256}
\end{figure}

Table~\ref{tab:late_coarsening_recovery} lists the fitted slope $a$, the wavenumber exponent $\beta$, and the collapse score $\mathcal{R}_{\mathrm{collapse}}$  for the two cases. 
For both cases, the spFNO values of $a$ and $\beta$ lie closer to the reference values than the bcFNO values. 
The collapse score increases for $(N,\delta t)=(128,0.2)$ and remains the same for $(N,\delta t)=(256,0.1)$. 
Figures~\ref{fig:sk_recovery_curves} and~\ref{fig:sk_recovery_curves_N256} show the corresponding curves of $L^3(t)$, $k_1(t)$, and the rescaled spectra.
For the collapse, the discrete rescaled spectral samples from all selected late-stage frames are pooled into 64 common $q$-bins. 
The solid curve joins the binwise means, and the shaded region gives the minimum and maximum samples in each bin. 
The width of this region includes variation across time and variation among the discrete $q$-values placed in the same bin. 
It need not vanish even when the spectra are close to a common rescaled profile, especially where that profile changes rapidly with $q$.

These results distinguish the roles of the two penalties. 
The bound-conforming term suppresses the large field-amplitude excursions seen in the FNO predictions. 
The structure-function term reduces errors near the low-wavenumber coarsening peak and changes the derived fits of $L^3(t)$ and $k_1(t)$. 
The spFNO therefore provides a spectral correction of bcFNO rather than a uniformly stronger long-time stabilization.

\FloatBarrier

\subsection{Online Cost of Trained Neural Operators}
\label{subsec:computational_cost}

We compare the online wall-clock cost of Fourier pseudospectral trajectory generation and neural-operator rollouts over the common physical-time interval $0\leq t\leq30000$.
All trajectories are stored at the same output interval $\Delta T=20$.
The recorded times exclude visualization, file output, and evaluation of derived quantities. The timing setup and warm-up treatment are detailed in~\ref{app:timing}.

The comparison isolates trajectory-generation cost after model training.
The one-time costs of generating numerical training data and optimizing the neural operators are excluded.
The results therefore characterize the computational workflow used in this study rather than the hardware-independent algorithmic speedup. 
All neural operators are evaluated from the same initial field for the representative $N=256$ cases.

\begin{table}[!t]
\centering
\caption{Online wall-clock time for trajectory generation over
$0\leq t\leq30000$ at $N=256$. Relative wall-clock time is normalized by
the FNO time.}
\label{tab:online_cost}
\small
\begin{tabular}{lccc}
\toprule
Method & Hardware & Time (s) & Relative time \\
\midrule
Solver, $\delta t=1$   & CPU & 156  & 19.5 \\
Solver, $\delta t=0.1$ & CPU & 1608 & 201  \\
FNO                    & GPU & 8.0  & 1.0  \\
bcFNO                  & GPU & 8.0  & 1.0  \\
spFNO                  & GPU & 8.0  & 1.0  \\
\bottomrule
\end{tabular}
\end{table}

Table~\ref{tab:online_cost} shows that, under the present computational workflow, each neural operator generates the trajectory over the full time interval in approximately $8.0$~s. 
The bound-conforming and structure-function terms modify the training objective only. Accordingly, bcFNO and spFNO retain the same inference architecture as FNO and incur the same online cost. 
The preceding results show that the resulting gains in long-time amplitude control and coarsening recovery do not require additional online evaluation cost.

\FloatBarrier

\section{Conclusion and Future Work}
\label{sec:conclusion}
In this work, we develop Fourier neural operators for long-time Cahn--Hilliard prediction under coarse physical-time supervision.
After training, such operators can complement numerical solvers 
when many trajectory evaluations are required.
Although the learned one-step map can be accurate, repeated application of the unregularized FNO produced free-energy growth, amplitude excursions outside the nominal phase-field range, and errors in late-time coarsening quantities.

The bound-conforming FNO (bcFNO) reduced the large amplitude excursions observed in these predictions without changing the online evaluation cost. 
Structure-function regularization was then added to the bcFNO training objective to form spFNO.
For the representative cases $(N=128,\delta t=0.2)$ and  $(N=256,\delta t=0.1)$, spFNO gave normalized structure functions closer to the numerical reference, which is the quantity directly constrained by the added loss. 
The fitted $L^3(t)$ slope and $k_1(t)$ exponent, derived from the prediction spectra, also moved closer to their reference values. 
The collapse score increased for $(128,0.2)$ and was unchanged for $(256,0.1)$.

These results support a hierarchy between amplitude control and statistical correction. 
When a long-time prediction develops large amplitude excursions and nonphysical free-energy growth, its spectral quantities do not provide a reliable characterization of Cahn--Hilliard coarsening. 
Controlling large excursions outside the nominal phase-field range and maintaining physically consistent free-energy behavior are therefore prerequisites for interpreting long-time coarsening statistics, although amplitude control and physically consistent free-energy behavior alone are insufficient to guarantee statistical accuracy.
The statistical correction in Stage~II is formulated for the late-time coarsening regime, where dynamic scaling provides a natural basis for structure-function supervision. 
Extending this approach to the full Cahn--Hilliard evolution remains an important direction for future work, since the early instability, structure-formation, and late-time coarsening regimes are
governed by distinct physical mechanisms. 
Such an extension may require stage-aware physical representations and multiscale neural-operator architectures that resolve the temporal and spatial scales relevant to each regime.

More broadly, nonlinear PDEs, including turbulent flows and other systems with self-similar evolution, may possess statistical structure that characterizes relevant aspects of their long-time behavior. 
Incorporating such system-specific information into operator training could complement pointwise state losses when long-horizon prediction is the objective.
Developing transferable physically meaningful constraints, and determining how equation-appropriate principal physical constraints should be adapted across equations and dynamical regimes, is a natural direction for future work.
In particular, accurate local one-step fitting alone need not ensure physically credible long-time predictions.

The Python source code is publicly available at  \url{https://github.com/Li-Yingli/Structure-preserving-FNO.git}.

\bibliography{ch}

\clearpage
\appendix
\section{Implementation Details}\label{app:implementation}
\subsection{Training Setup}

Table~\ref{tab:training_parameters} summarizes the loss parameters for the two representative source-free cases and the ranges used in the remaining experiments.

\begin{table}[!t]
\centering
\caption{Loss parameters for the representative source-free cases and their ranges across the remaining experiments.}
\label{tab:training_parameters}
\small
\begin{tabular}{@{}lccc@{}}
\toprule
Parameter & $(128,0.2)$ & $(256,0.1)$ & Range \\
\midrule
$\lambda_{\mathrm{bc}}$ & 0.18 & 0.18 & $[0.10, 0.33]$ \\
$\lambda_{\mathrm{sk}}$ & 0.22 & 0.20 & $[0.20, 0.80]$ \\
$w_{\mathrm{shape}}$ & 2.1 & 2.0 & $[1.2, 2.1]$ \\
$w_{\mathrm{null}}$ & 1.0 & 1.1 & $[0.6, 1.1]$ \\
$w_{\mathrm{peak}}$ & 1.6 & 1.5 & $[1.2, 2.0]$ \\
$\tau_{\rho}$ & 0.003 & 0.003 & $[0.002, 0.003]$ \\
$f_{\mathrm{peak}}$ & 0.30 & 0.24 & $[0.23, 0.30]$ \\
\bottomrule
\end{tabular}
\end{table}

\subsection{Reference-defined coarsening interval}
\label{app:regime_classification}

For each held-out numerical-reference trajectory, we identify a contiguous interval for coarsening evaluations. Let
\[
L(t)=\frac{1}{k_1(t)},\qquad
\gamma(t)=\frac{d\log L}{d\log t},\qquad
R(t)=\frac{L(t)}{\sqrt{\kappa}},
\]
where $\gamma(t)$ is evaluated using a 15-frame local window. A frame is a candidate coarsening frame when
\[
R(t)>8
\qquad\text{and}\qquad
\left|\gamma(t)-\frac13\right|<0.12.
\]
To suppress isolated changes in the local slope, entry into the coarsening interval requires five consecutive candidate frames. The interval ends when either $L(t)\ge0.9L_{\mathrm{box}}$, or $\gamma(t)<0.25$ and $R(t)\ge20$ persists for eight consecutive frames. 

The interval is determined from the numerical-reference trajectory only. 
For a fixed $(N,\delta t)$, its frame indices are then applied directly to the reference, bcFNO, and spFNO trajectories. The intervals used for the reported cases are listed in Table~\ref{tab:regime2}.
\begin{table}[t]
    \centering
    \caption{Reference-defined coarsening intervals used for the finite-regime fits. The same interval is used for the numerical reference, bcFNO, and spFNO.}    \label{tab:regime2}
    \begin{tabular}{ccc}
        \toprule
        $N$ & $\delta t$ & Coarsening-time interval $t$ \\
        \midrule
        128 & 1.0 & $[11500,\,29500]$ \\
        128 & 0.5 & $[6720,\,29500]$ \\
        128 & 0.2 & $[8460,\,29500]$ \\
        128 & 0.1 & $[8300,\,29500]$ \\
        \midrule
        256 & 1.0 & $[7500,\,29500]$ \\
        256 & 0.5 & $[6580,\,29500]$ \\
        256 & 0.2 & $[6400,\,29500]$ \\
        256 & 0.1 & $[6520,\,29500]$ \\
        \bottomrule
    \end{tabular}
\end{table}

\subsection{Specification of the Manufactured Forced Trajectory}
\label{app:manufactured_trajectory}

For reproducibility, we specify the exact trajectory used in
Sec.~\ref{subsec:manufactured_solution}. 
Let $q=2\pi/m$ and let $\mathcal K_+$ contain one representative from each conjugate Fourier pair with $2\leq |\mathbf n|\leq 8$. We prescribe
\begin{equation*}
  \phi^\star(\mathbf{x},t)
  =
  \bar\phi+
  \sum_{\mathbf n\in\mathcal K_+}
  a_{\mathbf n}(t)
  \cos\left(q\mathbf n\cdot\mathbf x+\vartheta_{\mathbf n}\right),
  \label{eq:app_manufactured_solution}
\end{equation*}
where the phases $\vartheta_{\mathbf n}$ are drawn once for each trajectory and remain fixed in time. 
The modal amplitudes are
\begin{equation*}
  a_{\mathbf n}(t)
  =
  \frac{\sqrt{2}A_{\rm rms}w_{\mathbf n}(t)}
  {\left(\sum_{\mathbf k\in\mathcal K_+}w_{\mathbf k}^2(t)\right)^{1/2}},
  \qquad
  w_{\mathbf n}(t)
  =
  \exp\left[
  -\frac{\left(|\mathbf n|-\nu_c(t)\right)^2}
  {2c_\sigma^2\nu_c^2(t)}
  \right],
  \label{eq:app_manufactured_coefficients}
\end{equation*}
with
\begin{equation*}
  \nu_c(t)=\nu_0\left(1+\frac{t}{t_0}\right)^{-1/3}.
  \label{eq:app_manufactured_center}
\end{equation*}
The factor $\sqrt{2}$ gives
\begin{equation*}
  \frac{1}{|\mathbb{T}_m^2|}
  \int_{\mathbb{T}_m^2}
  \left(\phi^\star-\bar\phi\right)^2\,d\mathbf{x}
  = A_{\rm rms}^2.
\end{equation*}
For each prescribed trajectory, $S^\star$ is evaluated analytically from
Eq.~\eqref{eq:manufactured_source}.

\begin{table}[t]
\centering
\caption{Parameters of the manufactured forced trajectory for
field-level prediction evaluation.}
\label{tab:manufactured_parameters}
\begin{tabular}{cccccccccc}
\toprule
$m$ & $N$ & $M$ & $\kappa$ & $\bar\phi$ &
$A_{\rm rms}$ & $c_\sigma$ & $\nu_0$ & $t_0$ & $\Delta T$ \\
\midrule
256 & 128 & 1 & 0.5 & 0.1 & 0.25 & 0.35 & 6.581 &
$15000/7$ & 20 \\
\bottomrule
\end{tabular}
\end{table}

Table~\ref{tab:manufactured_parameters} lists the parameters of the held-out manufactured trajectory used for field-level prediction evaluation; $m=256$ is the side length of $\mathbb T_m^2$. 
Its Fourier phases are fixed in time. The training trajectories use different phase realizations and may use different envelope parameters.

\subsection{Online Timing Setup}
\label{app:timing}

The spectral solver was implemented using NumPy/SciPy and executed on an Intel Xeon Gold~6132 CPU system.
The neural operators were implemented in PyTorch and evaluated on an NVIDIA Tesla V100 PCIe GPU with 16~GB of memory.
The software environment used NumPy~2.2.6, SciPy~1.15.3,
PyTorch~2.10.0, and CUDA~12.8.

Each timed evaluation generated one trajectory over
$0\leq t\leq30000$, with output states produced at intervals of
$\Delta T=20$.
The spectral-solver benchmark used no untimed warm-up
(\texttt{warmup=0}).
Before each timed neural-operator prediction, two untimed single-step forward passes were performed, following the command-line default
(\texttt{--warmup 2}); these warm-up passes did not generate complete
trajectories.
CPU wall-clock times were measured using \texttt{time.perf\_counter()}.
GPU execution times were measured using CUDA events, with
\texttt{torch.cuda.synchronize()} called before and after each timed prediction.
The reported values are the medians of five timed repetitions.

\section{Supplementary Results}\label{sec:appendix}

This section presents the training-seed comparison and the energy, amplitude, and mass-conservation results for source-free cases not shown in detail in the main text.

\subsection{Sensitivity to Training Initialization}
\label{app:multiseed}

We repeat the FNO and bcFNO experiments using three independent training seeds for $(N,\delta t)=(256,0.1)$ under Sampling~B.
The models use the same held-out trajectory, architecture, training setup, and prediction settings.

\begin{table}[!t]
\centering
\caption{Maximum field magnitude over held-out long-time rollouts for three training seeds at $(N,\delta t)=(256,0.1)$ under Sampling~B.}
\label{tab:multiseed_amplitude}
\small
\begin{tabular}{lccc}
\toprule
Model & Seed 0 & Seed 1 & Seed 2 \\
\midrule
FNO   & 3.65 & 27.6 & 2.53 \\
bcFNO & 1.15 & 1.17 & 1.23 \\
\bottomrule
\end{tabular}
\end{table}

For the three tested seeds, the FNO maximum field magnitude ranges from $2.53$ to $27.6$, whereas the bcFNO maximum ranges from $1.15$ to $1.23$. This comparison concerns maximum field magnitude over the held-out prediction.

\subsection{Additional Energy and Amplitude Predictions}\label{app:fno_failures}
    
\begin{figure}[!t]
  \centering
  \begin{minipage}{0.32\textwidth}
    \centering
    \includegraphics[width=\linewidth]{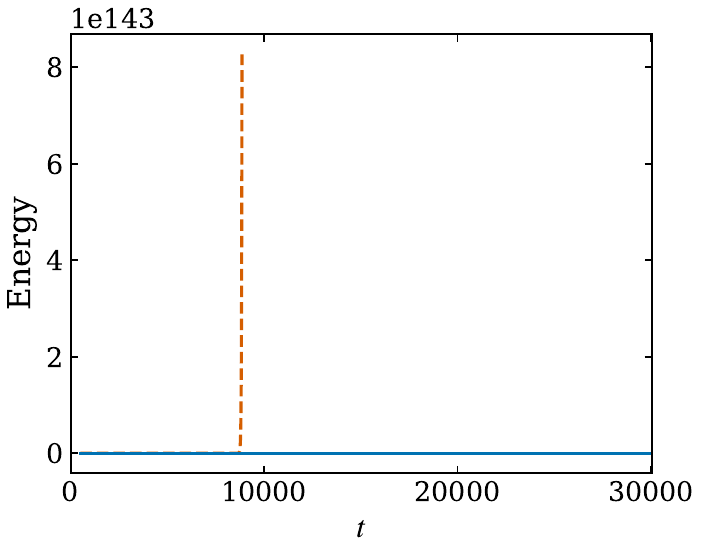}\\[0.25em]
    \includegraphics[width=\linewidth]{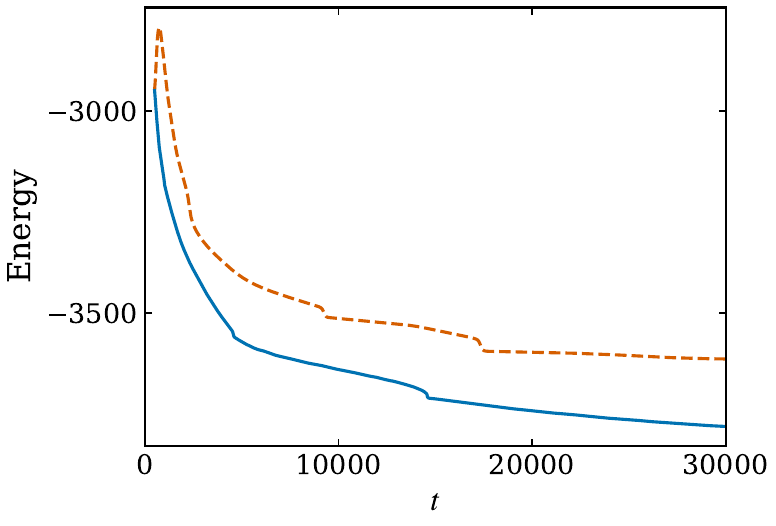}
  \end{minipage}
  \hfill
  \begin{minipage}{0.32\textwidth}
    \centering
    \includegraphics[width=\linewidth]{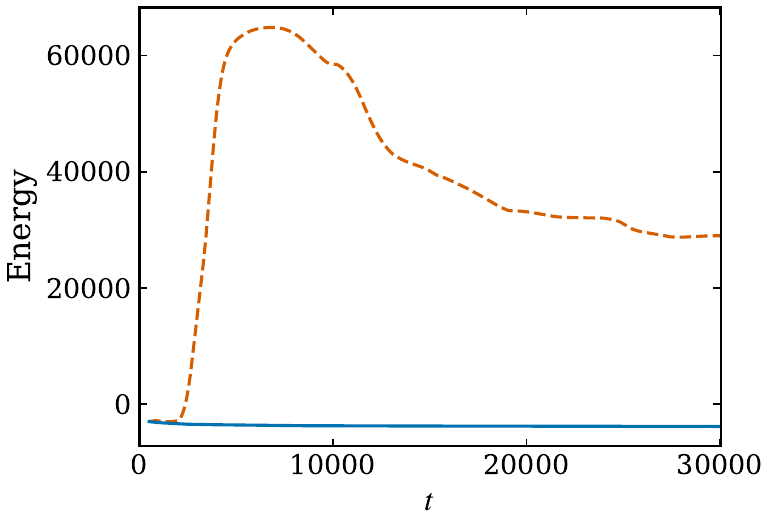}\\[0.25em]
    \includegraphics[width=\linewidth]{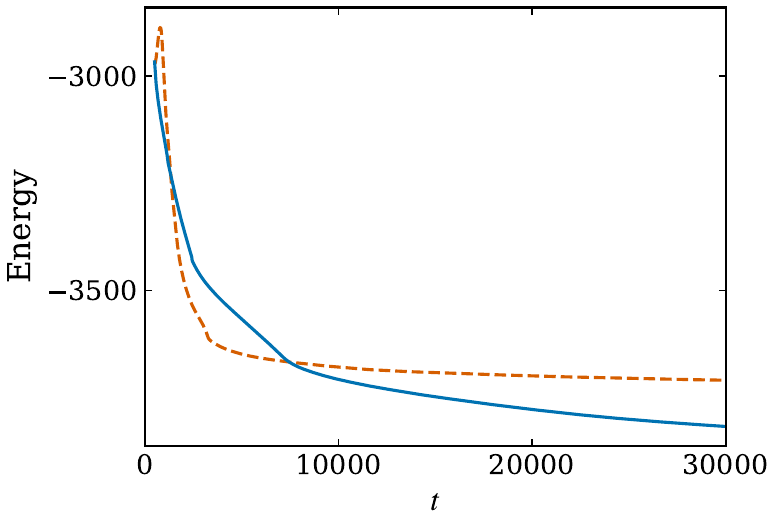}
  \end{minipage}
  \hfill
  \begin{minipage}{0.32\textwidth}
    \centering
    \includegraphics[width=\linewidth]{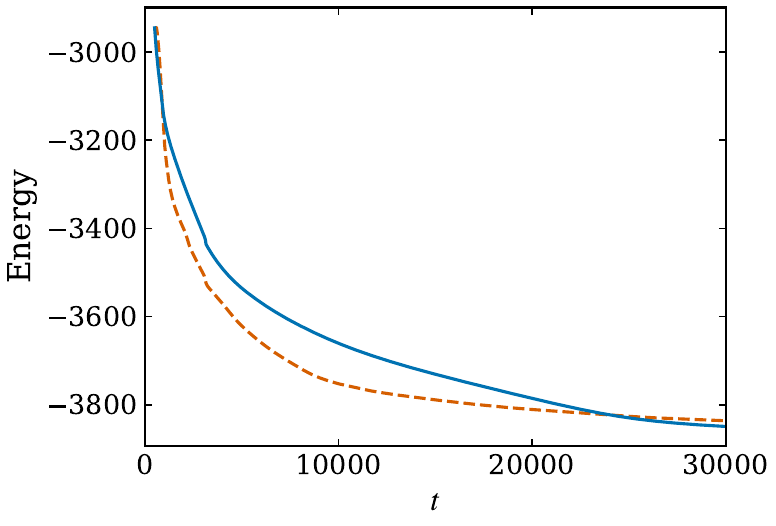}\\[0.25em]
    \includegraphics[width=\linewidth]{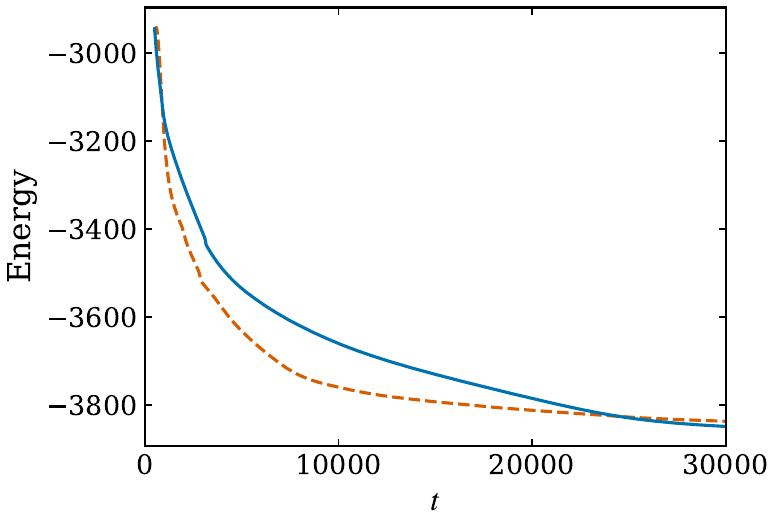}
  \end{minipage}
 \caption{Evolution of the total free energy $E(t)$ for FNO (top row) and bcFNO (bottom row). Blue and orange curves denote the reference trajectories and model predictions, respectively.  The columns correspond to $\delta t=1$, $0.5$, and $0.1$ (from left to right). $N=128$.}
\label{fig:app_fno_failure_energy_N128}
\end{figure}

\begin{figure}[!t]
  \centering
  \begin{minipage}{0.32\textwidth}
    \centering
    \includegraphics[width=\linewidth]{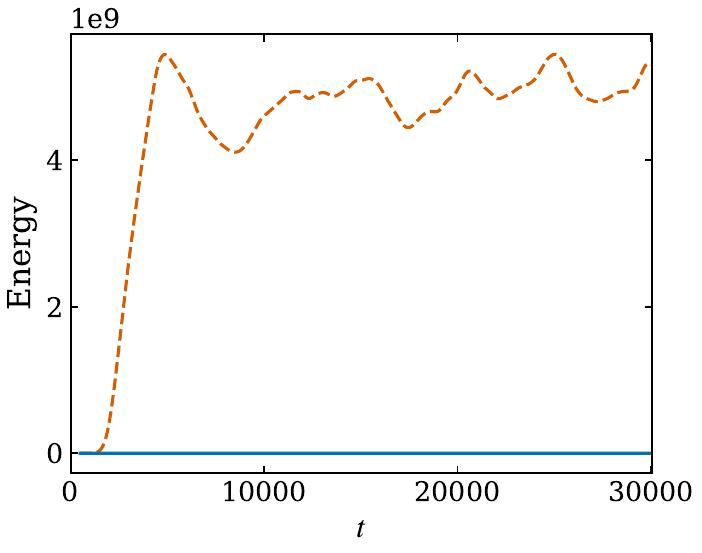}\\[0.25em]
    \includegraphics[width=\linewidth]{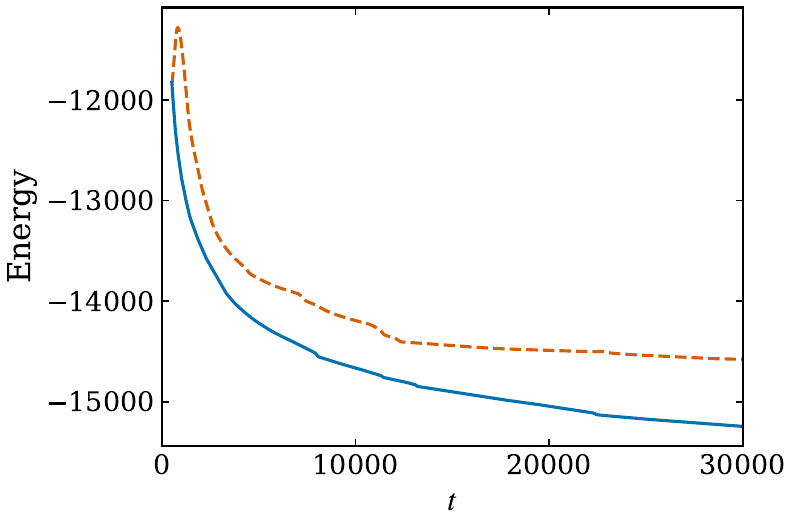}
  \end{minipage}
  \hfill
  \begin{minipage}{0.32\textwidth}
    \centering
    \includegraphics[width=\linewidth]{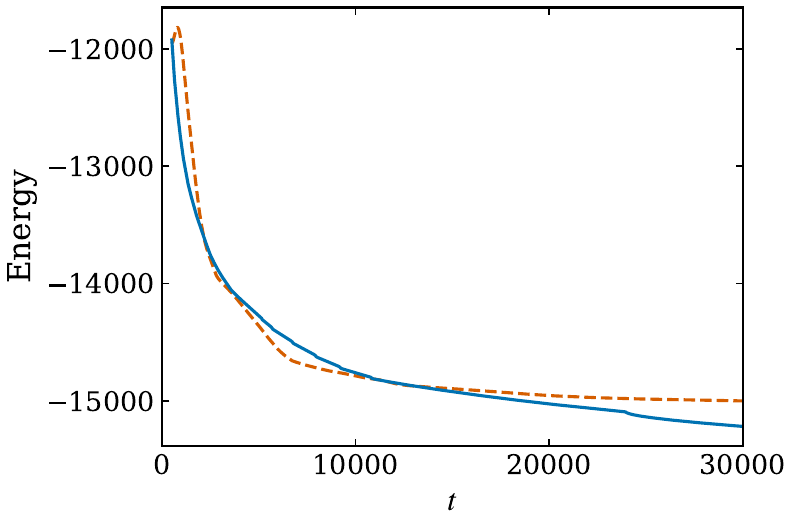}\\[0.25em]
    \includegraphics[width=\linewidth]{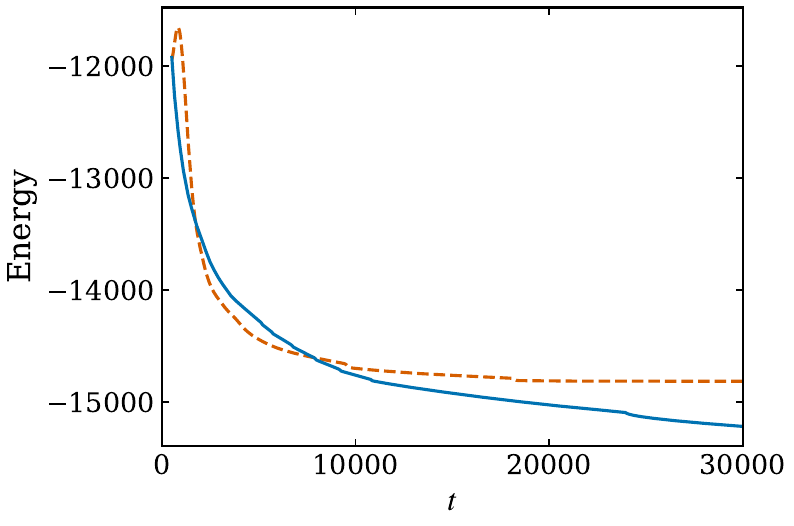}
  \end{minipage}
  \hfill
  \begin{minipage}{0.32\textwidth}
    \centering
    \includegraphics[width=\linewidth]{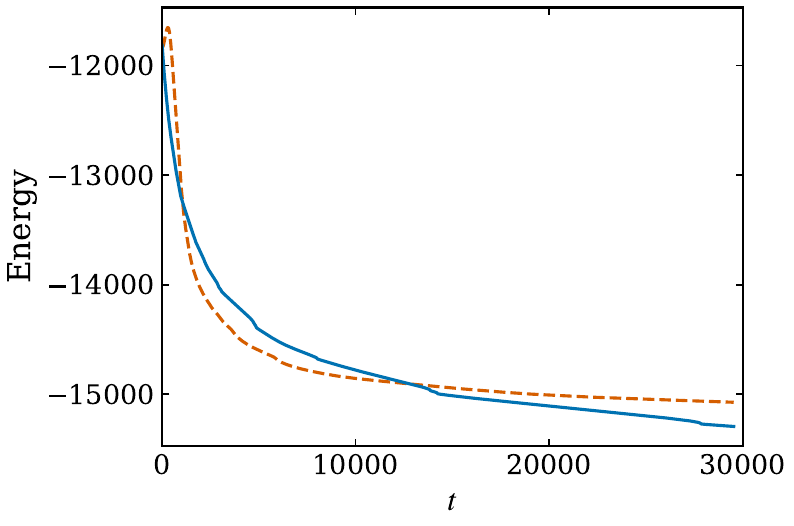}\\[0.25em]
    \includegraphics[width=\linewidth]{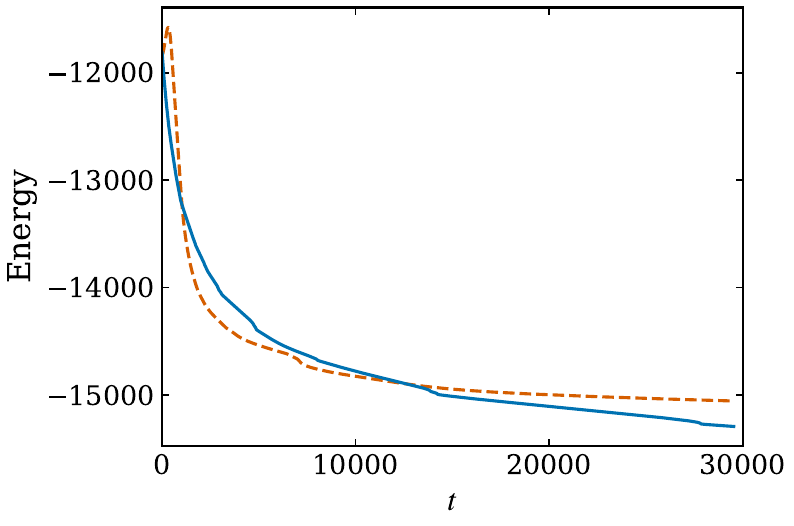}
  \end{minipage}
 \caption{Evolution of the free energy $E(t)$ for FNO (top row) and bcFNO (bottom row). Blue and orange curves denote the reference trajectories and model predictions, respectively. The columns correspond to $\delta t=1$, $0.5$, and $0.2$ (from left to right). $N=256$.}
\label{fig:app_fno_failure_energy_N256}
\end{figure}

The main text presents $(N,\delta t)=(128,0.2)$ and $(256,0.1)$ in detail.
Figures~\ref{fig:app_fno_failure_energy_N128}
and~\ref{fig:app_fno_failure_energy_N256} show the corresponding
free energy for the other six cases.

In several FNO predictions, the energy departs from the dissipative trend of the reference trajectory.
The corresponding bcFNO substantially reduces these deviations for the sampling distributions.
For the coarsest solver-generated targets, a short initial energy increase may remain, but it is distinct from the persistent late-time growth observed in the unstable FNO predictions.

\begin{figure}[!t]
  \centering
  \begin{subfigure}[t]{0.32\textwidth}
    \centering
    \includegraphics[width=\linewidth]{Figures/N128dt1/bcFNO/energy_N128_sdt1_bd0p33_lsk0_tw15252535_si20_phys500_strat.pdf}
  \end{subfigure}
  \hfill
  \begin{subfigure}[t]{0.32\textwidth}
    \centering
    \includegraphics[width=\linewidth]{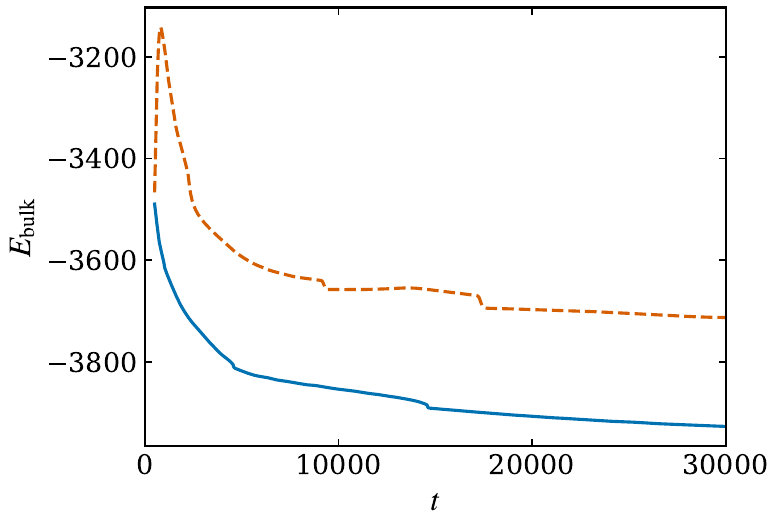}
  \end{subfigure}
  \hfill
  \begin{subfigure}[t]{0.32\textwidth}
    \centering
    \includegraphics[width=\linewidth]{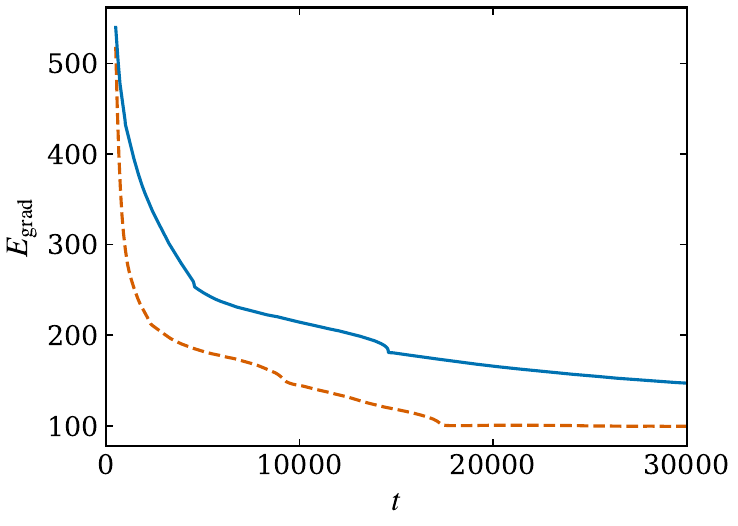}
  \end{subfigure}
  \par\vspace{0.5em}
  \begin{subfigure}[t]{0.32\textwidth}
    \centering
    \includegraphics[width=\linewidth]{Figures/N256dt1/bcFNO/energy_N256_sdt1_bd0p33_lsk0_nm48_tw15255015_si20_phys500_strat.pdf}
  \end{subfigure}
  \hfill
  \begin{subfigure}[t]{0.32\textwidth}
    \centering
    \includegraphics[width=\linewidth]{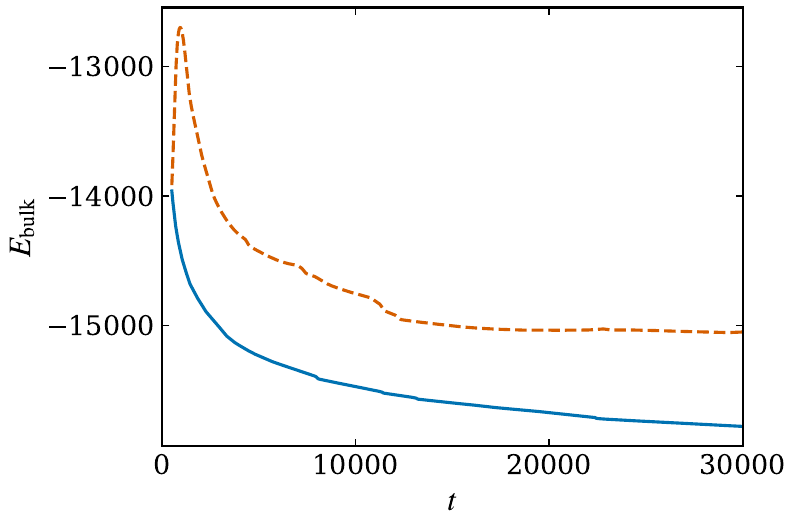}
  \end{subfigure}
  \hfill
  \begin{subfigure}[t]{0.32\textwidth}
    \centering
    \includegraphics[width=\linewidth]{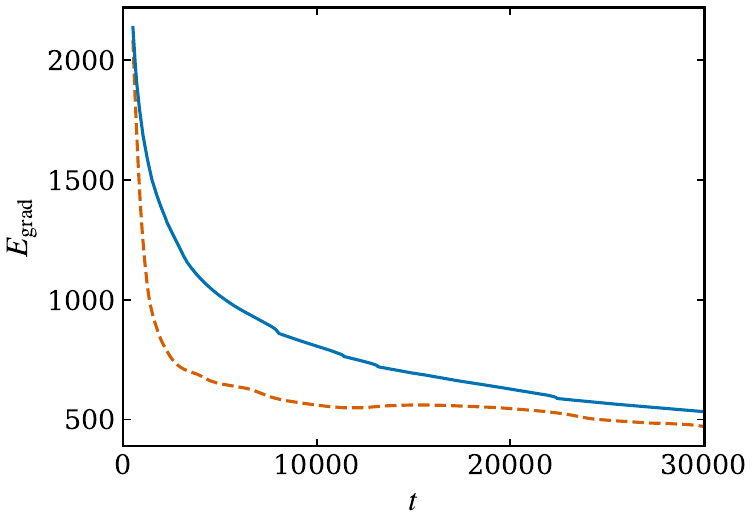}
  \end{subfigure}
    \caption{Decomposition of the bcFNO free energy at $\delta t=1$. The columns show the total energy $E(t)$, the bulk contribution $E_{\mathrm{bulk}}(t)$, and the interfacial-gradient contribution $E_{\mathrm{grad}}(t)$ (from left to right). The rows correspond to $N=128$ (top) and $N=256$ (bottom). Blue and orange curves denote the reference trajectories and bcFNO predictions, respectively.}

  \label{fig:app_bcfno_energy_decomposition}
\end{figure}

To examine the short initial energy increase at $\delta t=1$, we decompose the free energy as
\[
E(t)=E_{\mathrm{bulk}}(t)+E_{\mathrm{grad}}(t),
\]
where $E_{\mathrm{bulk}}$ is the double-well contribution and
$E_{\mathrm{grad}}$ is the interfacial-gradient contribution.
As shown in Fig.~\ref{fig:app_bcfno_energy_decomposition}, the initial
increase in the bcFNO total energy is produced by an increase in
$E_{\mathrm{bulk}}$ that exceeds the concurrent decrease in
$E_{\mathrm{grad}}$.
After this initial adjustment, both contributions decrease and the total
energy resumes a dissipative trend.
This transient is therefore qualitatively distinct from unphysical late-time energy growth.

\begin{figure}[!t]
  \centering
  \begin{minipage}{0.32\textwidth}
    \centering
    \includegraphics[width=\linewidth]{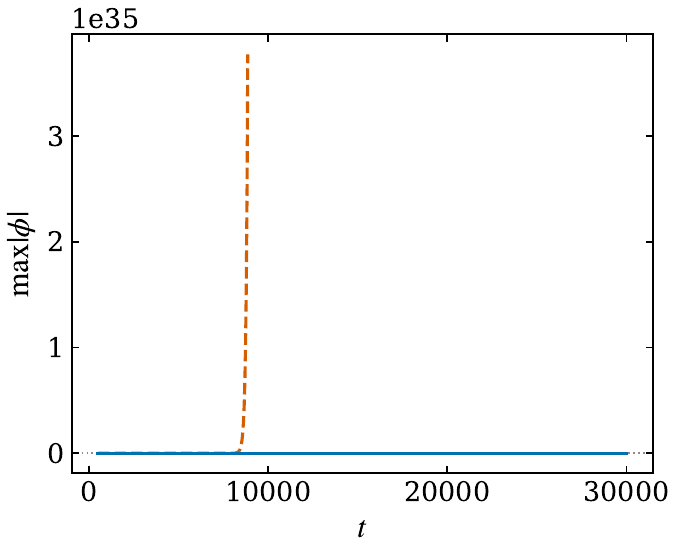}\\[0.25em]
    \includegraphics[width=\linewidth]{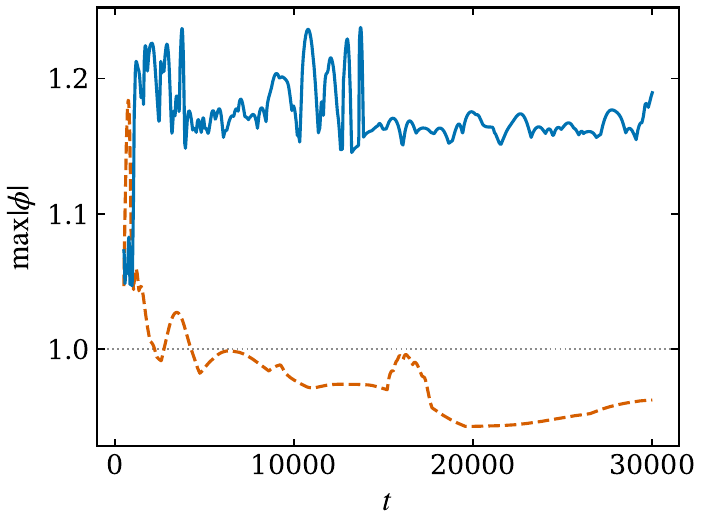}
  \end{minipage}
  \hfill
  \begin{minipage}{0.32\textwidth}
    \centering
    \includegraphics[width=\linewidth]{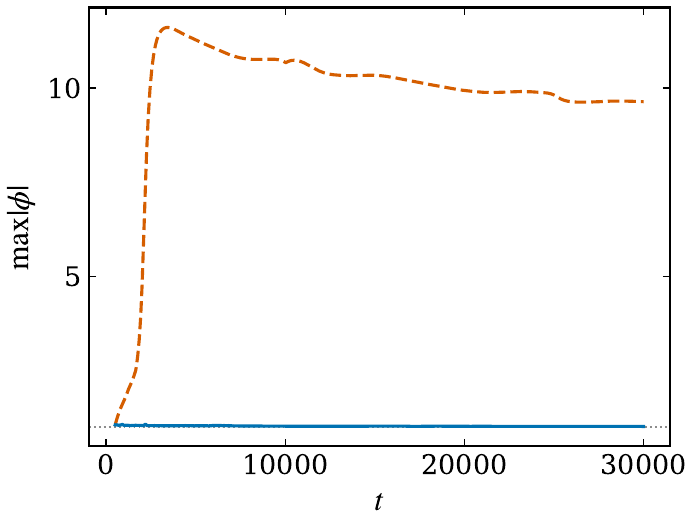}\\[0.25em]
    \includegraphics[width=\linewidth]{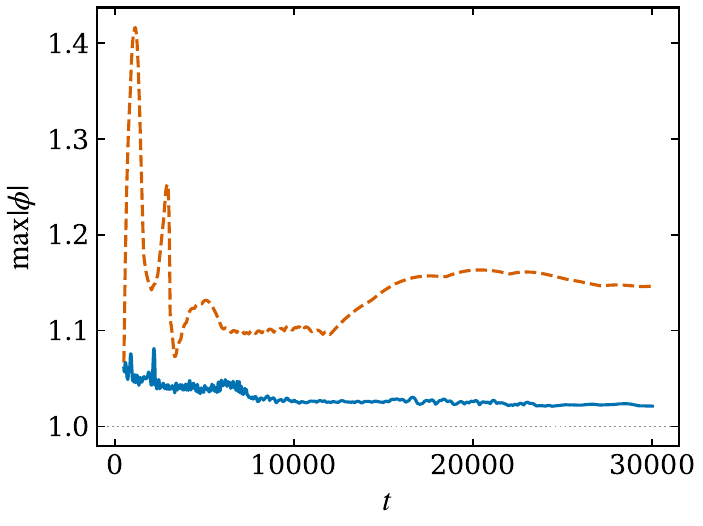}
  \end{minipage}
  \hfill
  \begin{minipage}{0.32\textwidth}
    \centering
    \includegraphics[width=\linewidth]{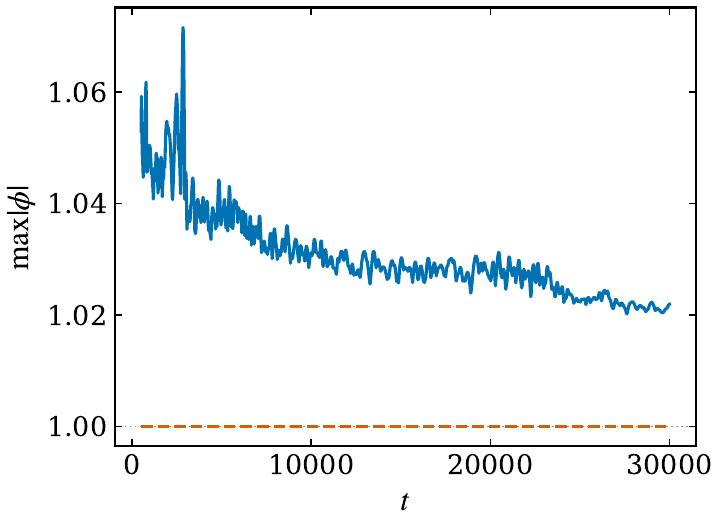}\\[0.25em]
    \includegraphics[width=\linewidth]{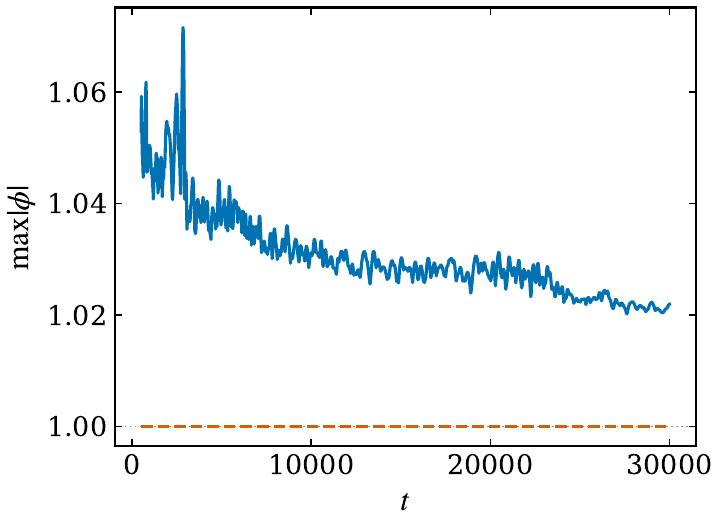}
  \end{minipage}
\caption{Maximum field magnitude $\max|\phi|$ for FNO (top) and bcFNO (bottom). Blue and orange curves denote the reference trajectories and model predictions, respectively. The columns correspond to $\delta t=1$, $0.5$, and $0.1$ (from left to right). $N=128$.}
\label{fig:app_fno_phi_max_N128}
\end{figure}

\begin{figure}[!t]
  \centering
  \begin{minipage}{0.32\textwidth}
    \centering
    \includegraphics[width=\linewidth]{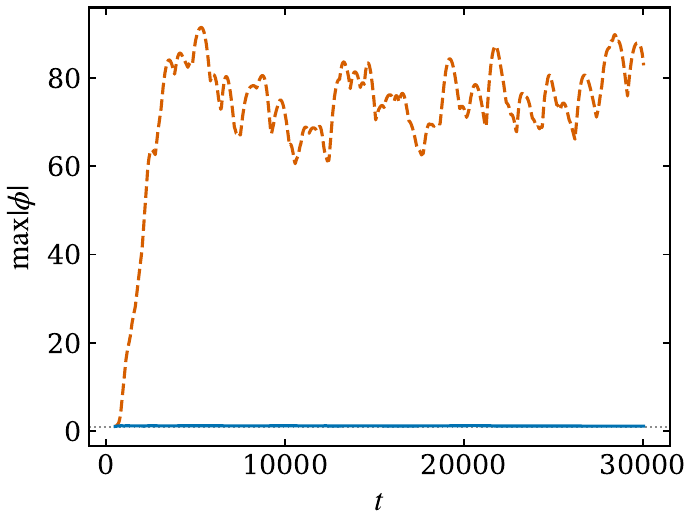}\\[0.25em]
    \includegraphics[width=\linewidth]{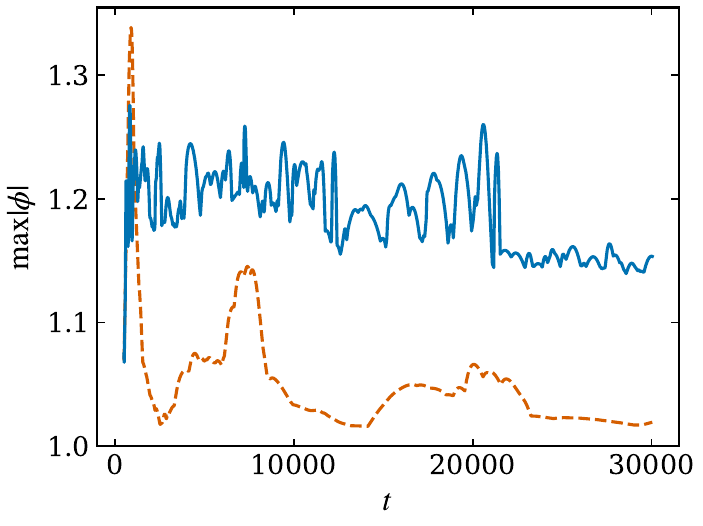}
  \end{minipage}
  \hfill
  \begin{minipage}{0.32\textwidth}
    \centering
    \includegraphics[width=\linewidth]{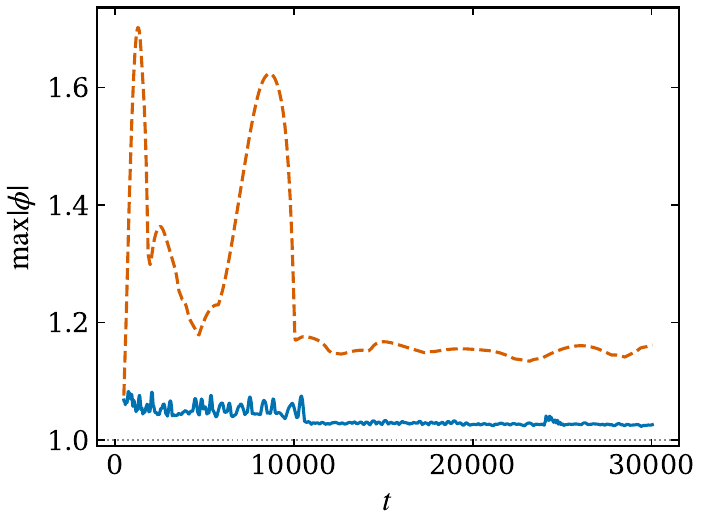}\\[0.25em]
    \includegraphics[width=\linewidth]{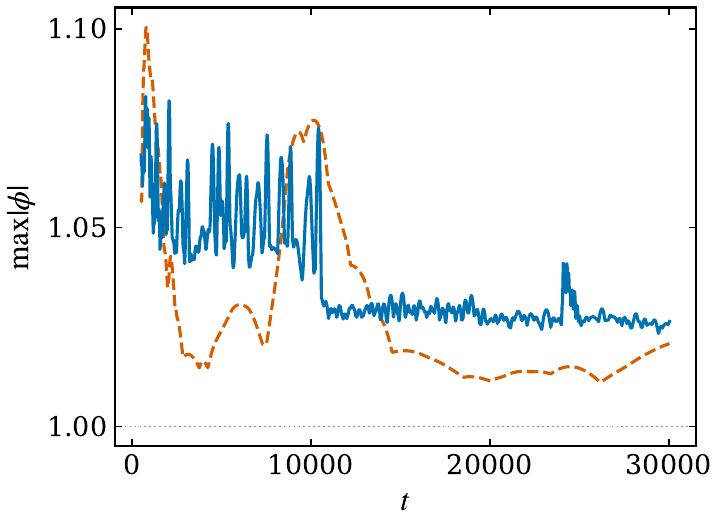}
  \end{minipage}
  \hfill
  \begin{minipage}{0.32\textwidth}
    \centering
    \includegraphics[width=\linewidth]{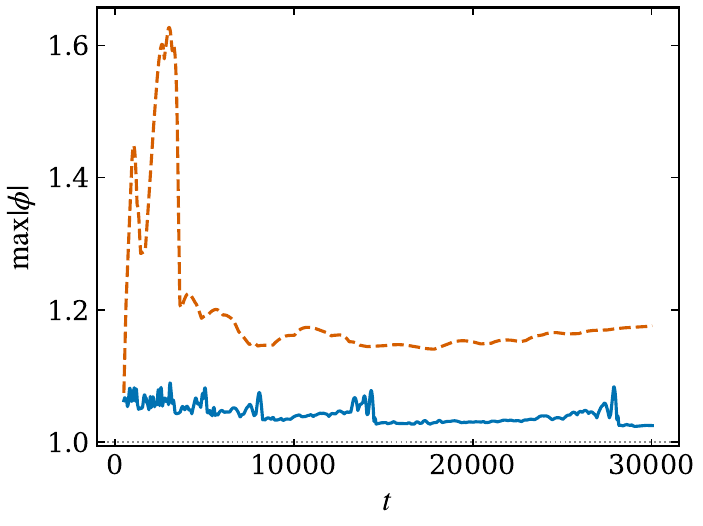}\\[0.25em]
    \includegraphics[width=\linewidth]{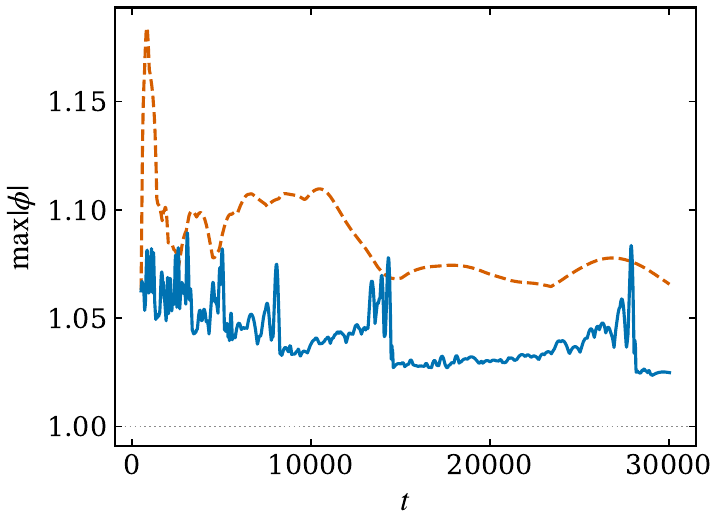}
  \end{minipage}
    \caption{Maximum field magnitude $\|\phi(\cdot,t)\|_\infty$ for FNO (top) and bcFNO (bottom) at $N=256$.
    The columns correspond to $\delta t=1$, $0.5$, and $0.2$ from left to right. Blue and orange denote the reference trajectories and model predictions, respectively.}
    \label{fig:app_fno_phi_max_N256}
\end{figure}

Figures~\ref{fig:app_fno_phi_max_N128}
and~\ref{fig:app_fno_phi_max_N256} show the maximum field magnitude for the same six cases.
Several FNO predictions exhibit sustained excursions outside the nominal
phase-field range.
The bcFNO suppresses the large excursions across the cases shown, although small initial overshoots may remain.

In summary, the energy and amplitude evaluations show that the
bound-conforming formulation consistently suppresses large amplitude
excursions and reduces the associated non-dissipative energy behavior for the cases and sampling distributions considered here.

\FloatBarrier

\subsection{Mass-Conservation Evaluations}
\label{app:mass}

Mass conservation is evaluated using
\[
\delta M(t)
=
\overline{\phi}(t)-\overline{\phi}(t_0),
\]
where $t_0$ is the initial prediction time.
This quantity measures conservation within each predicted trajectory rather than its difference from the reference trajectory.
Exact preservation of the Cahn--Hilliard spatial average corresponds to $\delta M(t)=0$.

\begin{figure}[!t]
  \centering
  \begin{minipage}{0.32\textwidth}
    \centering
    \includegraphics[width=\linewidth]{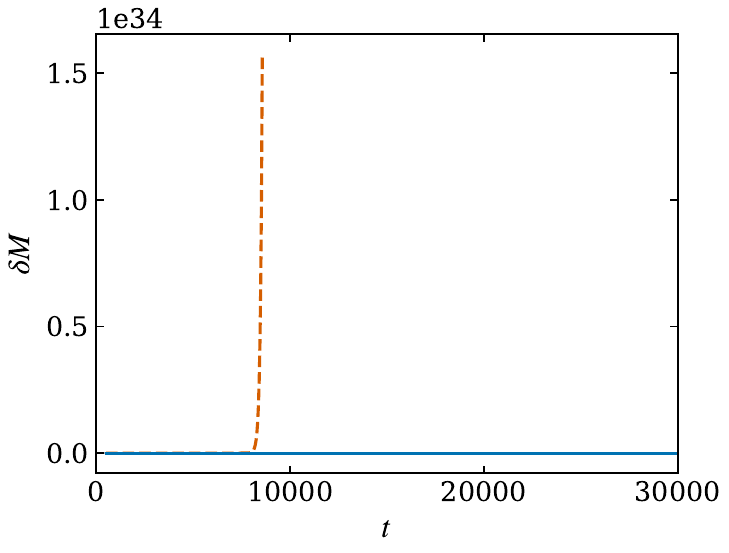}\\[0.25em]
    \includegraphics[width=\linewidth]{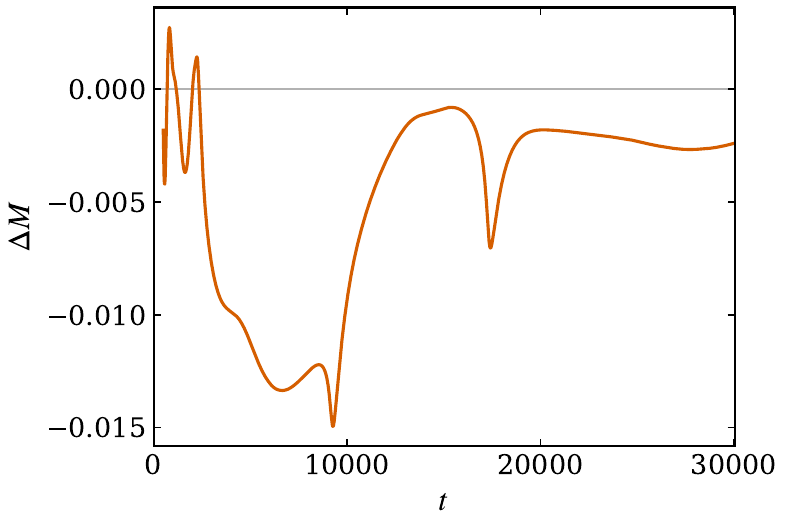}
  \end{minipage}
  \hfill
  \begin{minipage}{0.32\textwidth}
    \centering
    \includegraphics[width=\linewidth]{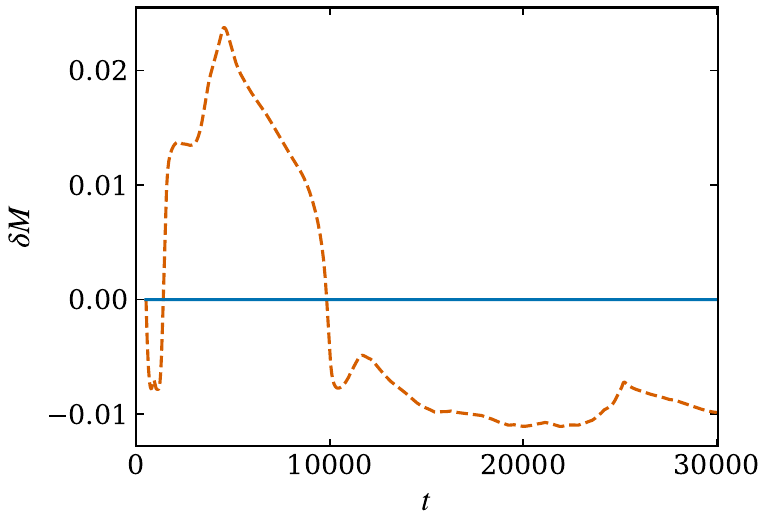}\\[0.25em]
    \includegraphics[width=\linewidth]{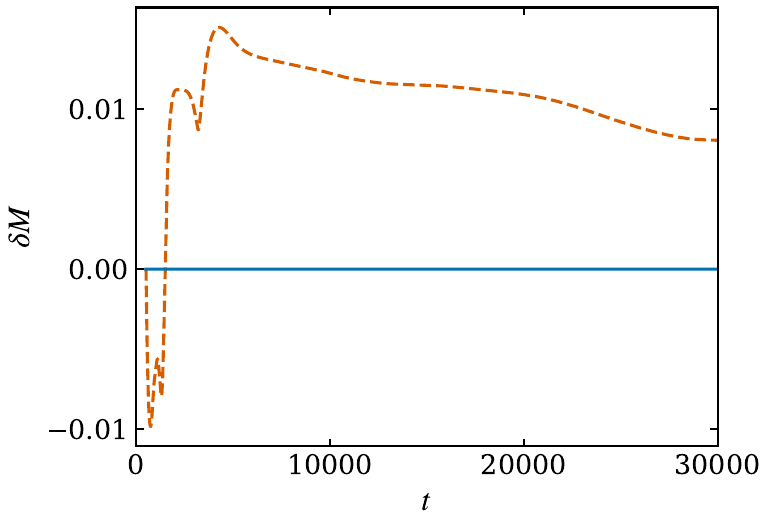}
  \end{minipage}
  \hfill
  \begin{minipage}{0.32\textwidth}
    \centering
    \includegraphics[width=\linewidth]{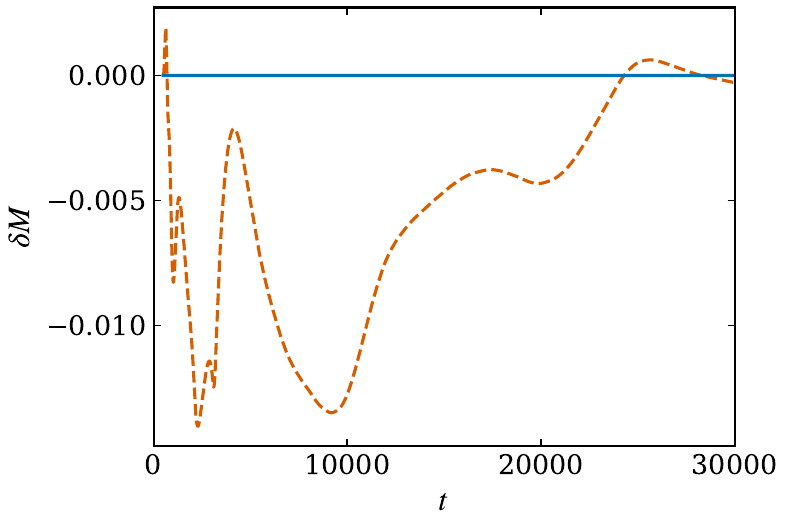}\\[0.25em]
    \includegraphics[width=\linewidth]{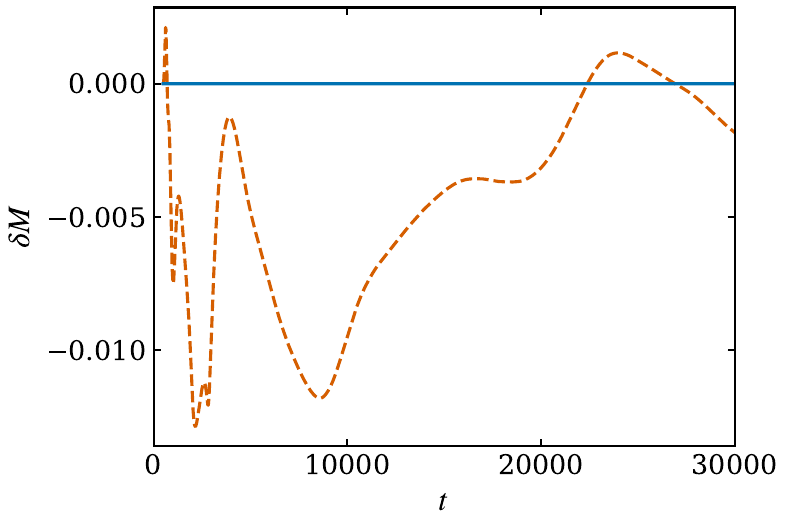}
  \end{minipage}
    \caption{Mass drift $\delta M(t)=\bar{\phi}(t)-\bar{\phi}(t_0)$ for FNO (top) and bcFNO (bottom) at $N=128$. The columns correspond to $\delta t=1$, $0.5$, and $0.1$ from left to right.
Blue and orange denote the reference trajectories and model predictions,
respectively.}
  \label{fig:app_fno_mass_conservation_N128}
\end{figure}

\begin{figure}[!t]
  \centering
  \begin{minipage}{0.32\textwidth}
    \centering
    \includegraphics[width=\linewidth]{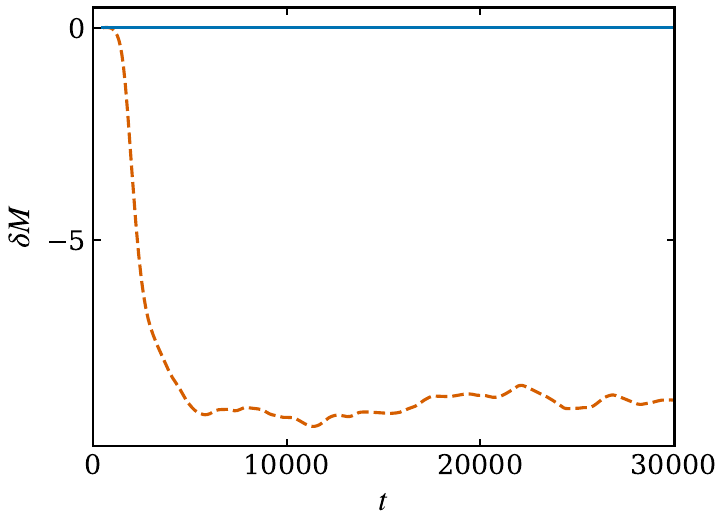}\\[0.25em]
    \includegraphics[width=\linewidth]{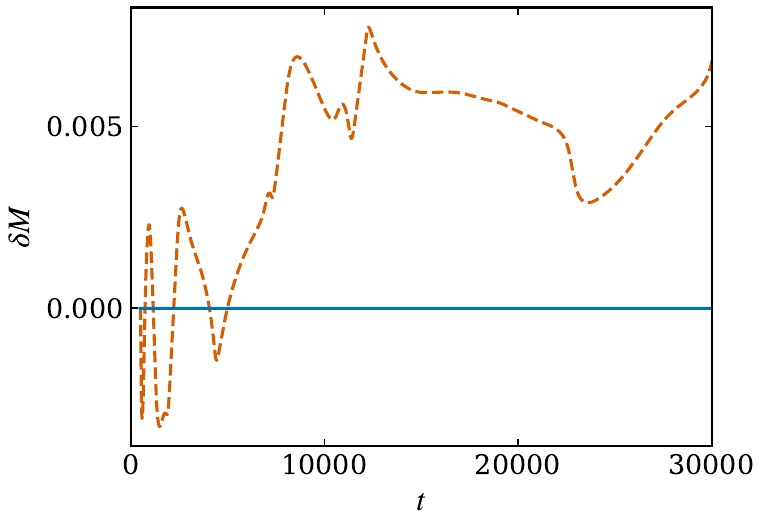}
  \end{minipage}
  \hfill
  \begin{minipage}{0.32\textwidth}
    \centering
    \includegraphics[width=\linewidth]{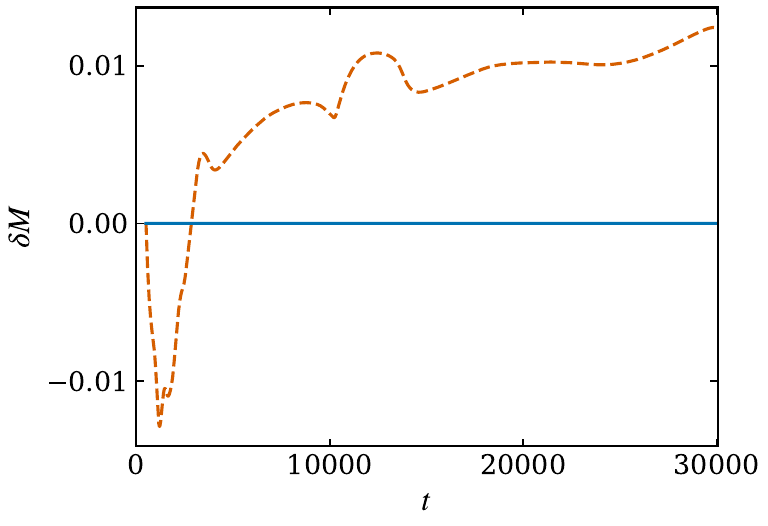}\\[0.25em]
    \includegraphics[width=\linewidth]{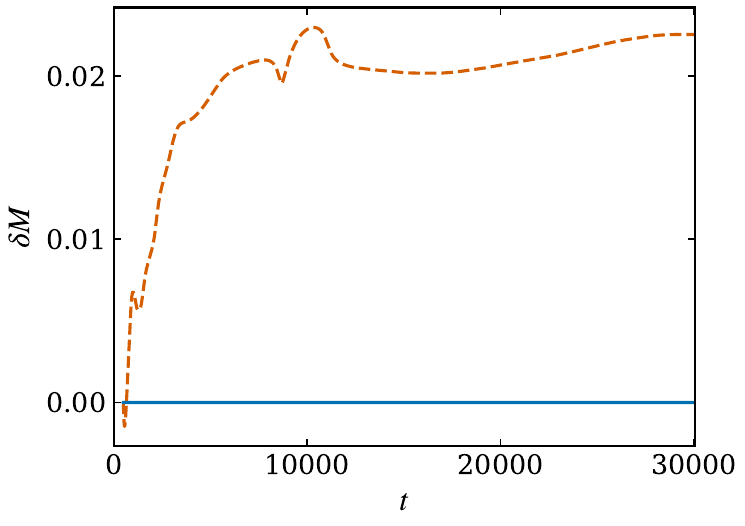}
  \end{minipage}
  \hfill
  \begin{minipage}{0.32\textwidth}
    \centering
    \includegraphics[width=\linewidth]{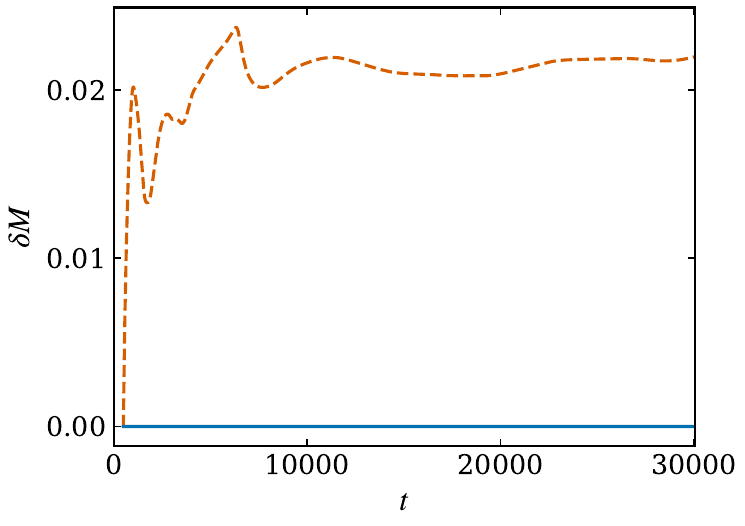}\\[0.25em]
    \includegraphics[width=\linewidth]{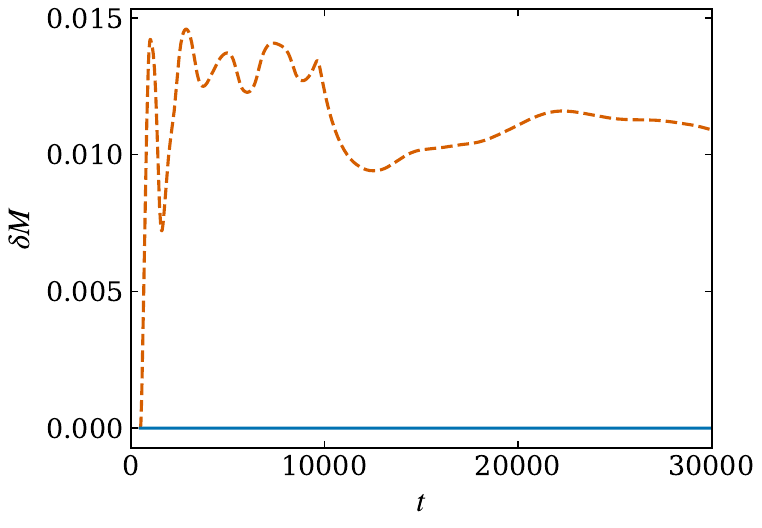}
  \end{minipage}
    \caption{Evolution of the mass drift $\delta M(t)=\bar{\phi}(t)-\bar{\phi}(t_0)$ for FNO (top row) and bcFNO (bottom row) at $N=256$. Blue and orange curves denote the reference trajectories and model predictions, respectively. The columns correspond to $\delta t=1$, $0.5$, and $0.2$ (from left to right).}
  \label{fig:app_fno_mass_conservation_N256}
\end{figure}

Figures~\ref{fig:app_fno_mass_conservation_N128}
and~\ref{fig:app_fno_mass_conservation_N256} show that appreciable mass drift is concentrated primarily in cases generated using the coarsest internal time step.
Neither FNO nor bcFNO enforces the spatial mean explicitly, and bcFNO does not provide a uniform reduction in mass drift across all cases.
Any observed improvement should therefore be interpreted as empirical rather than as a conservation property of the bound-conforming formulation.
These results show that the present models do not explicitly enforce mass conservation.

\FloatBarrier

\subsection{Results for Additional Discretizations}
\label{app:robustness}

This subsection reports bcFNO and spFNO results for the six cases not examined in detail in Sec.~\ref{subsec:exp_sk}:
$N=128$ with $\delta t\in\{1,0.5,0.1\}$ and
$N=256$ with $\delta t\in\{1,0.5,0.2\}$.
Each model is trained independently with the same coarse prediction interval $\Delta T=20$.

\begin{figure}[!t]
  \centering
  \begin{minipage}{0.32\textwidth}
    \centering
    \includegraphics[width=\linewidth]{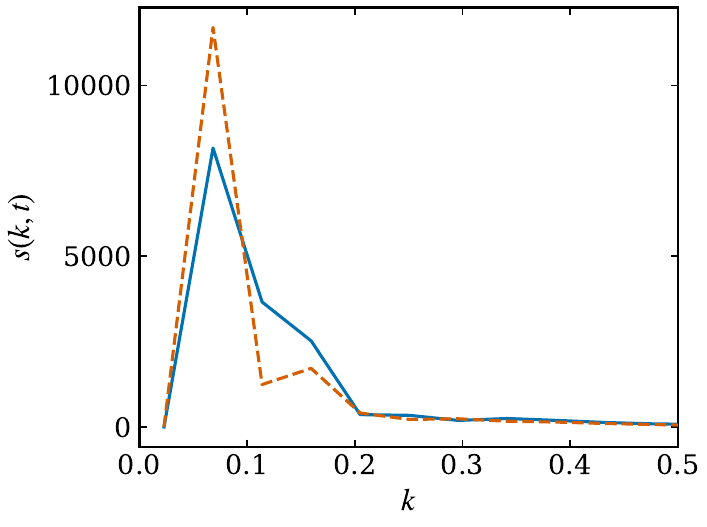}\\[0.2em]
    \includegraphics[width=\linewidth]{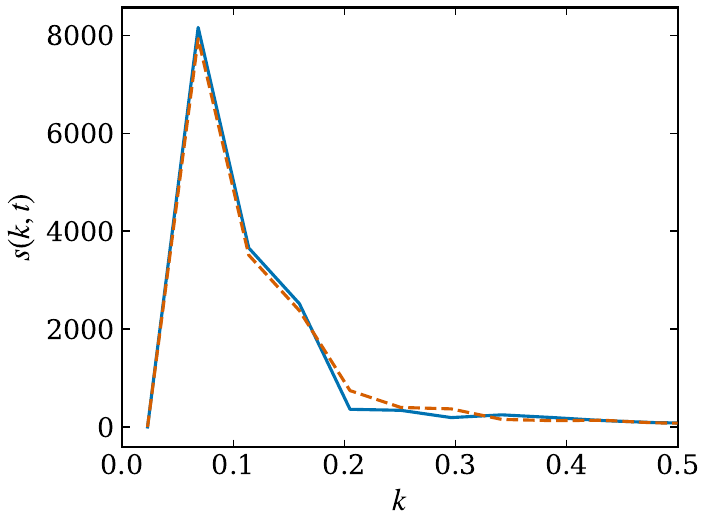}
  \end{minipage}
  \hfill
  \begin{minipage}{0.32\textwidth}
    \centering
    \includegraphics[width=\linewidth]{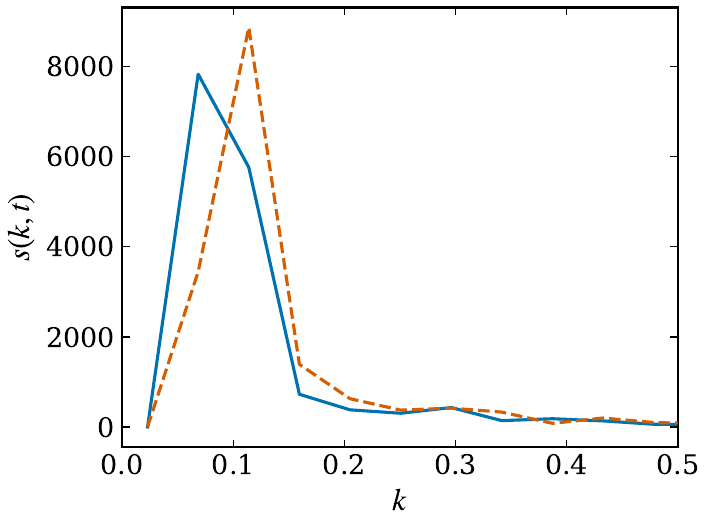}\\[0.2em]
    \includegraphics[width=\linewidth]{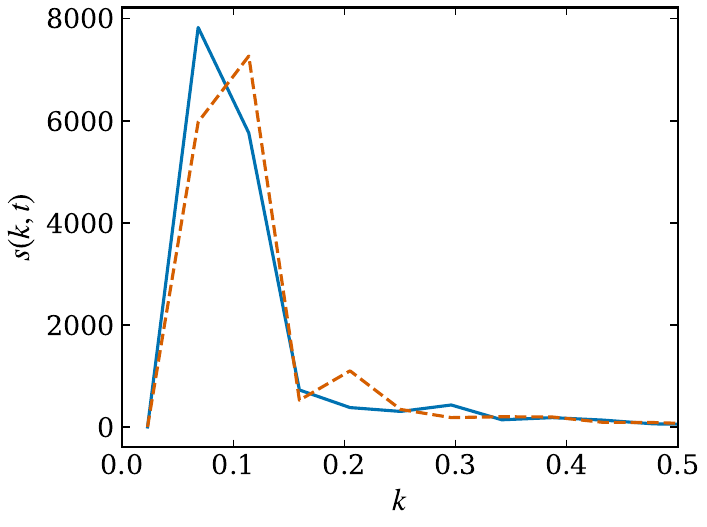}
  \end{minipage}
  \hfill
  \begin{minipage}{0.32\textwidth}
    \centering
    \includegraphics[width=\linewidth]{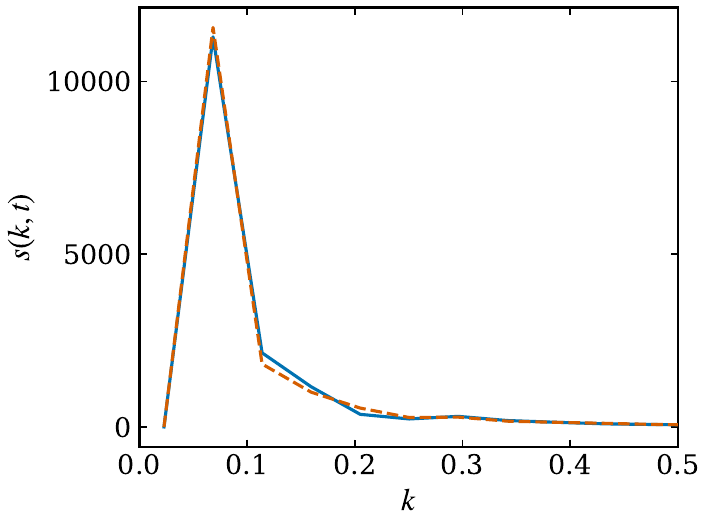}\\[0.2em]
    \includegraphics[width=\linewidth]{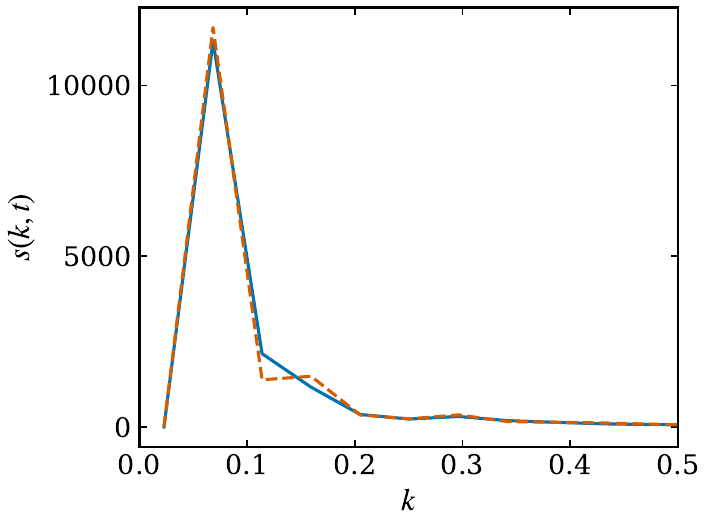}
  \end{minipage}
   \caption{Late-time normalized structure functions $s(k,t)$ at $t=29000$ for bcFNO
(top) and spFNO (bottom) at $N=128$.
The columns correspond to $\delta t=1$, $0.5$, and $0.1$ from left to right.
Blue and orange denote the reference trajectories and model predictions,
respectively.}
\label{fig:app_sk_late_N128}
\end{figure}

\begin{figure}[!t]
  \centering
  \begin{minipage}{0.32\textwidth}
    \centering
    \includegraphics[width=\linewidth]{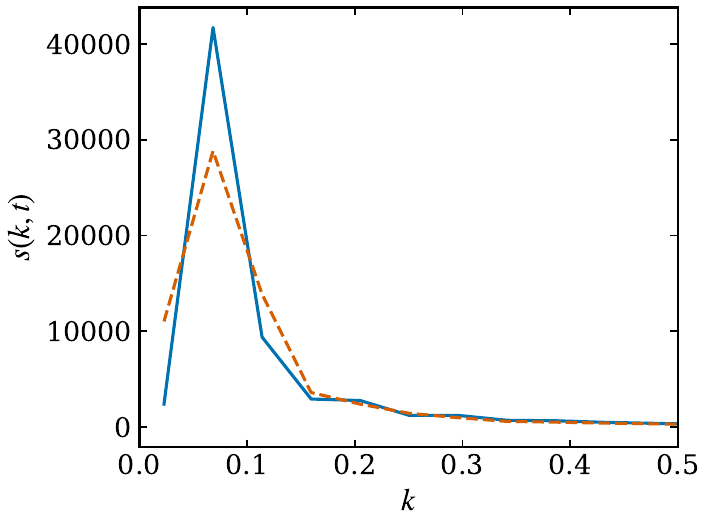}\\[0.2em]
    \includegraphics[width=\linewidth]{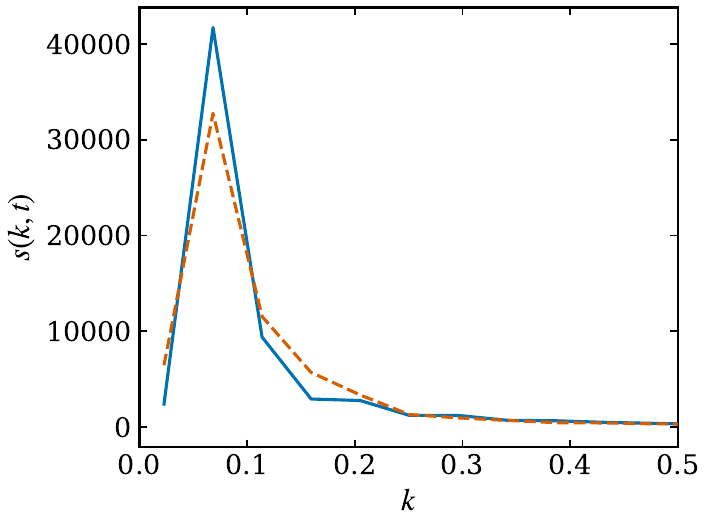}
  \end{minipage}
  \hfill
  \begin{minipage}{0.32\textwidth}
    \centering
    \includegraphics[width=\linewidth]{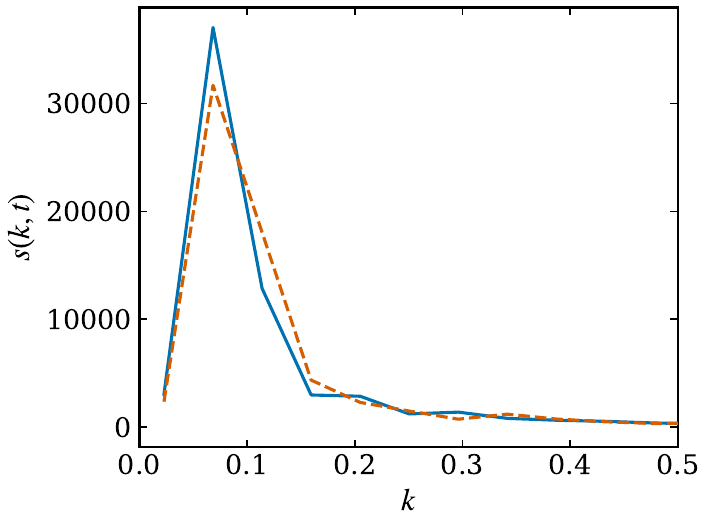}\\[0.2em]
    \includegraphics[width=\linewidth]{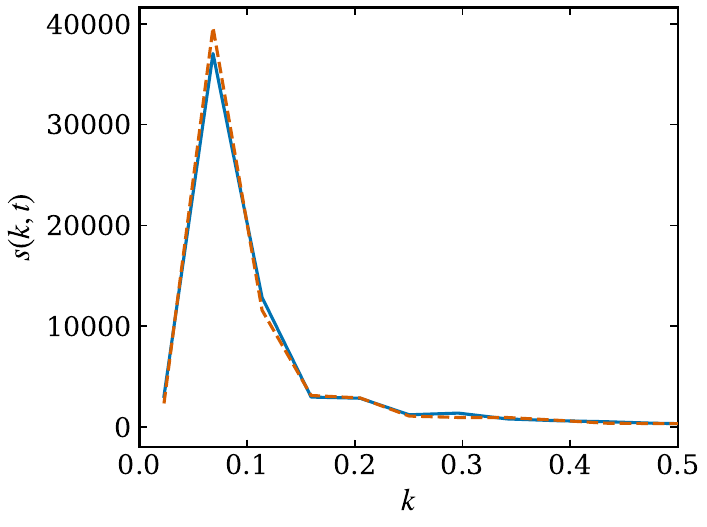}
  \end{minipage}
  \hfill
  \begin{minipage}{0.32\textwidth}
    \centering
    \includegraphics[width=\linewidth]{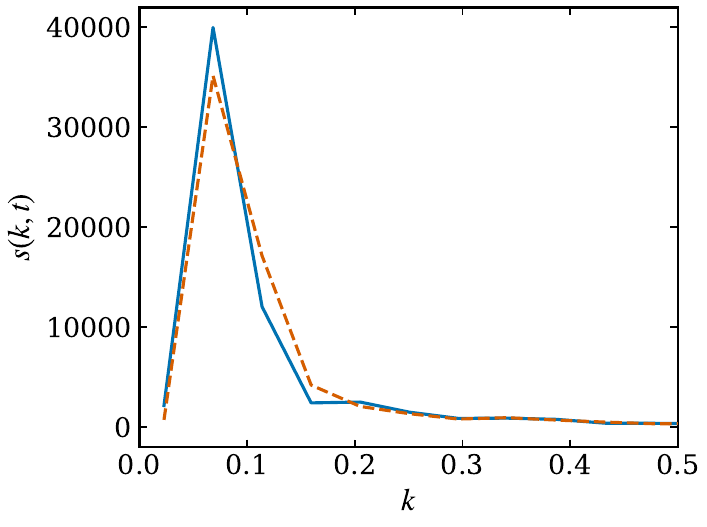}\\[0.2em]
    \includegraphics[width=\linewidth]{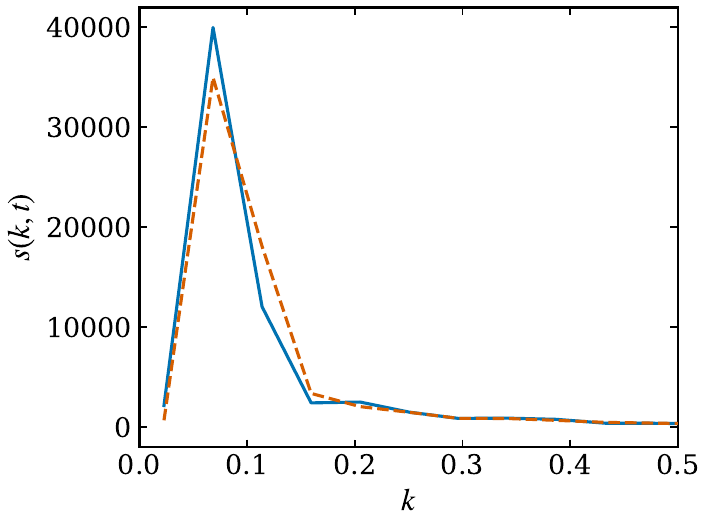}
  \end{minipage}
   \caption{Late-time normalized structure functions $s(k,t)$ at $t=29000$ for bcFNO
(top) and spFNO (bottom) at $N=256$.
The columns correspond to $\delta t=1$, $0.5$, and $0.2$ from left to right.
Blue and orange denote the reference trajectories and model predictions, respectively.}
\label{fig:app_sk_late_N256}
\end{figure}

Figs.~\ref{fig:app_sk_late_N128}
and~\ref{fig:app_sk_late_N256} compare the normalized structure
functions at $t=29000$.
For these cases, the spFNO normalized structure functions are closer to the reference curves than those of bcFNO, particularly near the
low-wavenumber peak that sets the domain scale.
This is the quantity directly targeted by the additional spFNO regularization.

Table~\ref{tab:app_coarsening_summary} reports the fitted slope $a$, the wavenumber exponent $\beta$, and the collapse score for the additional cases. 
For each case, $a$ and $\beta$ are fitted over the reference-defined coarsening interval, which is used unchanged for the numerical reference, bcFNO, and spFNO trajectories. 
The collapse score is computed from the coarsening-interval spectra. 
Its reference value is one by construction and is omitted from the table; larger values indicate closer agreement with the reference master curve.

\begin{table*}[!t]
\centering
\caption{Finite-regime coarsening results for the six additional cases.}
\label{tab:app_coarsening_summary}
\small
\begin{tabular}{c ccc ccc cc}
\toprule
& \multicolumn{3}{c}{$a$ in $L^3(t)=at+b$}
& \multicolumn{3}{c}{$\beta$ in $k_1(t)\sim t^\beta$}
& \multicolumn{2}{c}{$\mathcal R_{\mathrm{collapse}}$} \\
\cmidrule(lr){2-4}\cmidrule(lr){5-7}\cmidrule(lr){8-9}
$(N,\delta t)$
& $a_{\rm ref}$ & $a_{\rm bc}$ & $a_{\rm sp}$
& $\beta_{\rm ref}$  & $\beta_{\rm bc}$  & $\beta_{\rm sp}$
& $\mathcal R_{\rm bc}$ & $\mathcal R_{\rm sp}$ \\
\midrule
$(128,1.0)$
& $1.1\times10^{-2}$ & $1.5\times10^{-2}$ & $7.1\times10^{-3}$
& $-0.24$ & $-0.22$ & $-0.17$
& 0.35 & 0.40 \\
$(128,0.5)$
& $1.2\times10^{-2}$ & $2.3\times10^{-3}$ & $8.9\times10^{-3}$
& $-0.21$ & $-0.05$ & $-0.14$
& 0.55 & 0.93 \\
$(128,0.1)$
& $2.4\times10^{-2}$ & $1.5\times10^{-2}$ & $1.9\times10^{-2}$
& $-0.38$ & $-0.19$ & $-0.27$
& 0.86 & 0.74 \\
\addlinespace
$(256,1.0)$
& $1.8\times10^{-2}$ & $1.2\times10^{-2}$ & $7.4\times10^{-3}$
& $-0.28$ & $-0.14$ & $-0.09$
& 0.33 & 0.33 \\
$(256,0.5)$
& $1.4\times10^{-2}$ & $6.3\times10^{-3}$ & $8.4\times10^{-3}$
& $-0.22$ & $-0.09$ & $-0.10$
& 0.86 & 0.93 \\
$(256,0.2)$
& $1.6\times10^{-2}$ & $7.7\times10^{-3}$ & $7.5\times10^{-3}$
& $-0.24$ & $-0.11$ & $-0.10$
& 0.79 & 0.78 \\
\bottomrule
\end{tabular}
\end{table*}

For $(N,\delta t)=(128,0.5)$, spFNO brings both $a$ and $\beta$
closer to their reference values and increases the collapse score.
For $(128,0.1)$, both fitted quantities move closer to the reference, whereas the collapse score decreases. No simultaneous improvement of these three quantities is observed for $(128,1.0)$, $(256,1.0)$, or $(256,0.2)$.

\begin{figure}[!t]
  \centering
  \begin{minipage}{0.32\textwidth}
    \centering
    \includegraphics[width=\linewidth]{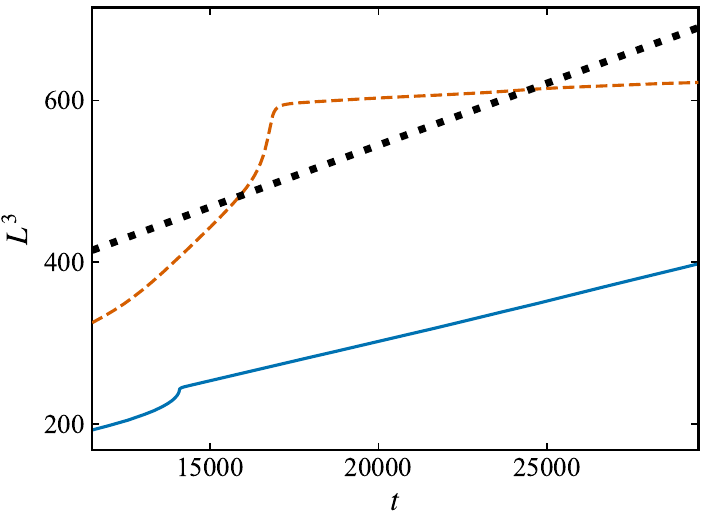}\\[0.25em]
    \includegraphics[width=\linewidth]{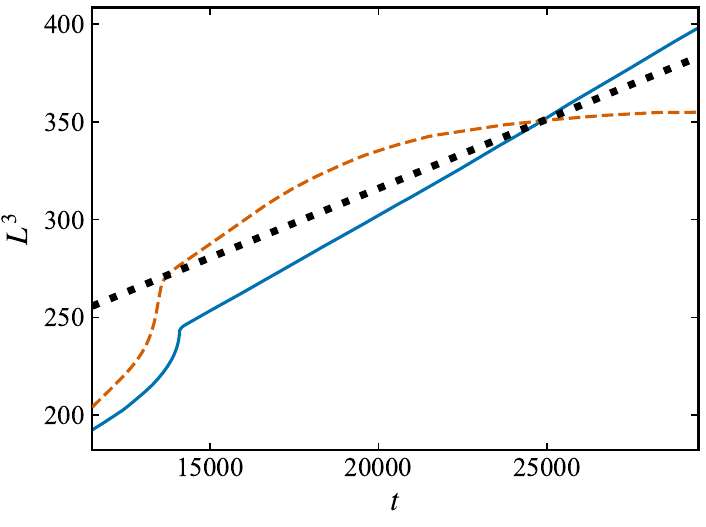}
  \end{minipage}
  \hfill
  \begin{minipage}{0.32\textwidth}
    \centering
    \includegraphics[width=\linewidth]{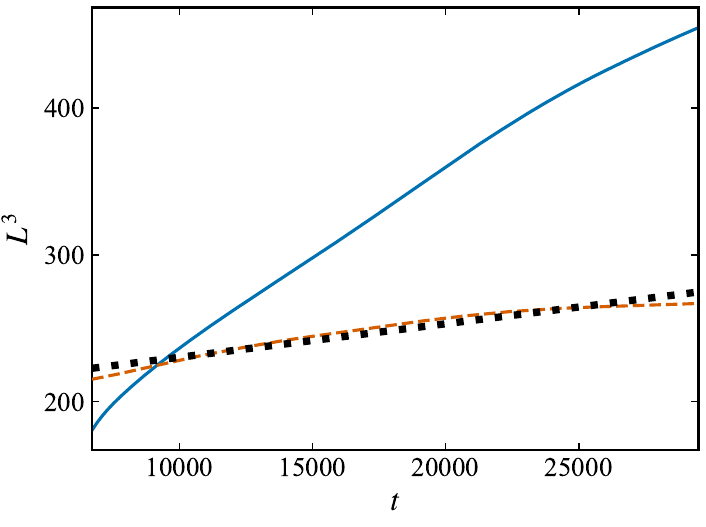}\\[0.25em]
    \includegraphics[width=\linewidth]{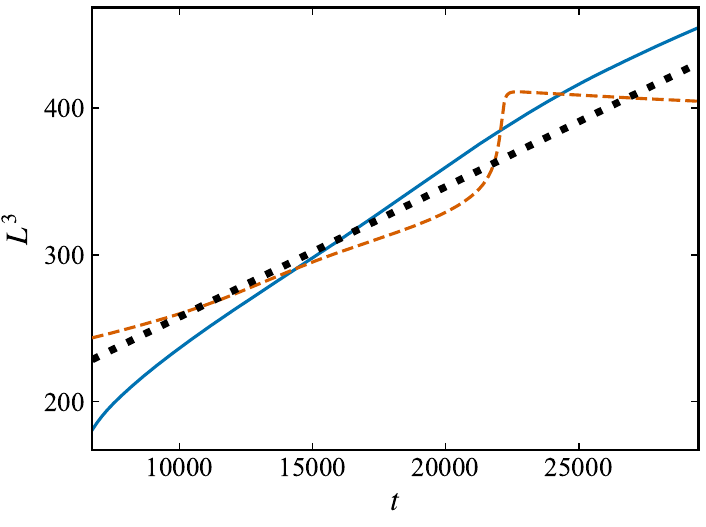}
  \end{minipage}
  \hfill
  \begin{minipage}{0.32\textwidth}
    \centering
    \includegraphics[width=\linewidth]{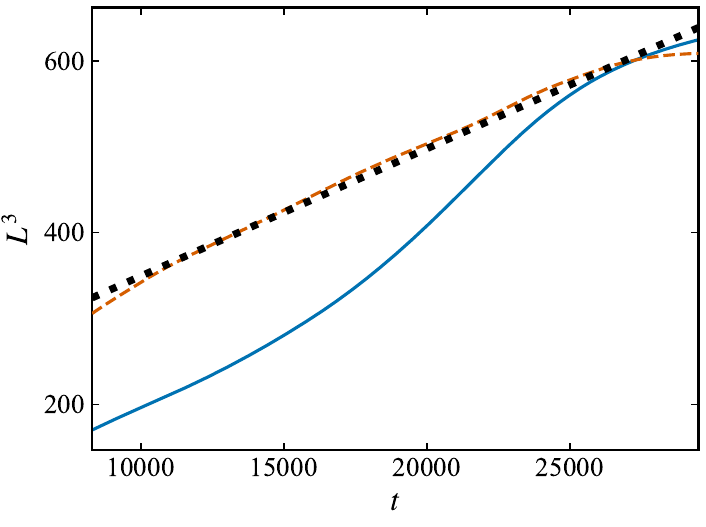}\\[0.25em]
    \includegraphics[width=\linewidth]{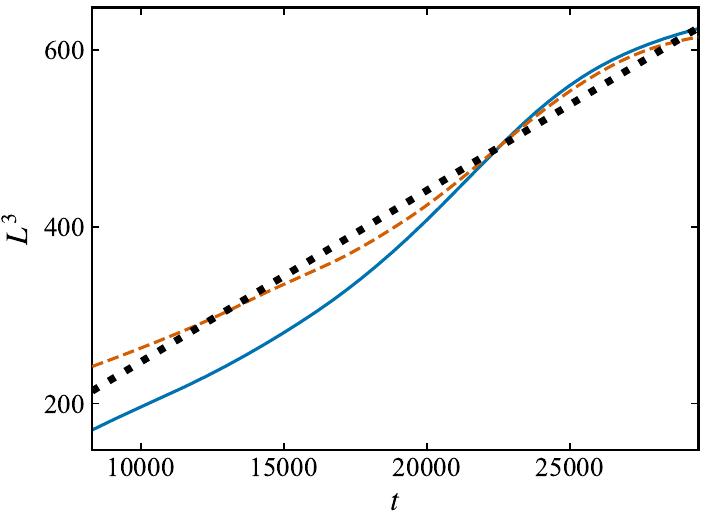}
  \end{minipage}
  \caption{Evolution of $L^3(t)$ for bcFNO (top) and spFNO (bottom) at
$N=128$.
The columns correspond to $\delta t=1$, $0.5$, and $0.1$ from left to right.
Blue and orange denote the reference trajectories and model predictions,
respectively.}
  \label{fig:app_L3_robustness_N128}
\end{figure}

\begin{figure}[!t]
  \centering
  \begin{minipage}{0.32\textwidth}
    \centering
    \includegraphics[width=\linewidth]{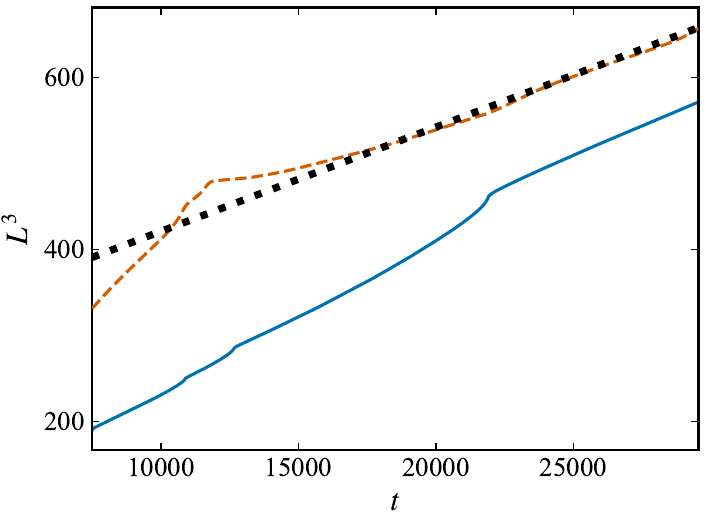}\\[0.25em]
    \includegraphics[width=\linewidth]{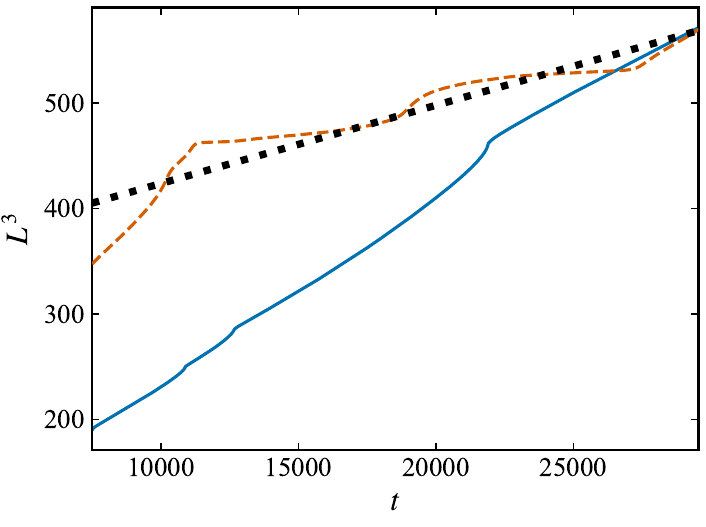}
  \end{minipage}
  \hfill
  \begin{minipage}{0.32\textwidth}
    \centering
    \includegraphics[width=\linewidth]{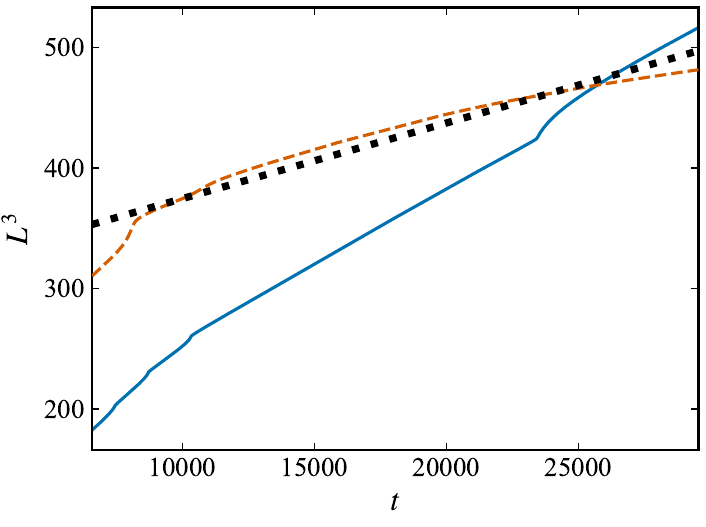}\\[0.25em]
    \includegraphics[width=\linewidth]{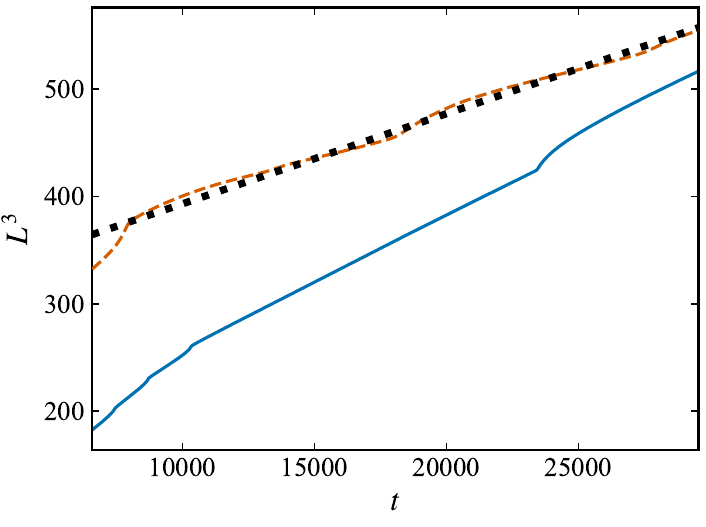}
  \end{minipage}
  \hfill
  \begin{minipage}{0.32\textwidth}
    \centering
    \includegraphics[width=\linewidth]{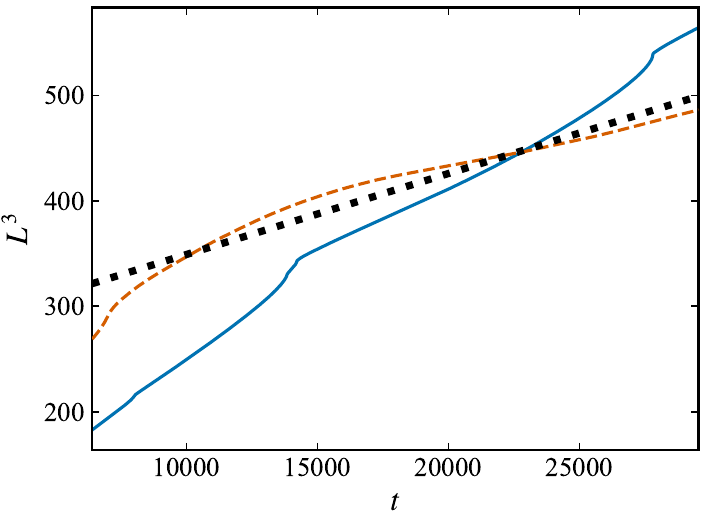}\\[0.25em]
    \includegraphics[width=\linewidth]{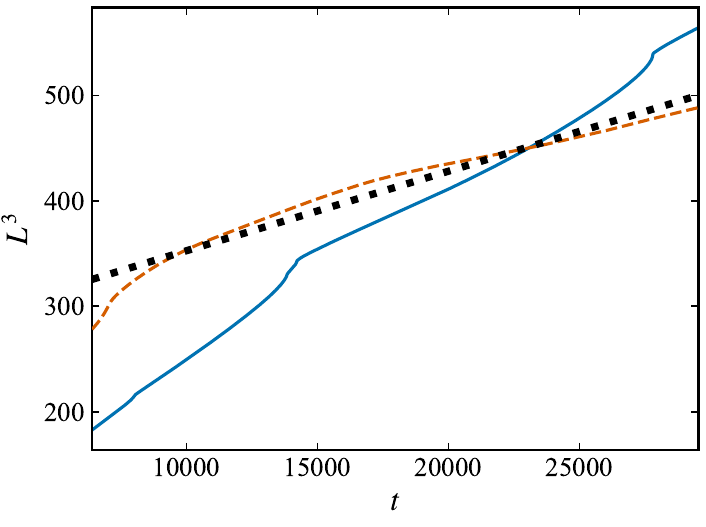}
  \end{minipage}
  \caption{Evolution of $L^3(t)$ for bcFNO (top) and spFNO (bottom) at $N=256$. The columns correspond to $\delta t=1$, $0.5$, and $0.2$ (from left to right). Blue and orange curves denote the reference trajectories and model predictions, respectively. }
  \label{fig:app_L3_robustness_N256}
\end{figure}

\begin{figure}[!t]
  \centering
  \begin{minipage}{0.32\textwidth}
    \centering
    \includegraphics[width=\linewidth]{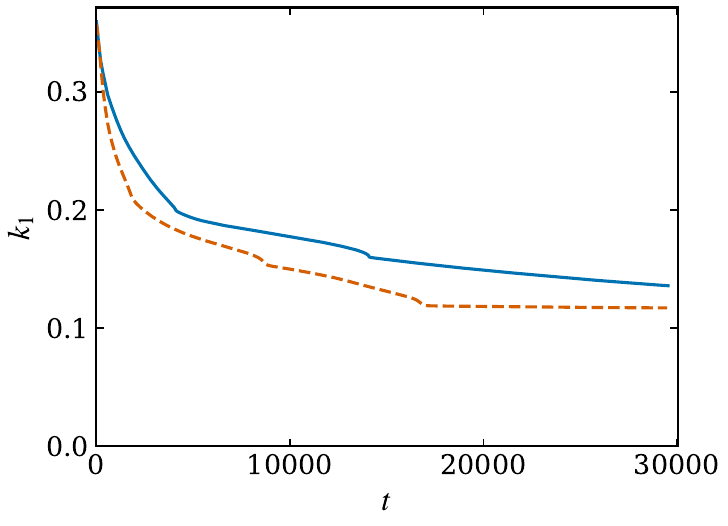}\\[0.25em]
    \includegraphics[width=\linewidth]{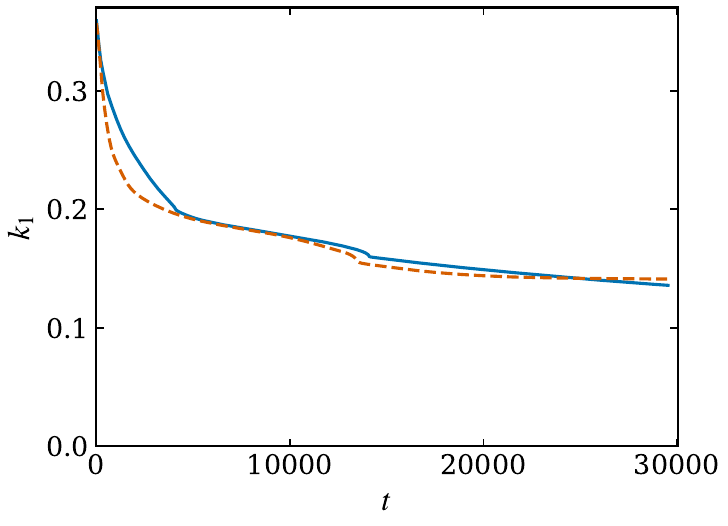}
  \end{minipage}
  \hfill
  \begin{minipage}{0.32\textwidth}
    \centering
    \includegraphics[width=\linewidth]{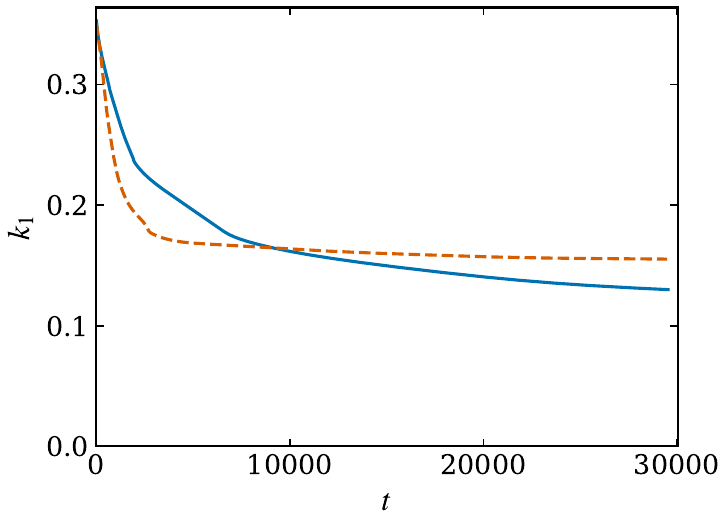}\\[0.25em]
    \includegraphics[width=\linewidth]{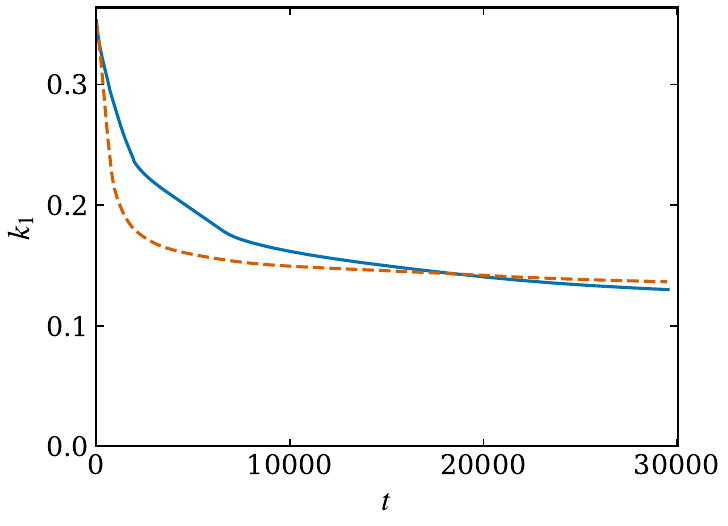}
  \end{minipage}
  \hfill
  \begin{minipage}{0.32\textwidth}
    \centering
    \includegraphics[width=\linewidth]{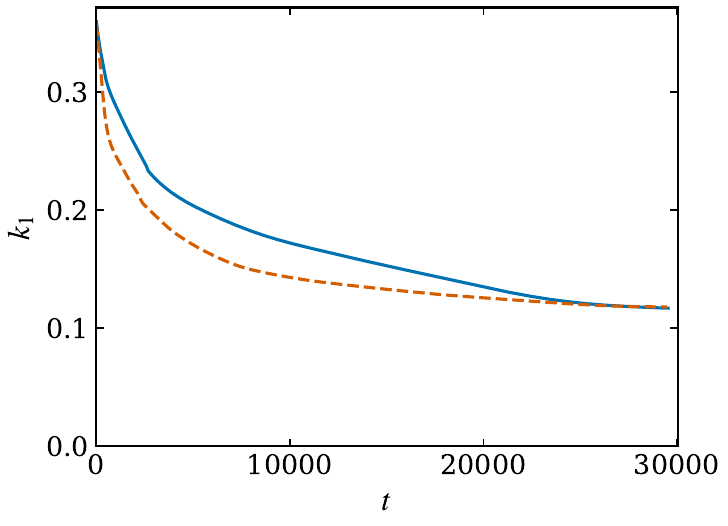}\\[0.25em]
    \includegraphics[width=\linewidth]{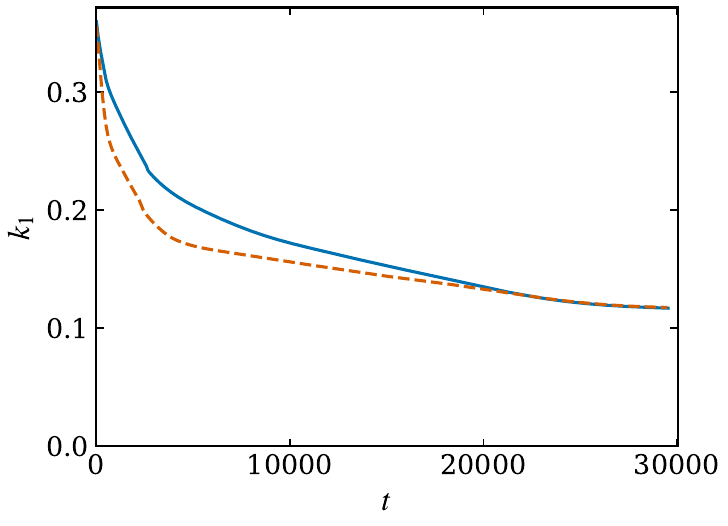}
  \end{minipage}
  \caption{Evolution of $k_1(t)$ for bcFNO (top) and spFNO (bottom) at $N=128$. The columns correspond to $\delta t=1$, $0.5$, and $0.1$ (from left to right). Blue and orange curves denote the reference trajectories and model predictions, respectively.}
  \label{fig:app_k1_robustness_N128}
\end{figure}

\begin{figure}[!t]
  \centering
  \begin{minipage}{0.32\textwidth}
    \centering
    \includegraphics[width=\linewidth]{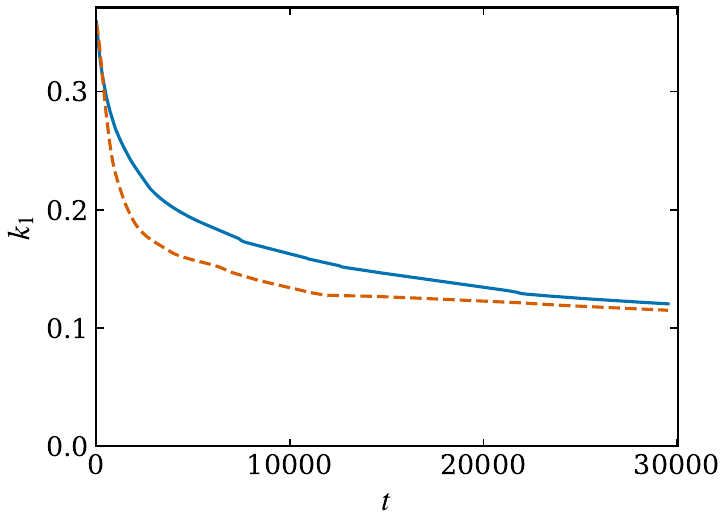}\\[0.25em]
    \includegraphics[width=\linewidth]{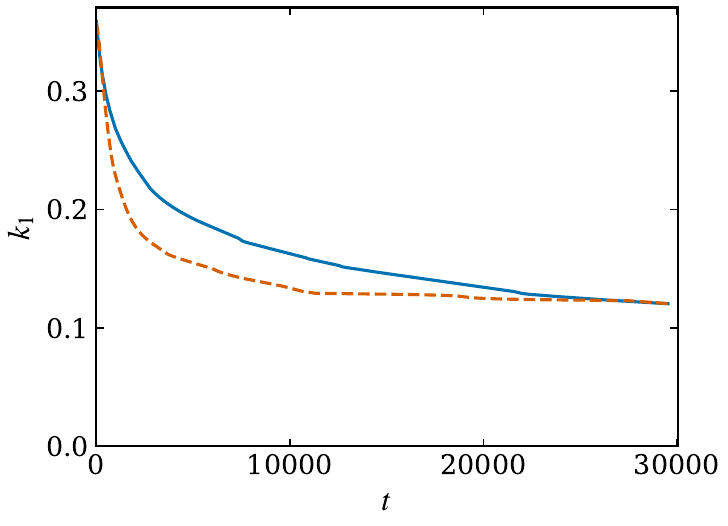}
  \end{minipage}
  \hfill
  \begin{minipage}{0.32\textwidth}
    \centering
    \includegraphics[width=\linewidth]{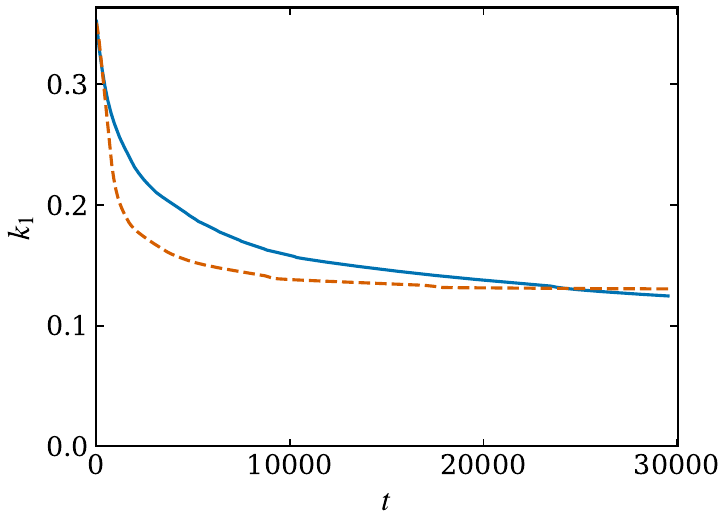}\\[0.25em]
    \includegraphics[width=\linewidth]{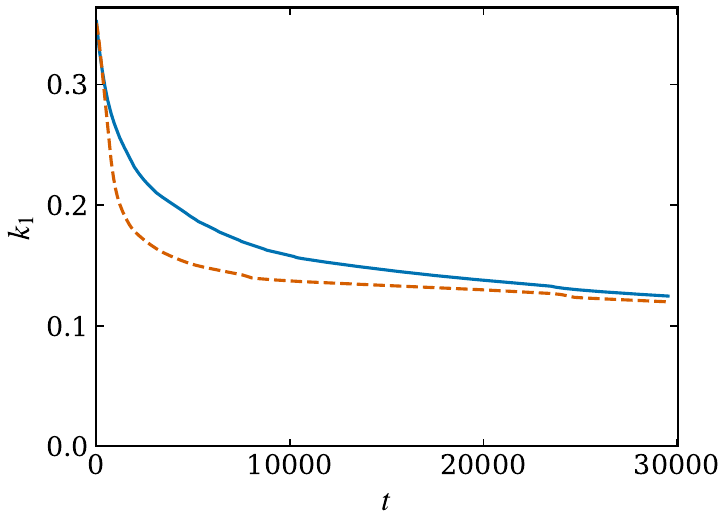}
  \end{minipage}
  \hfill
  \begin{minipage}{0.32\textwidth}
    \centering
    \includegraphics[width=\linewidth]{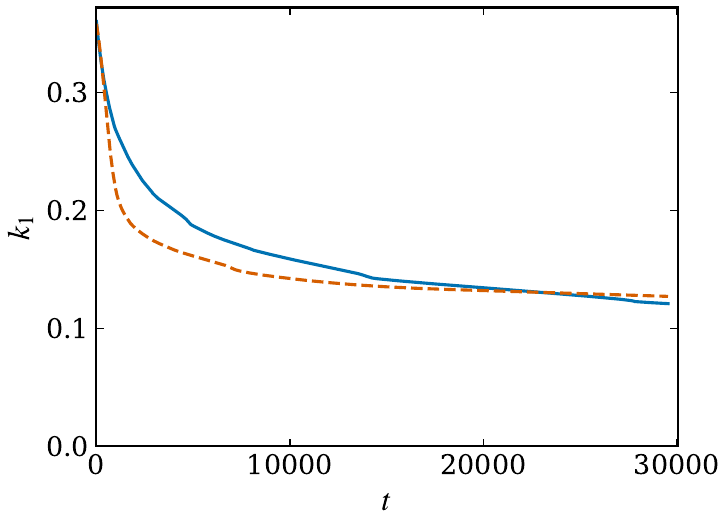}\\[0.25em]
    \includegraphics[width=\linewidth]{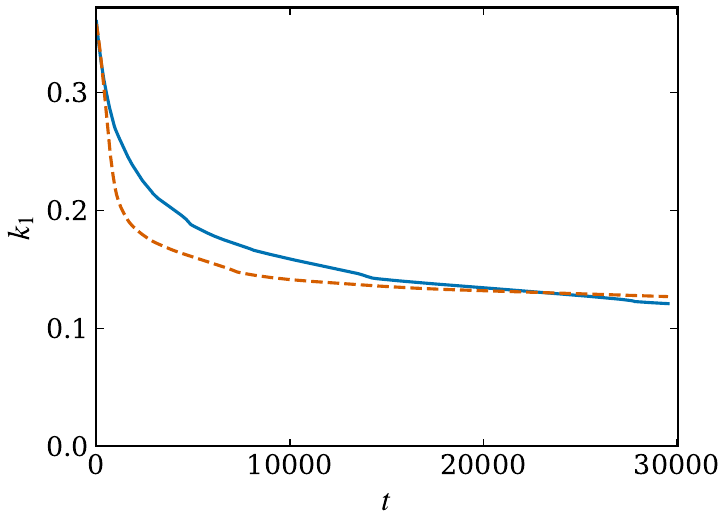}
  \end{minipage}
  \caption{Evolution of $k_1(t)$ for bcFNO (top) and spFNO (bottom) at $N=256$. The columns correspond to $\delta t=1$, $0.5$, and $0.2$ (from left to right). Blue and orange curves denote the reference trajectories and model predictions, respectively. }
  \label{fig:app_k1_robustness_N256}
\end{figure}

\begin{figure}[!t]
  \centering
  \begin{minipage}{0.32\textwidth}
    \centering
    \includegraphics[width=\linewidth]{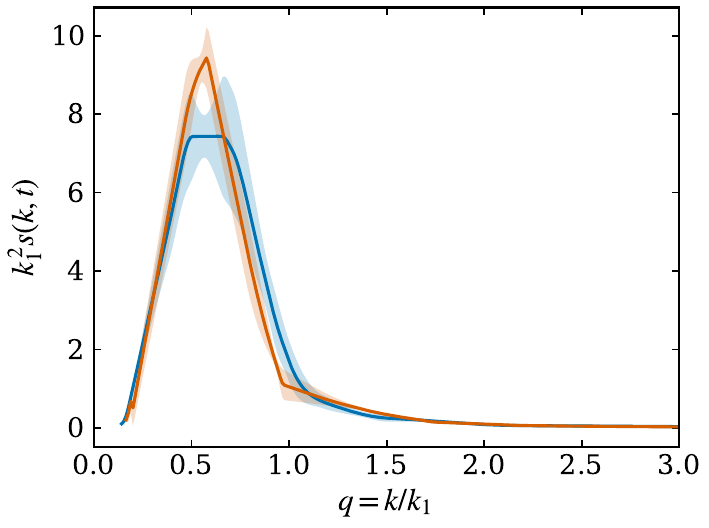}\\[0.25em]
    \includegraphics[width=\linewidth]{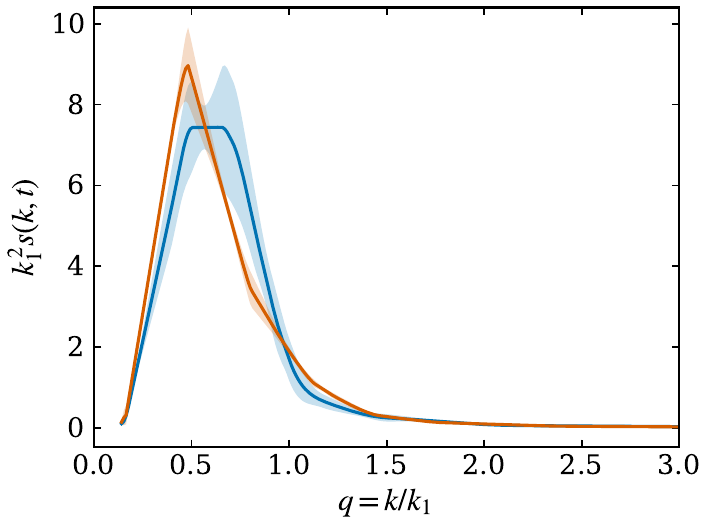}
  \end{minipage}
  \hfill
  \begin{minipage}{0.32\textwidth}
    \centering
    \includegraphics[width=\linewidth]{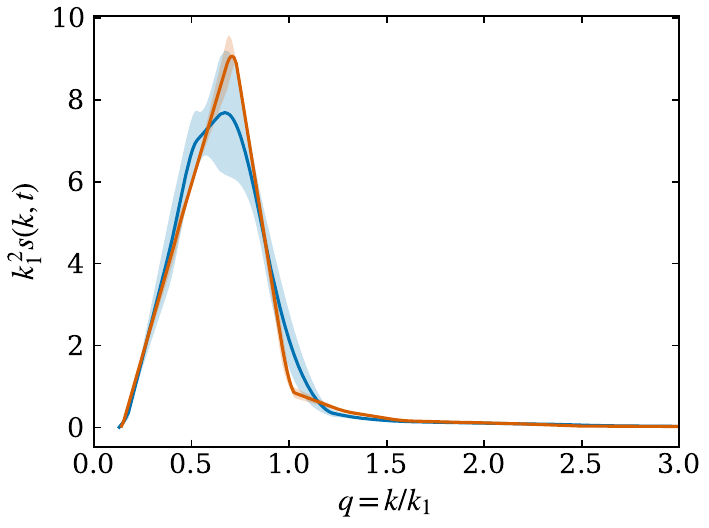}\\[0.25em]
    \includegraphics[width=\linewidth]{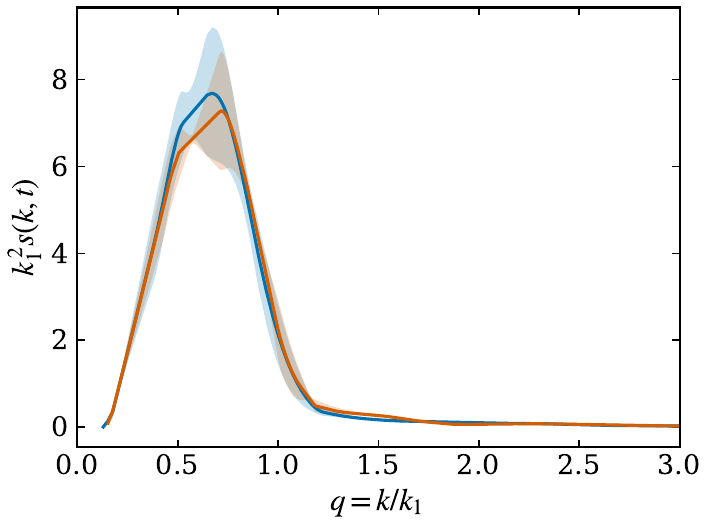}
  \end{minipage}
  \hfill
  \begin{minipage}{0.32\textwidth}
    \centering
    \includegraphics[width=\linewidth]{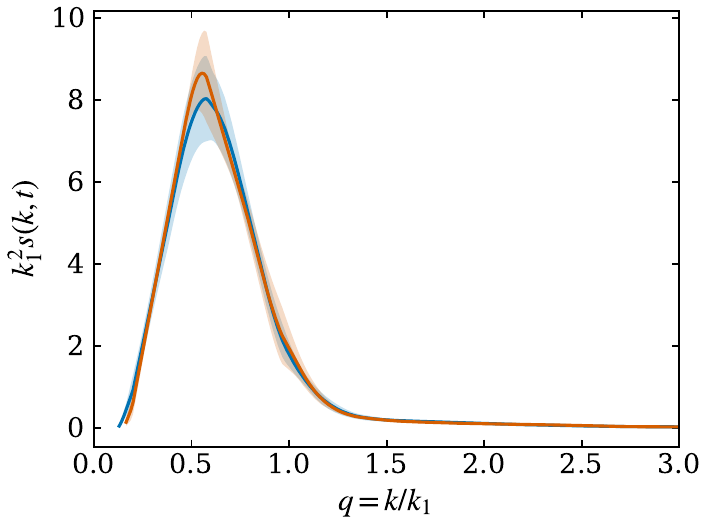}\\[0.25em]
    \includegraphics[width=\linewidth]{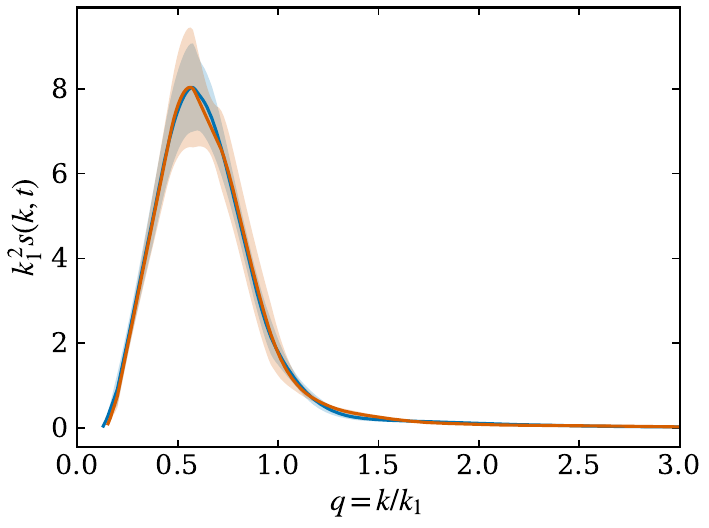}
  \end{minipage}
  \caption{Dynamic-scaling collapse for bcFNO (top) and spFNO (bottom) at $N=128$. The columns correspond to $\delta t=1$, $0.5$, and $0.1$ (from left to right). Blue and orange curves denote the reference trajectories and model predictions, respectively. The mean curves and shaded envelopes are constructed as described in
Sec.~\ref{subsec:exp_sk}.}
  \label{fig:app_robustness_N128}
\end{figure}

\begin{figure}[!t]
  \centering
  \begin{minipage}{0.32\textwidth}
    \centering
    \includegraphics[width=\linewidth]{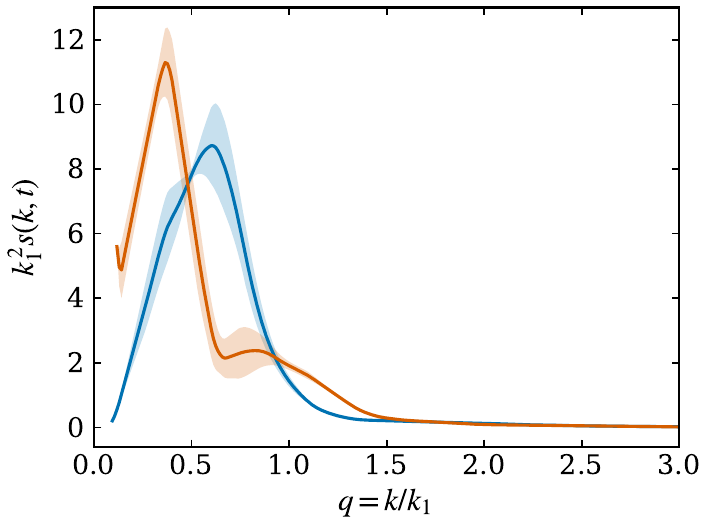}\\[0.25em]
    \includegraphics[width=\linewidth]{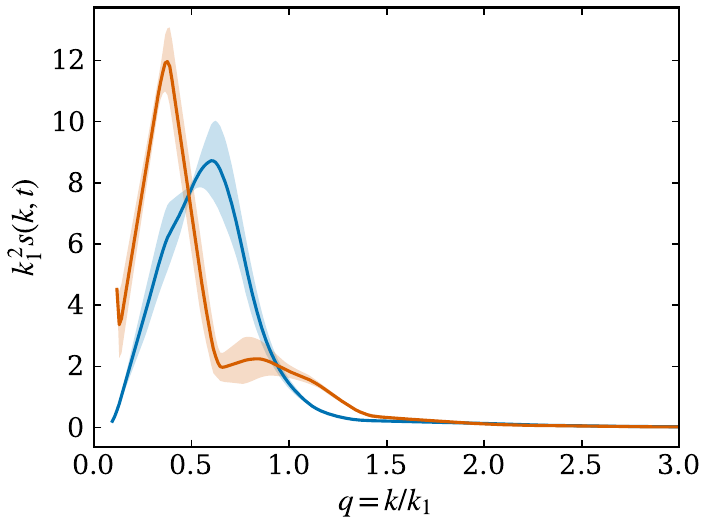}
  \end{minipage}
  \hfill
  \begin{minipage}{0.32\textwidth}
    \centering
    \includegraphics[width=\linewidth]{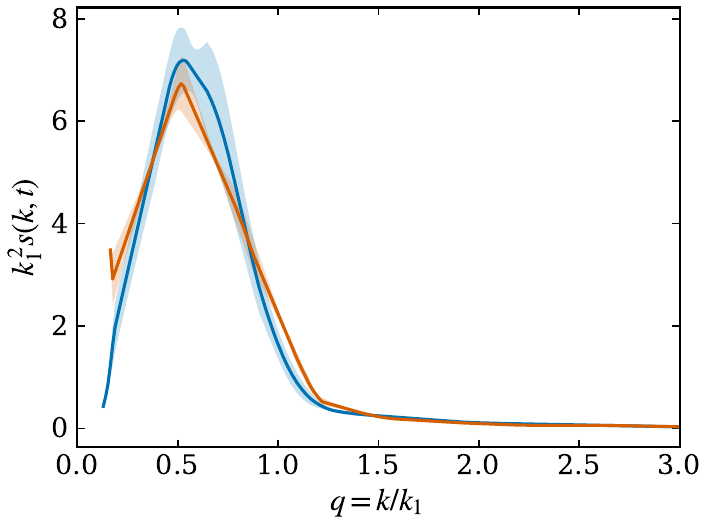}\\[0.25em]
    \includegraphics[width=\linewidth]{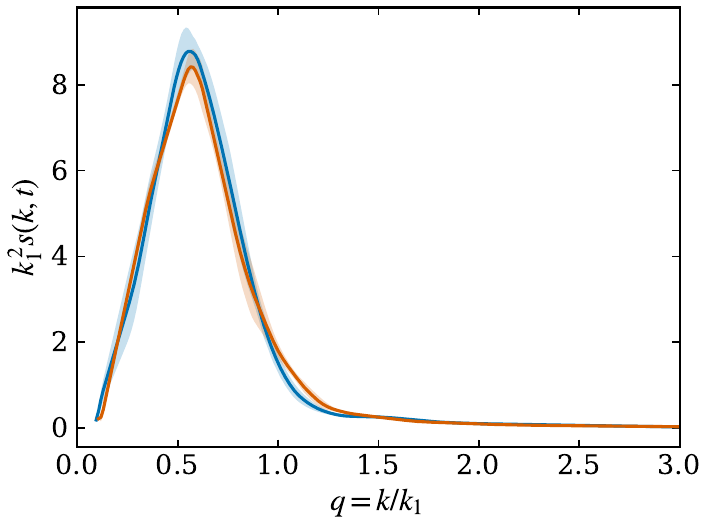}
  \end{minipage}
  \hfill
  \begin{minipage}{0.32\textwidth}
    \centering
    \includegraphics[width=\linewidth]{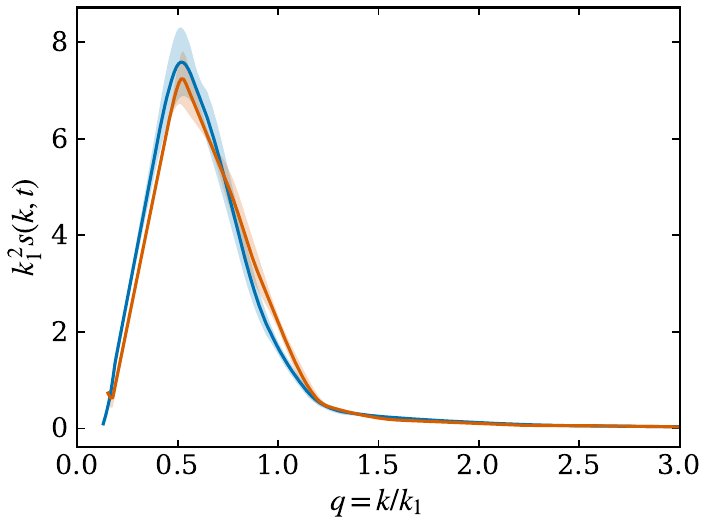}\\[0.25em]
    \includegraphics[width=\linewidth]{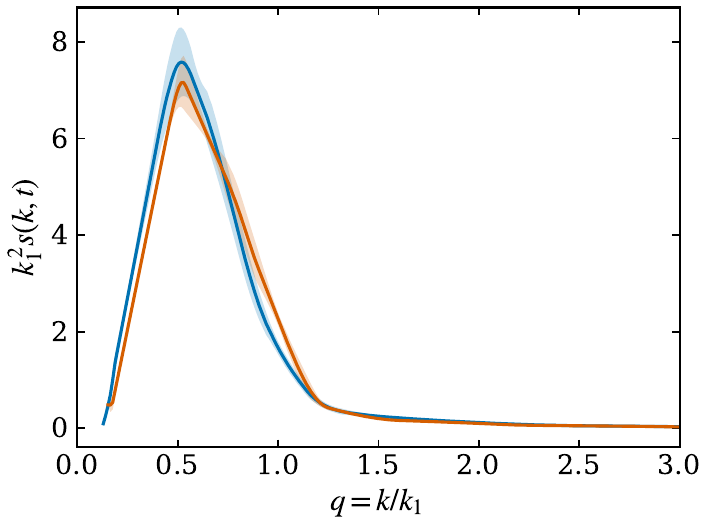}
  \end{minipage}
  \caption{Dynamic-scaling collapse for bcFNO (top) and spFNO (bottom) at $N=256$. The columns correspond to $\delta t=1$, $0.5$, and $0.2$ (from left to right). Graphical conventions follow Fig.~\ref{fig:app_robustness_N128}.}
  \label{fig:app_robustness_N256}
\end{figure}

The spFNO compares $s(k,t)$ at individual prediction times.
By contrast, $a$ is obtained from an affine fit of $L^3(t)$ over the coarsening interval, $\beta$ from a power-law fit of $k_1(t)$, and $\mathcal R_{\mathrm{collapse}}$ from the rescaled spectra. 
Thus, across the tested cases, reductions in the late-stage discrepancy of $s(k,t)$ do not translate uniformly into fitted values closer to the reference for all three scaling quantities.
Figs.~\ref{fig:app_L3_robustness_N128}--\ref{fig:app_k1_robustness_N256} show the trajectories used for these fits.
Figures~\ref{fig:app_robustness_N128}
and~\ref{fig:app_robustness_N256} show the rescaled structure functions used to compute the collapse scores. The binning and averaging follow Sec.~\ref{subsec:exp_sk}.

\end{document}